\documentclass{article}
\usepackage[utf8]{inputenc}
\usepackage[margin=1in]{geometry}
\usepackage{amsmath}
\usepackage[dvipsnames]{xcolor}
\usepackage[section]{placeins} % keep floats within section they are written
\usepackage{lineno}
\usepackage{natbib}
\usepackage[title]{appendix}
\usepackage{graphicx}
\usepackage{parskip} 
\usepackage{xcite}
\usepackage{xr}
\usepackage{caption}
\usepackage{subcaption}
\makeatletter
\newcommand*{\addFileDependency}[1]{
  \typeout{(#1)}
  \@addtofilelist{#1}
  \IfFileExists{#1}{}{\typeout{No file #1.}}
}
\makeatother

\providecommand{\keywords}[1]
{
  \small	
  \textbf{\textbf{Keywords --}} #1
}
\usepackage{titling}
\predate{}
\postdate{}
\date{}

\title{Structural sensitivity in the functional responses of predator-prey models}
\author{Sarah K. Wyse$^{\text{a},\star}$, Maria M. Martignoni$^{\text{b}}$, May Anne Mata$^{\text{c}}$,\\ Eric Foxall$^{\text{a}}$, and Rebecca C. Tyson$^{\text{a}}$}

\begin{document}

\maketitle
%\linenumbers

$^{\text{a}}$ Department of Computer Science, Mathematics, Physics and Statistics, Irving K. Barber Faculty of Science, University of British Columbia Okanagan, SCI354, 1177 Research Road, Kelowna V1V 1V7, BC, Canada\\
$^{\text{b}}$ Department of Mathematics and Statistics, Memorial University of Newfoundland, 230 Elizabeth Ave, St. John's, A1C 5S7, NL, Canada\\
$^{\text{c}}$ Department of Math, Physics, and Computer Science, University of the Philippines Mindanao, Mintal, Davao City, Philippines\\
$^{\star}$ Corresponding author. \\
\textit{Email addresses}: sarah.wyse@alumni.ubc.ca (S.K. Wyse), mmartignonim@mun.ca (M.M. Martignoni),\\ memata@up.edu.ph (M.A. Mata), efoxall@mail.ubc.ca (E. Foxall), rebecca.tyson@ubc.ca (R.C. Tyson).\\
\textit{ORCID IDs}: 0000-0002-5500-2711 (S.K. Wyse), 0000-0002-0192-8228 (M.M. Martignoni), 0000-0002-2967-344X (M.A. Mata), 0000-0003-1610-6662 (E. Foxall), 0000-0002-7373-4473 (R.C. Tyson).

\begin{abstract}
    In mathematical modeling, several different functional forms can often be used to fit a data set equally well, especially if the data is sparse.  In such cases, these mathematically different but similar looking functional forms are typically considered interchangeable. Recent work, however, shows that similar functional responses may nonetheless result in significantly different bifurcation points for the Rosenzweig-MacArthur predator-prey system. Since the bifurcation behaviours include destabilising oscillations, predicting the occurrence of such behaviours is clearly important. Ecologically, different bifurcation behaviours mean that different predictions may be obtained from the models. These predictions can range from stable coexistence to the extinction of both species, so obtaining more accurate predictions is also clearly important for conservationists. Mathematically, this difference in bifurcation structure given similar functional responses is called structural sensitivity. We extend the existing work to find that the Leslie-Gower-May predator-prey system is also structurally sensitive to the functional response. Using the Rosenzweig-MacArthur and Leslie-Gower-May models, we then aim to determine if there is some way to obtain a functional description of data so that different functional responses yield the same bifurcation structure, i.e., we aim to describe data such that our model is not structurally sensitive. We first add stochasticity to the functional responses and find that better similarity of the resulting bifurcation structures is achieved. Then, we analyze the functional responses using two different methods to determine which part of each function contributes most to the observed bifurcation behaviour. We find that prey densities around the coexistence steady state are most important in defining the functional response. Lastly, we propose a procedure for ecologists and mathematical modelers to increase the accuracy of model predictions in predator-prey systems.
\end{abstract}

\keywords{Structural sensitivity, bifurcation structure, functional response, coexistence steady state, stochasticity, prey density domain}

\textbf{Funding:} This work was supported by The University of British Columbia Okanagan [Undergraduate Research Award] and Natural Sciences and Engineering Research Council of Canada [grant numbers RGPIN-2016-0577, RGPIN-2018-04480]. The funding sources were not involved in the conduction of research or preparation of this article.

\newpage

\section{Introduction} \label{Sect:Introduction}
%--------------------------------------------------------------------------------
Mathematical modeling of biological phenomena necessitates making assumptions about the biological system so that the underlying mechanisms can be expressed mathematically. Indeed, the modeler must decide which features of a system to include and how to realistically formulate them \citep{krkosek:2007, barclay:2001, molnar:2014, tyson:1995, tyson:2016, tyson:2020}. This process becomes more complicated when there are many apparently equivalent mathematical functions that may be used \citep{fussmann:2005}. In this paper, we are interested in systems where the data is typically sparse and/or noisy, making it impossible to select one best-fitting functional curve \citep{krebs:2001}. We focus here on the case of predator-prey systems and, in particular, on the predation interaction. The mathematical expression of this interaction is often referred to as the functional response.

A common assumption in modeling is that, if mathematical functions possess similar features and are appropriately parametrized to resemble one another, then the models utilizing these functional responses will produce essentially the same behaviour. Hence, many mathematicians typically choose a convenient functional response that fits the predation data. It was found, however, that similar looking functional responses may produce very different results \citep{fussmann:2005}. Indeed, \citet{fussmann:2005} found that, for a given prey density, their model can produce high amplitude oscillations, bistability, or stable coexistence at low densities, depending on the form of the functional response, even if these responses are visually indistinguishable. This difference in model behaviour means that the assumption that similar curves can be used interchangeably does not always hold. This behaviour is called structural sensitivity. A structurally sensitive model has the property that a small perturbation in model functions, here the functional response, can lead to a different bifurcation structure. Ecologically, predictions made using models with different bifurcation structures could range from stable coexistence to extinction of both species.

The predator-prey system in \citet{fussmann:2005} is a Rosenzweig-MacArthur model that makes use of three similar functional response curves that are visually indistinguishable with respect to data fitting. The three functional responses are: Holling Type II, which is derived by considering how the predator's activity is split between searching for and handling prey; Ivlev, which is based on the maximum digestion rate of the predator; and a hyperbolic tangent (trigonometric), which is a phenomenological curve that has no biological basis but is used in some population models \citep{aldebert:2018}. Although the theory behind each of these functional responses is different, they are considered interchangeable because it is difficult to determine if a predator is limited by its handling time or its digestion time. \citet{jeschke:2002} have developed a model that includes both of these limitations on predation together but we do not consider it in this manuscript. These functional responses have five specific characteristics in common: zero predation at zero prey density, monotonically increasing predation with increasing prey density, a negative second derivative, saturating behaviour at high prey density, and smoothness \citep{fussmann:2005, seoWolkowicz:2018}.

\citet{fussmann:2005} found that the Rosenzweig-MacArthur model with these three functional response curves is structurally sensitive. That is, using the carrying capacity as the bifurcation parameter, this model produces different bifurcation diagrams for all three very similar looking functional responses. The model has a Hopf bifurcation of the coexistence steady state for all three functional response curves, but the critical value of the carrying capacity at which the Hopf bifurcation occurs varies by as much as 2000\% from the smallest to largest critical transition value. For large values of the carrying capacity, the model that uses the Holling functional response produces large oscillations with limit cycles that approach very close to the prey-only and extinction steady states. On the other hand, the model using the trigonometric functional response initially produces smaller oscillations with limit cycles further from the prey-only and extinction steady states. The model using the Ivlev functional response produces oscillations of intermediate size. Additionally, the trigonometric model is the only model that has bistability. The difference in these results is problematic because the model disagreements make it impossible to determine, for example, if a real system is nearing a critical transition, such as the rapid emergence of large oscillations or extinction \citep{scheffer:2009}. This dilemma leads us to ask if there is a better way to mathematically describe data such that there is better similarity of the bifurcation structures obtained from different functional responses with respect to the bifurcation points, the coexistence steady state, and the branches of the limit cycles.

In recent literature, \citet{adamson:2014, adamson:2013} have shown the existence of structural sensitivity in a range of mathematical models and have developed a test to determine whether or not a model is structurally sensitive. This test, however, does not address the question of how to alter the functional responses in a model to reduce or remove structural sensitivity, it only suggests that care must be taken in interpreting model predictions.

In this paper, we investigate three different approaches for translating a description of the data into a functional response curve. Our goal is to determine if any of our approaches can decrease or even remove structural sensitivity. First, we repeat the work of \citet{fussmann:2005} on another predator-prey model, the Leslie-Gower-May model, and find that it is also structurally sensitive to its functional response. We then base our study of structural sensitivity on the behaviour of these two different predator-prey models. Next, we study the effect of stochasticity of various types and amplitudes on the system and observe the resulting changes in structural sensitivity. Lastly, we investigate what part(s) of the functional response might be related to structural sensitivity. 

In Section \ref{Model}, we introduce the models and functional responses. In Section \ref{Methods}, we explain the methods we use to investigate structural sensitivity in the functional response. In Section \ref{Results}, we state our results. In Section \ref{Discussion}, we discuss and relate our findings to previous work and suggest a method to improve predictions obtained from structurally sensitive models.

\section{Models}
%-------------------------------------------------------------------------------------
\label{Model}
In this paper, we consider two classic predator-prey models: the Rosenzweig-MacArthur model and the Leslie-Gower-May model. 

In the Rosenzweig-MacArthur model, prey density, $x$, grows via logistic growth with a growth rate $r$ and carrying capacity $K$. Predators and prey interact via a predation term, $f_i(x)$, which is discussed in more detail below. Lastly, predator density, $y$, increases through growth related to predation and decreases through natural death with mortality rate $m$. The Rosenzweig-MacArthur model is written \citep{turchin:2013}:
\begin{subequations}
\label{RM}
\begin{align}
	\frac{dx}{dt} &= rx\left(1-\frac{x}{K}\right)  - f_i(x)y, \label{RM1}\\
    \frac{dy}{dt} &= f_i(x)y - my. \label{RM2}
\end{align}
\end{subequations}

The parameter values given in \citet{fussmann:2005} are used throughout this paper (for the Rosenzweig-MacArthur model) unless specified otherwise.

The prey density equation in the Leslie-Gower-May model is the same as that in the Rosenzweig-MacArthur model. The predator density equation, on the other hand, has logistic growth of the predator with a growth rate $s$ and carrying capacity $x/q$ that decreases as prey density decreases. Here, the parameter $q$ is the equilibrium density ratio \citep{turchin:1997}. The Leslie-Gower-May model is written \citep{turchin:2013}: 
\begin{subequations}
\label{May}
\begin{align}
	\frac{dx}{dt} &= rx\left(1-\frac{x}{K}\right) - f_i(x)y, \label{May1}\\
	\frac{dy}{dt} &= sy\left(1-\frac{qy}{x}\right). \label{May2}
\end{align}
\end{subequations}

The parameter values given in \citet{vitense:2016} are used throughout this paper (for the Leslie-Gower-May model) unless specified otherwise.

The predation functional response, $f_i(x)$, is our main focus. Following \citet{fussmann:2005}, we consider the following three functional responses: Holling, Ivlev, and trigonometric. These functional responses are described, respectively, by the following three functions:

\begin{subequations}
\label{functionalResponses}
\begin{align}
    f_H(x) &= \frac{a_Hx}{1+b_Hx} \label{Holling}, \\
    f_I(x) &= a_I(1-\text{exp}(b_Ix)) \label{Ivlev}, \\
    f_T(x) &= a_T\text{tanh}(b_Tx) \label{trigonometric},
\end{align}
\end{subequations}
where $a_i$ and $b_i$ are the two parameters through which the curves can be fit to data, or to each other. The Rosenzweig-MacArthur functional responses have units of prey/predator/day and the Leslie-Gower-May functional responses have units of prey/predator/year. Note that in application, the three functional responses fit data equally well. Further, since each functional response has two parameters, it is not possible to distinguish the models using criteria such as the Akaike Information Criteria or Schwarz/Bayesian Information Criteria as these criteria rely greatly on the number of parameters used in each functional response \citep{Akaike1974}.

\section{Methods}
%------------------------------------------------------------------------------------
\label{Methods}
Structural sensitivity has not previously been established in the Leslie-Gower-May model. In this paper, we first show that model \eqref{May} is structurally sensitive. Then we study the structural sensitivity of the Rosenzweig-MacArthur and Leslie-Gower-May models by altering the functional responses \eqref{functionalResponses} in the following three ways: 
\begin{enumerate}
    \item Adding various types of stochasticity (Section \ref{Sect:stochMethods}).
    \item Varying the parameters $a_i$ and $b_i$ to fit the functional responses over a specified subset of the prey density domain (Section \ref{Sect:CoexSSMethods}).
    \item Using continuous curves defined piecewise from two functional responses, switching from one to the other at a point where the functional response curves intersect (Section \ref{Sect:SensRegionsMethods}). 
\end{enumerate}

Method 1 uses stochasticity to create more overlap in the functional responses and hypothesizes that if the functional responses are more similar, then the resulting model behaviours will also become more similar. Methods 2 and 3 consist of focusing on different intervals of prey densities by applying a different weight to each part of the functional response. Our goal with these methods is to determine if there is a portion of the functional response that dictates the bifurcation structure.

Before investigating these methods to determine if we can decrease or remove structural sensitivity, we first fit our functional responses. To do this fitting, one of the functional responses is fixed and the other two are fit to the first response via nonlinear least squares to solve for the unknown $a_i$ and $b_i$ values. We use the Nelder-Mead algorithm \citep{Nelder:1965} to apply nonlinear least squares and chose the interval of prey densities to be [0,4] to match the previous work by \citet{fussmann:2005}.

\begin{figure}[ht]
    \centering
    \begin{subfigure}{0.45\linewidth}
        \includegraphics[width=\linewidth]{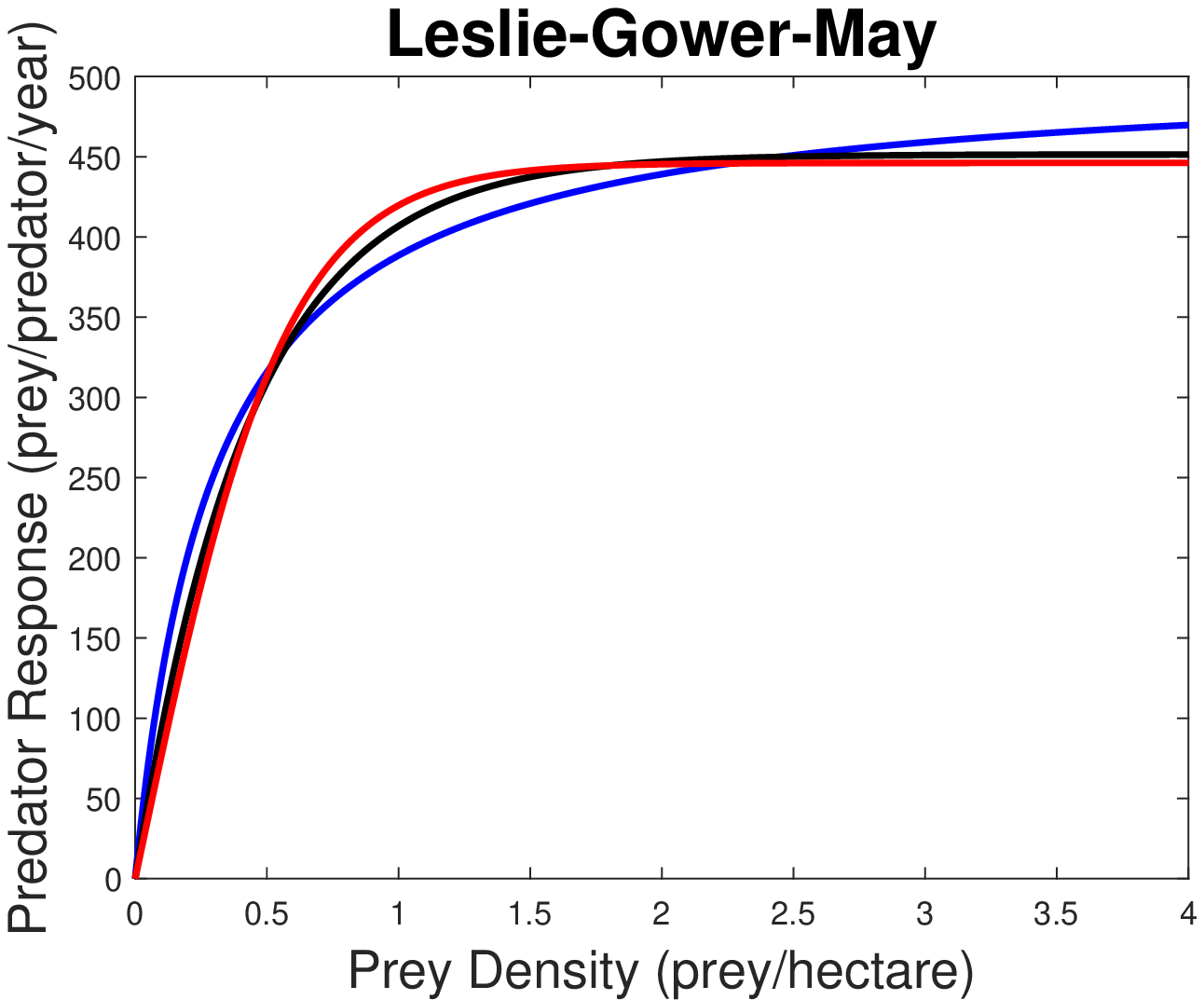}
        \subcaption{Original Functional Responses}
        \label{fig:CurvesOrigMay}
    \end{subfigure}
    \begin{subfigure}{0.45\linewidth}
        \includegraphics[width=\linewidth]{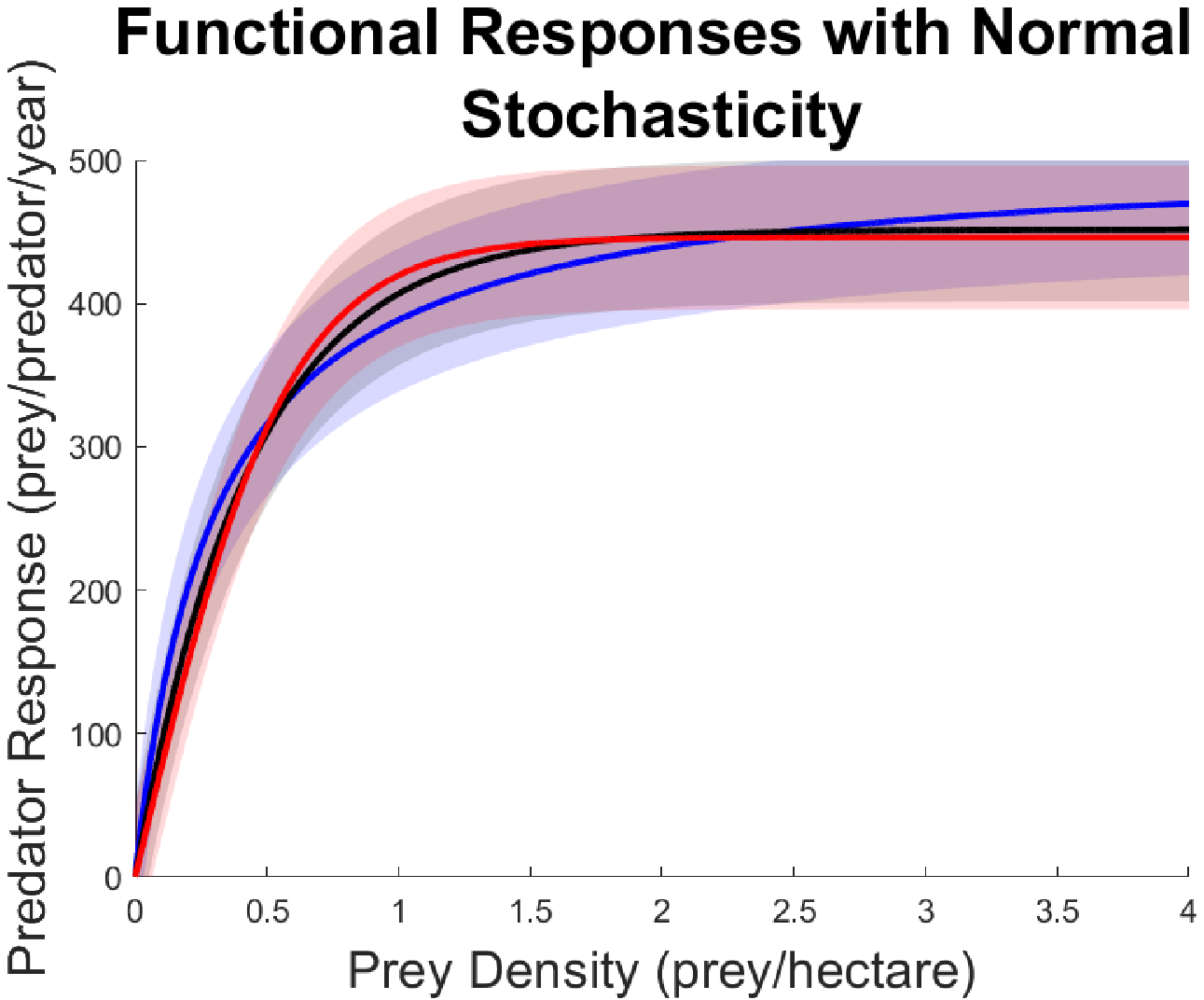}
        \subcaption{Stochastic Functional Responses}
        \label{fig:StochCurves}
    \end{subfigure}
    \\
    \begin{subfigure}{0.45\linewidth}
        \includegraphics[width=\linewidth]{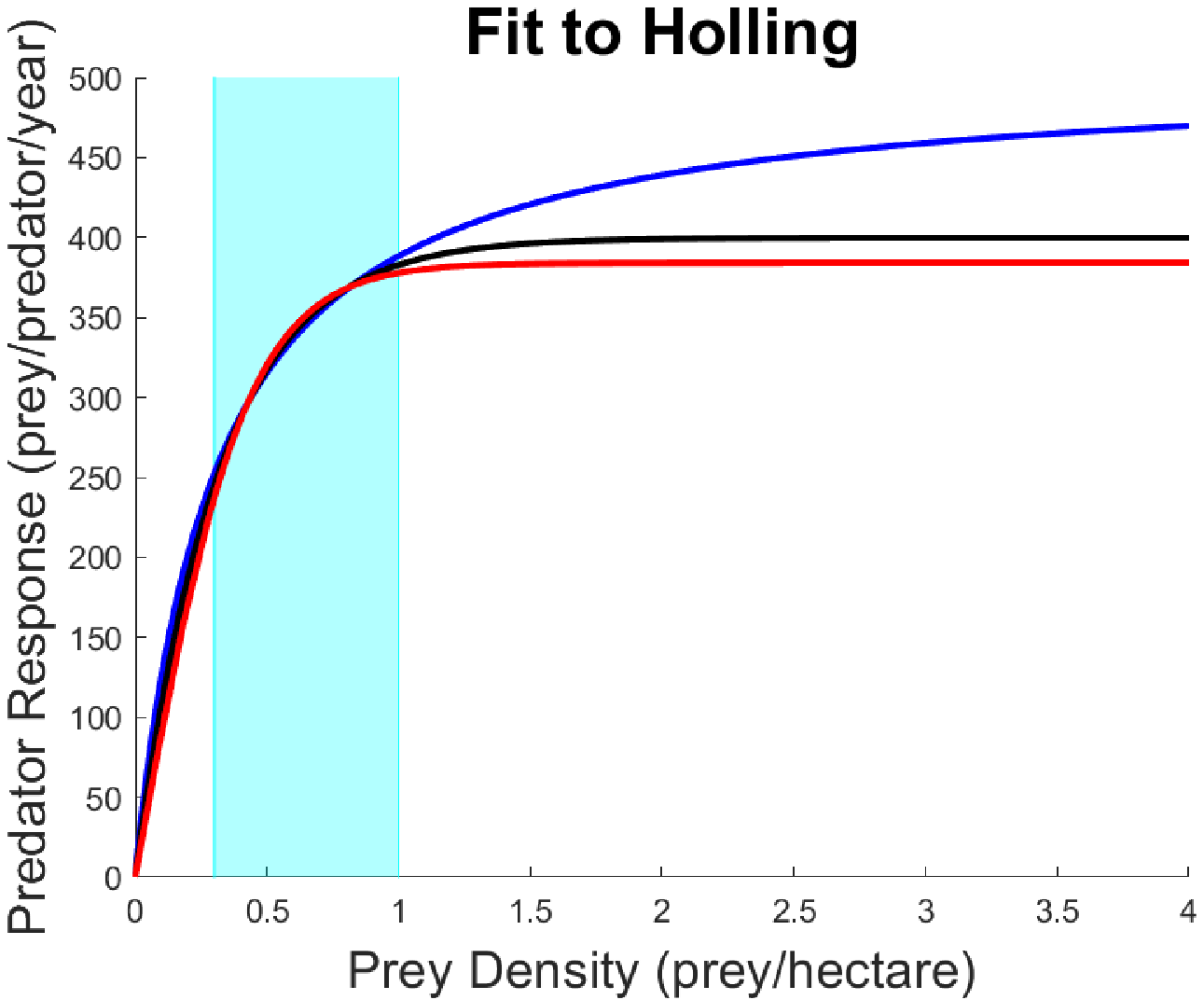}
        \subcaption{Fitted Functional Responses}
        \label{fig:SSCurvesMay}
    \end{subfigure}
    \begin{subfigure}{0.45\linewidth}
        \includegraphics[width=\linewidth]{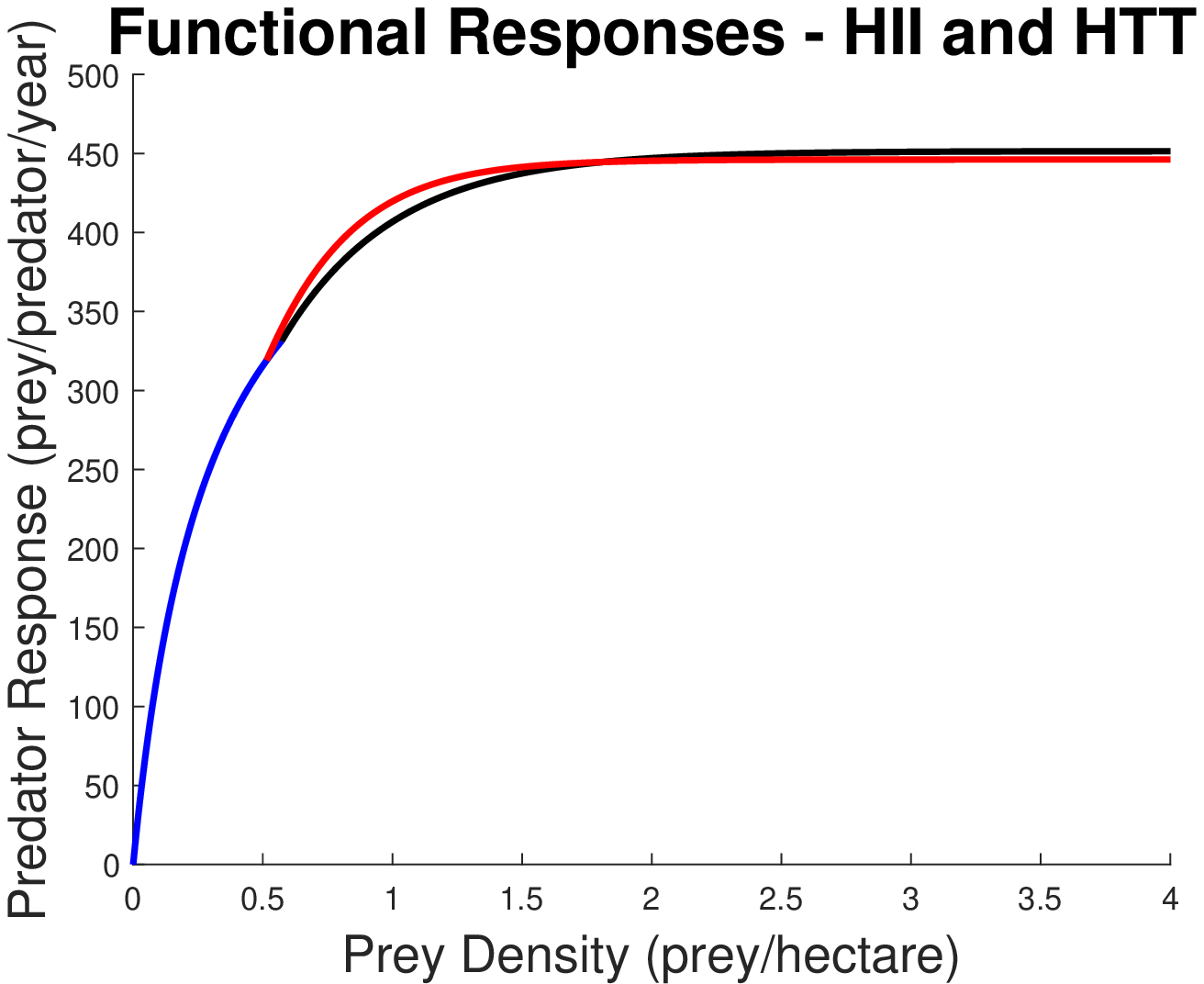}
        \subcaption{Piecewise Identical Functional Responses}
        \label{fig:SwitchingCurvesMay}
    \end{subfigure}
    \\
    \includegraphics[width=0.7\linewidth]{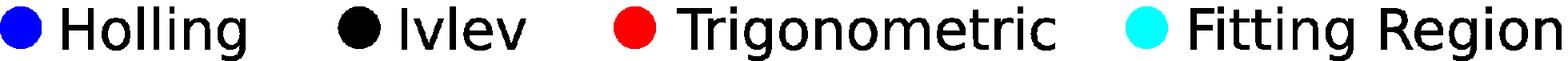}
    \caption{Functional responses used in the Leslie-Gower-May models: a) the original model; b) stochastic functional responses with $\sigma = 50$ where the shaded region represents one standard deviation; c) the functional responses with nonlinear least squares applied over prey densities in [0.3,1] shown by the cyan fitting region; and d) sample piecewise identical functional responses. Parameter values are given in Table \ref{tab:paramValues} and intersection values for the piecewise identical responses are given in Table \ref{tab:intersectionValues}. The Rosenzweig-MacArthur functional responses are similar with a difference in time scale and are given in the Supplementary Material (Figure 1). }
    \label{fig:MayCurves}
\end{figure}

In the case of the Rosenzweig-MacArthur model, the fixed functional response is the Ivlev response and we fit the Holling and trigonometric functional responses to the fixed Ivlev functional response. As expected, we recover the same parameter values reported in \citet{fussmann:2005}. In the case of the Leslie-Gower-May model, we use parameter values from \citet{vitense:2016} to define and fix the Holling functional response. Then, we fit the Ivlev and trigonometric functional responses to the fixed Holling functional response. The resulting Rosenzweig-MacArthur and Leslie-Gower-May functional responses appear in Figures S1a and \ref{fig:CurvesOrigMay}.

\begin{table}[ht]
    \centering
    \begin{tabular}{|c|c|c|c|c|c|c|}
    \hline 
   & \multicolumn{3}{|c|}{Rosenzweig-MacArthur} & \multicolumn{3}{|c|}{Leslie-Gower-May} \\
    \hline
    & Holling & Ivlev & Trigonometric & Holling & Ivlev & Trigonometric \\
    \hline
    $a_i$ & 3.05 & 1 & 0.99 & 1683.333 & 451.447 & 446.182 \\
    $b_i$ & 2.68 & 2 & 1.48 & 3.333 & 2.313 & 1.743\\
    \hline
    \end{tabular}
    \caption{Values of the parameters $a_i$ and $b_i$ used to fit the original functional responses. Note that the Rosenzweig-MacArthur model uses a functional response with units prey/predator/day and the Leslie-Gower-May model uses a functional response with units prey/predator/year.}
    \label{tab:paramValues}
\end{table}

\subsection{Structural Sensitivity}
%------------------------------------------------------------------------------------
To determine if the Leslie-Gower-May model is structurally sensitive, we plot a bifurcation diagram (shown in Figure \ref{fig:BifurMay}) using the carrying capacity , K, as the bifurcation parameter for each of the sub-models produced using the three functional responses \eqref{functionalResponses} in model \eqref{May}, and note the differences and similarities between them.

We find that the Leslie-Gower-May model is structurally sensitive: the Hopf bifurcation occurs at a very different carrying capacity for each sub-model. Figure \ref{fig:BifurMay} shows that the coexistence steady state in the model using the Holling functional response destabilizes at $K=5.55$, in the Ivlev model at $K=12.09$, and in the trigonometric model at $K=26.36$. Additionally, the Ivlev and trigonometric models have saddle node bifurcations of limit cycles at $K=10.51$ and $K=12.94$, respectively, and this bifurcation is not present in the Holling model. Following the formation of these limit cycles, bistability is observed in the Ivlev and trigonometric models. The Ivlev model is bistable over prey density values $K \in [10.51,12.09]$ and the trigonometric model is bistable over prey density values $K \in [12.09,26.36]$.

\subsection{Stochastic Model Methods} \label{Sect:stochMethods}
%------------------------------------------------------------------------------------

The addition of stochasticity to the functional responses makes them more like real data and increases their overlap. We add stochasticity to models \eqref{RM} and \eqref{May} by rewriting the model equations in the following form:
\begin{subequations}
\label{eq:stochastic-f}
\begin{align}
    d\xi &= -\xi(t)dt + \sigma dW, \label{SDE} \\
	\frac{dx}{dt} &= F(x, y, \xi(t)), \\
	\frac{dy}{dt} &= G(x, y, \xi(t)),
\end{align}
\end{subequations}
i.e., $\xi(t)$ is an Ornstein-Uhlenbeck process and $(x, y)$ solves a system
of ODEs driven by $\xi$. We use $\sigma$ to alter the amount of noise included in each time step. $F$ is the right-hand-side of equations \ref{RM1} and \ref{May1} (depending on the model in question) with $f_i(x)+\xi(t)$ in place of $f_i(x)$. Similarly, $G$ is the right-hand-side of equation \ref{RM2} with $f_i(x)+\xi(t)$ in place of $f_i(x)$. Note that $f_i(x)$ does not appear in equation \ref{May2}, so it is unchanged. Hence, model \ref{eq:stochastic-f} is a system of one stochastic differential equation (SDE) and two ordinary differential equations (ODEs). At each time step, we sample a Gaussian random variable, $dW = \sqrt{dt}Z$ where $Z \sim N(0,1)$, and solve the SDE (equation \ref{SDE}) using the Euler-Maruyama scheme to obtain $\xi(t)$. Then, we plug $\xi(t)$ into the system of ODEs and numerically solve the system using Euler's method with a step size small enough to achieve stable solutions. Note that if the solutions result in prey or predator densities less than zero due to the stochastic element, then the prey or predator density is set to zero. For the Rosenzweig-MacArthur model, we use $\sigma = 0.01$ and for the Leslie-Gower-May model, we use $\sigma = 50$ unless specified otherwise. Sample stochastic functional responses are given in Figures \ref{fig:StochCurves} and S1b. The corresponding results are given in Section \ref{Sect:StochResults}.

\subsection{Fitted Functional Response Methods} \label{Sect:CoexSSMethods}
%------------------------------------------------------------------------------------
The behaviour of both the Rosenzweig-MacArthur and the Leslie-Gower-May models is stable coexistence for values of the carrying capacity, $K$, less than the destabilizing bifurcation. We hypothesize that if the model behaviour is more similar leading up to the destabilizing bifurcation, then perhaps the destabilizing bifurcation will occur for more similar values of the carrying capacity. If our hypothesis is correct, it suggests that a key portion of the functional response is in the vicinity of the coexistence steady state. 

We choose an interval containing the prey density values of the coexistence steady state to perform the fit. We call this interval the fitting region and choose it such that it is smaller than the original $[0,4]$ fitting region. Note that we choose this region regardless of whether the coexistence steady state is stable or unstable. In models that do not use data, we choose a functional response to be fixed. We will refer to the sub-model where Holling is chosen as the fixed functional response as the fixed Holling case and the other sub-models follow the same naming scheme. Note that in the applied case of this method, we do not have to choose a fixed functional response since we can fit the functional responses to data. We then apply nonlinear least squares to the functional responses fit to either the fixed functional response or to data over the fitting region. The average prey density values of the coexistence steady states in the Rosenzweig-MacArthur and Leslie-Gower-May models are $x=0.05$ and $x=0.65$, respectively. So, the corresponding fitting regions are taken to be $x\in [0,0.1]$ and $x\in [0.3,1]$. The parameter values obtained from this fit to the fixed Holling functional response are given in Table \ref{table:ParametersLS}. The fixed Ivlev and fixed trigonometric functional responses are given in Figures S1c and S4. The corresponding parameter values are given in Table S1. The Holling fitted functional response results are given in Section \ref{sect:LSResults}.

\begin{table}[ht]
\centering
\begin{tabular}{|c|c|c|c|c|c|c|}
    \hline 
   & \multicolumn{3}{|c|}{Rosenzweig-MacArthur} & \multicolumn{3}{|c|}{Leslie-Gower-May} \\
    \hline
     & Holling & Ivlev & Trigonometric & Holling & Ivlev & Trigonometric \\
    \hline
    $a_i$ & 3.05 & 0.6282 & 0.3707 & 1683.333 & 399.8009 & 384.1465\\
    $b_i$ & 2.68 & 4.8204 & 7.6281 & 3.333 & 3.1656 & 2.4028\\
    \hline
\end{tabular}
    \caption{Values of the parameters obtained from the fitted functional response where the Holling functional response is fixed and the other functional responses are fit to the Holling response via the parameters $a_i$ and $b_i$.}
    \label{table:ParametersLS}
\end{table}

\subsection{Piecewise Identical Functional Response Methods} \label{Sect:SensRegionsMethods}
%------------------------------------------------------------------------------------
Here, rather than fitting the functional responses to each other, we modify the functional responses so that they are exactly the same over portions of the prey density domain $[0,4]$. Each functional response has three intersection points with each of the other two functional responses. All three functional responses intersect at (0,0), since there is zero predation at zero prey density. The values for the non-trivial intersections of the Rosenzweig-MacArthur and the Leslie-Gower-May functional response are given in Table \ref{tab:intersectionValues}. As a result, we can consider piecewise functional responses made up of contributions from two functional responses. These piecewise functional responses give us another way to ask if there is a key part of the functional response that determines the model behaviour. That is, we can examine how each portion of the functional response, (i.e., low, medium, and high prey density) affects the model's bifurcation structure. 

We consider three cases for each functional response for each model (18 cases total). To form each case we choose two functional responses and interchange the functional response we use for each of the three regions. The intersections noted earlier become the transition points between the regions. Sample piecewise functional responses are given in Figures \ref{fig:SwitchingCurvesMay} and S1d. In these figures, we introduce notation that is used from here forward to describe the piecewise identical sub-models. This notation is a three letter code that corresponds to each of the three prey density regions we have chosen to investigate, with Holling, Ivlev, and trigonometric functions each represented by their first letter. For example, HII represents the functional response that uses the Holling functional response for low prey densities, and uses the Ivlev functional response for medium and high prey densities. The corresponding results are given in Section \ref{Sect:SensRegionsResults}.

\begin{table}[ht]
    \centering
    \begin{tabular}{|c|c|c|c|}
    \hline
    \textbf{Rosenzweig-MacArthur} & \textbf{Holling/Ivlev} & \textbf{Ivlev/Trig} & \textbf{Trig/Holling}\\
    \hline
    Extinction & 0 & 0 & 0\\
    Low to medium prey density & 0.6191 & 0.5651 & 0.5942\\
    Medium to high prey density & 2.5799 & 2.1597 & 2.4695\\
    \hline
    \textbf{Leslie-Gower-May} & \textbf{Holling/Ivlev} & \textbf{Ivlev/Trig} & \textbf{Trig/Holling}\\
    \hline
    Extinction & 0 & 0 & 0 \\
    Low to medium prey density & 0.5726 & 0.4458 & 0.5152 \\
    Medium to high prey density & 2.4486 & 2.2611 & 1.8077 \\
    \hline
    \end{tabular}
    \caption{Prey density values for the intersections in the functional responses used to form the piecewise identical functional responses.}
    \label{tab:intersectionValues}
\end{table}

\section{Results}
%------------------------------------------------------------------------------------
\label{Results}
\subsection{Stochastic Model Results}
\label{Sect:StochResults}
%------------------------------------------------------------------------------------

In the case of the Rosenzweig-MacArthur model, we added Gaussian stochasticity with $\sigma = 0.01$. We choose this amplitude because stochasticity of smaller amplitude had virtually no effect on the bifurcation structure. $\sigma = 0.01$ stochasticity adds some variance to the time series but does not significantly alter the qualitative behaviour for carrying capacities less than the Hopf bifurcation in the deterministic model. That is, given initial conditions that lead to the stable coexistence steady state in the deterministic case, the stochastic model also results in trajectories that converge to the coexistence steady state. For values of the carrying capacity greater than that of the bifurcation in the deterministic model, the lower side of the limit cycle reaches very low prey densities of approximately $10^{-60}$ at times ($10^{-58}$ times smaller than the stochasticity). As a result, the stochasticity causes extinction within 100 time steps in most of the time series of the stochastic model for carrying capacities greater than that of the deterministic model's destabilizing bifurcation. These results mean that we are unable to clearly observe the stochastic behaviour and we cannot conclude whether or not stochasticity reduces structural sensitivity in the Rosenzweig-MacArthur model.

In the case of the Leslie-Gower-May model, we added Gaussian stochasticity with $\sigma = 50$. As in the case of the Rosenzweig-MacArthur model, we chose this amplitude of stochasticity because it is large enough to affect the bifurcation structure, but not so large that all of trajectories rapidly terminate in extinction of the predator or both populations. Results obtained from both smaller and larger amplitudes of stochasticity are provided in the Supplementary Material (Figure S3).

\begin{figure}[ht]
    \centering
    \begin{subfigure}{\linewidth}
        \includegraphics[width=0.32\linewidth]{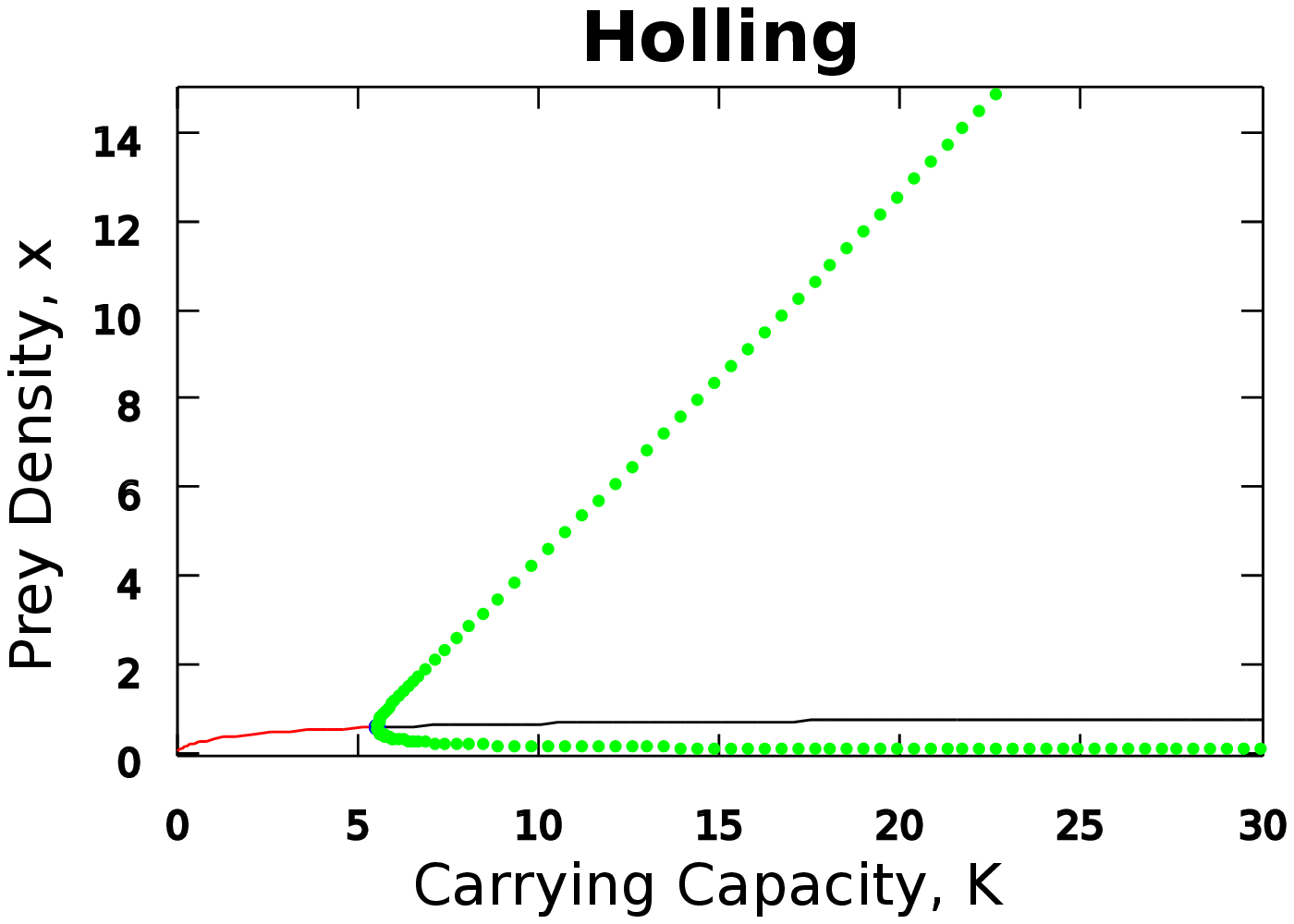}
        \includegraphics[width=0.32\linewidth]{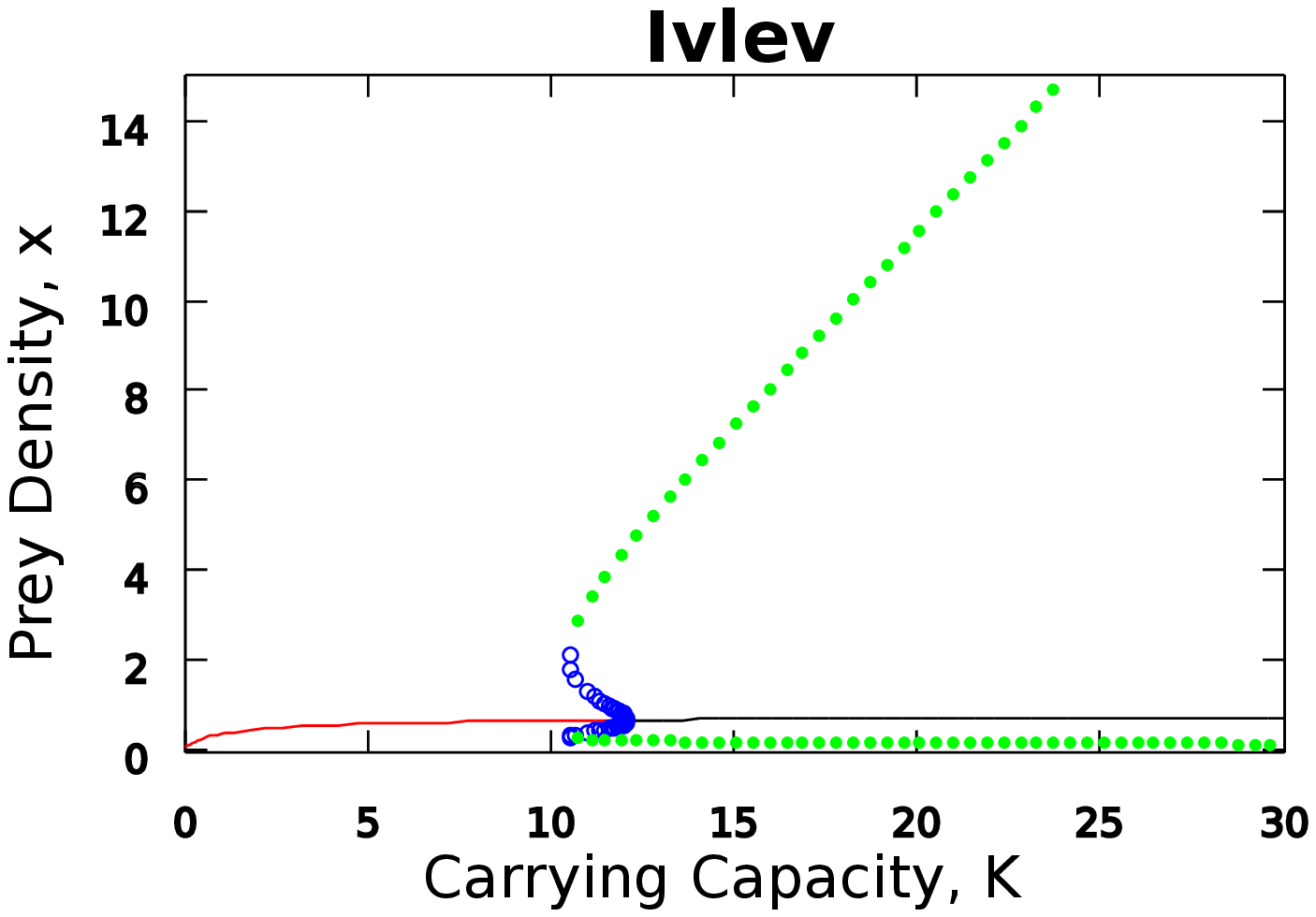}
        \includegraphics[width=0.32\linewidth]{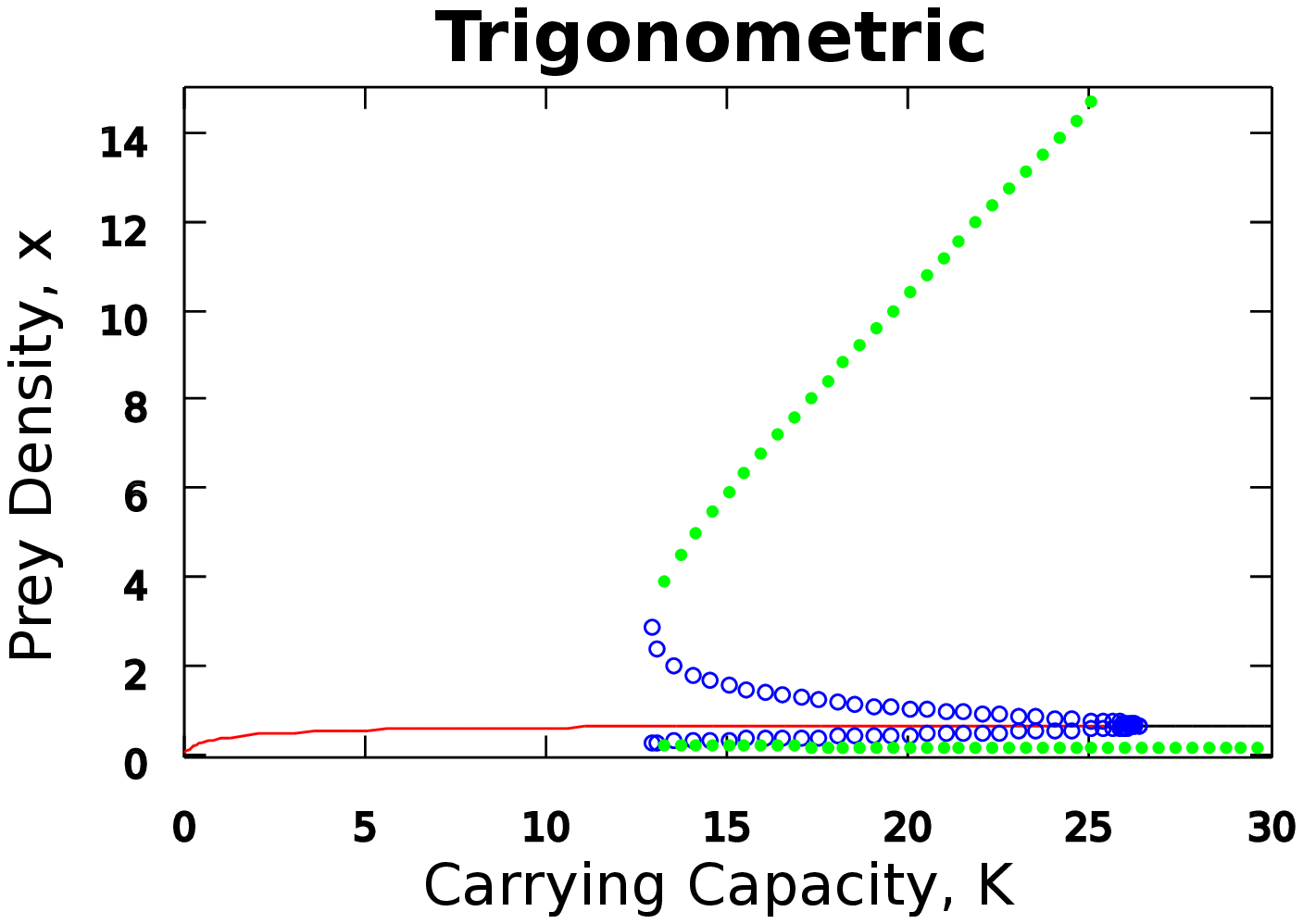}
        \subcaption{Original Functional Response}
        \label{fig:BifurMay}
    \end{subfigure}
    \\
    \begin{subfigure}{\linewidth}
        \includegraphics[width=0.32\linewidth]{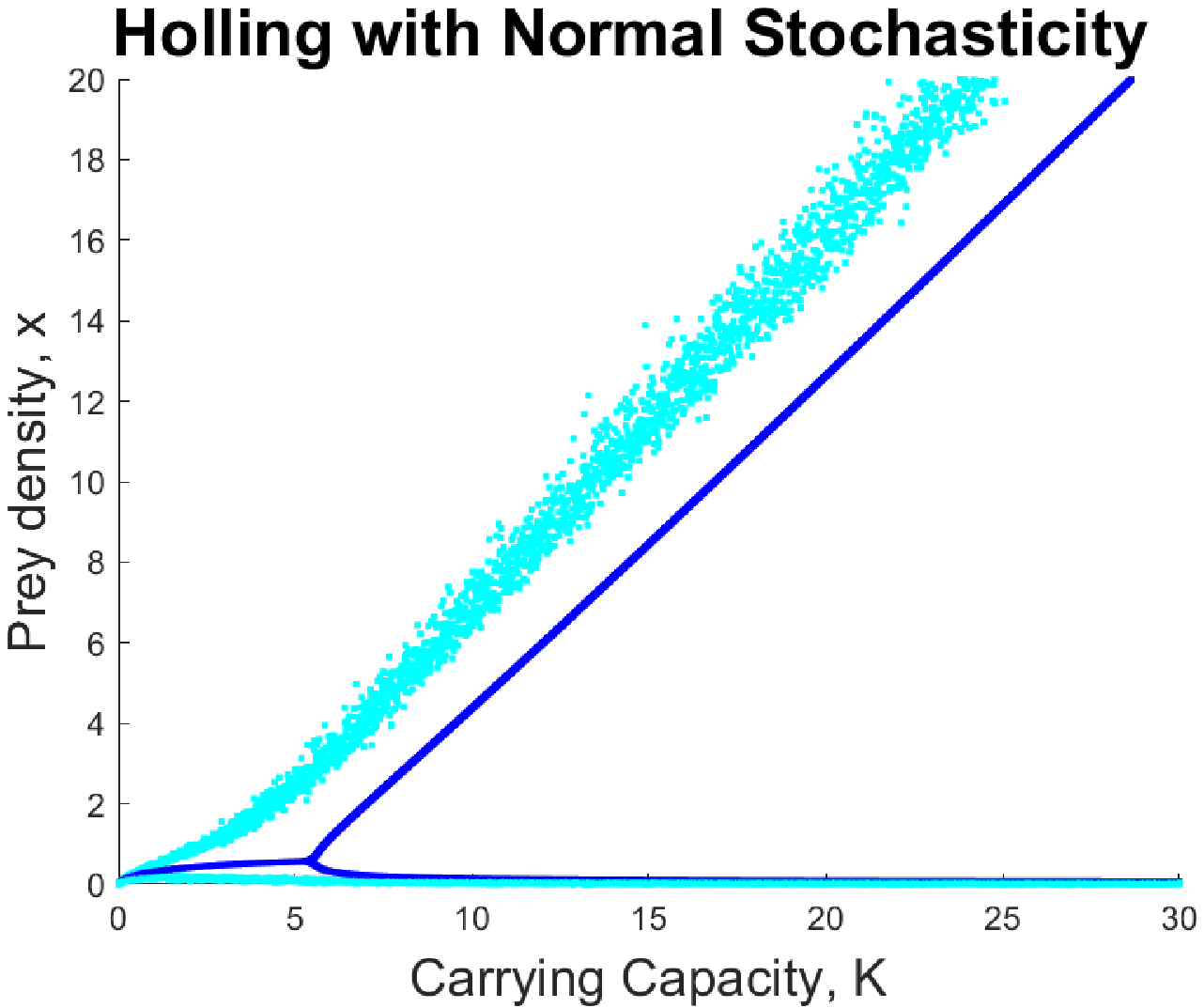}
        \includegraphics[width=0.32\linewidth]{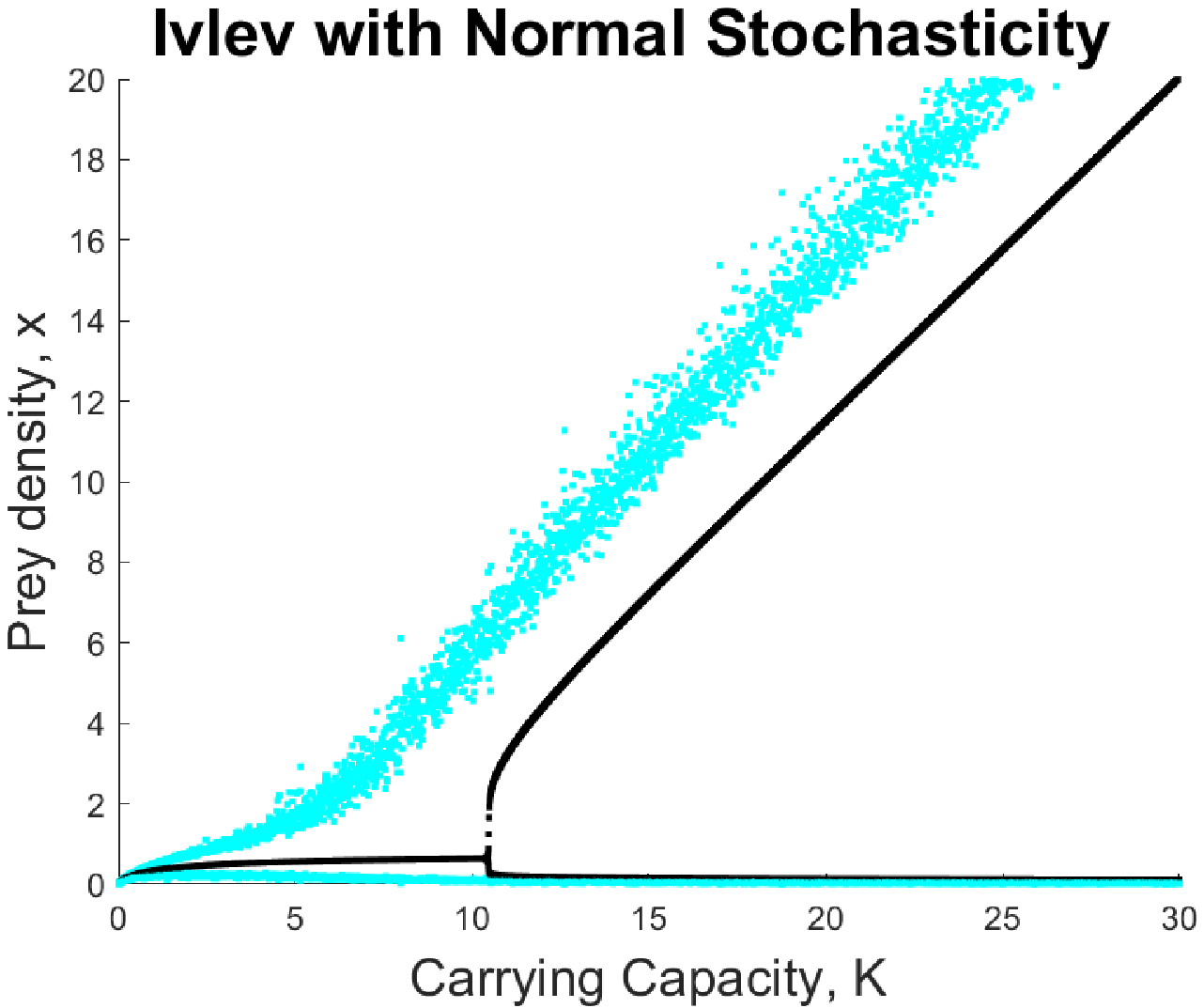}
        \includegraphics[width=0.32\linewidth]{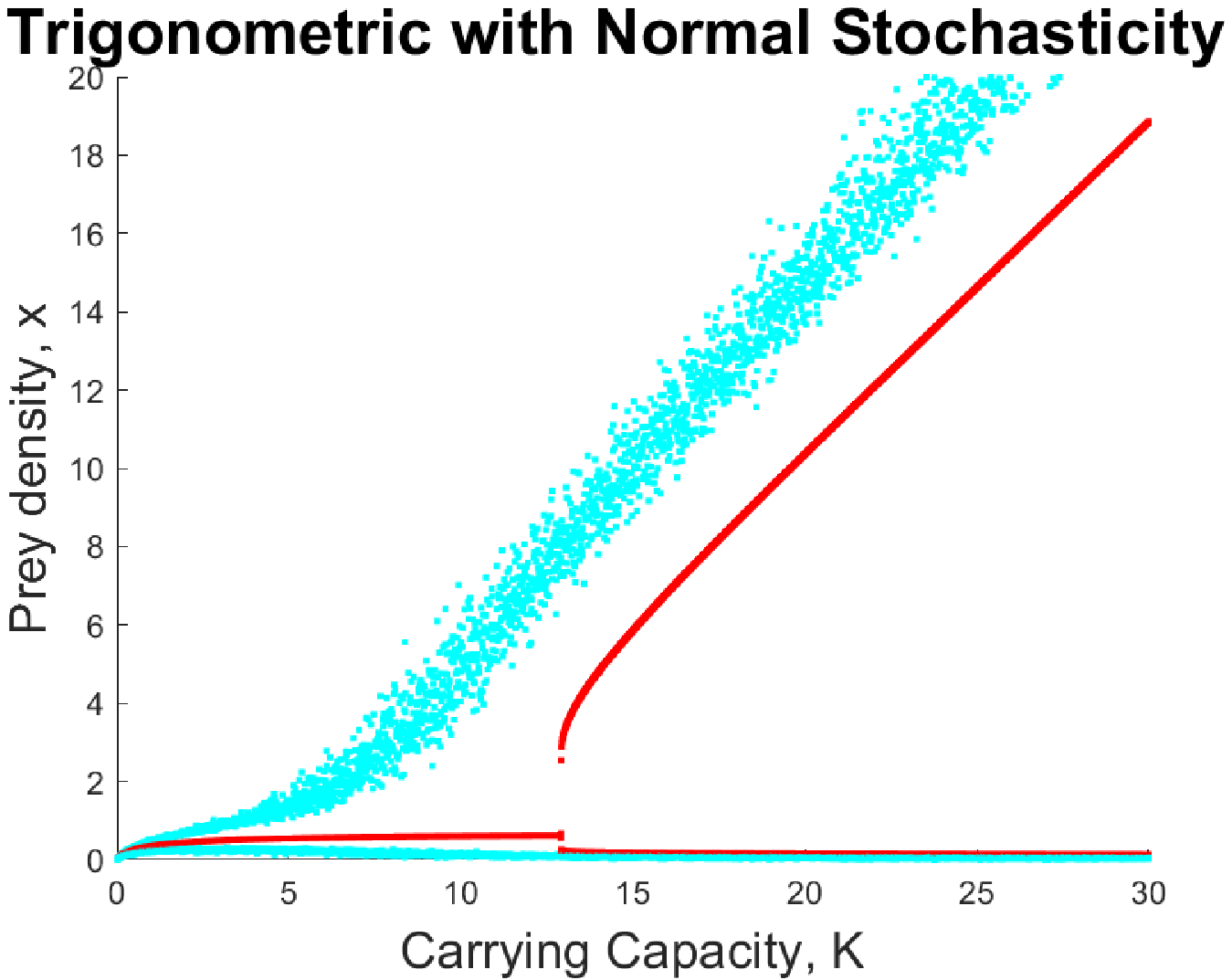}
        \\ 
        \includegraphics[width=\linewidth]{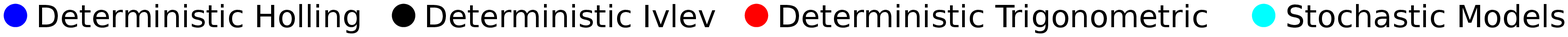}
        \subcaption{Stochastic Functional Response}
        \label{fig:StochBifur}
    \end{subfigure}
    \\
    \begin{subfigure}{\linewidth}
        \includegraphics[width=0.32\linewidth]{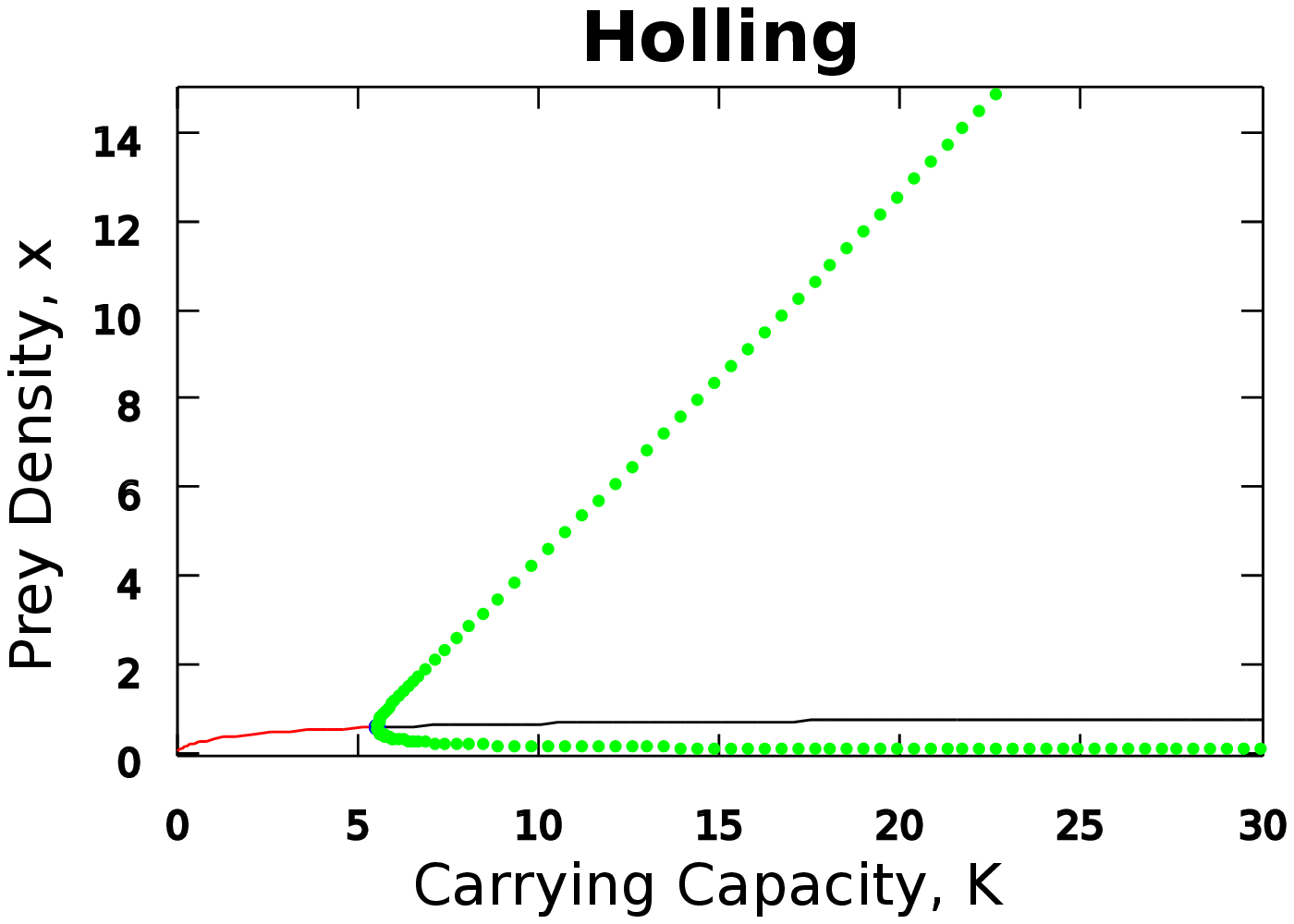}
        \includegraphics[width=0.32\linewidth]{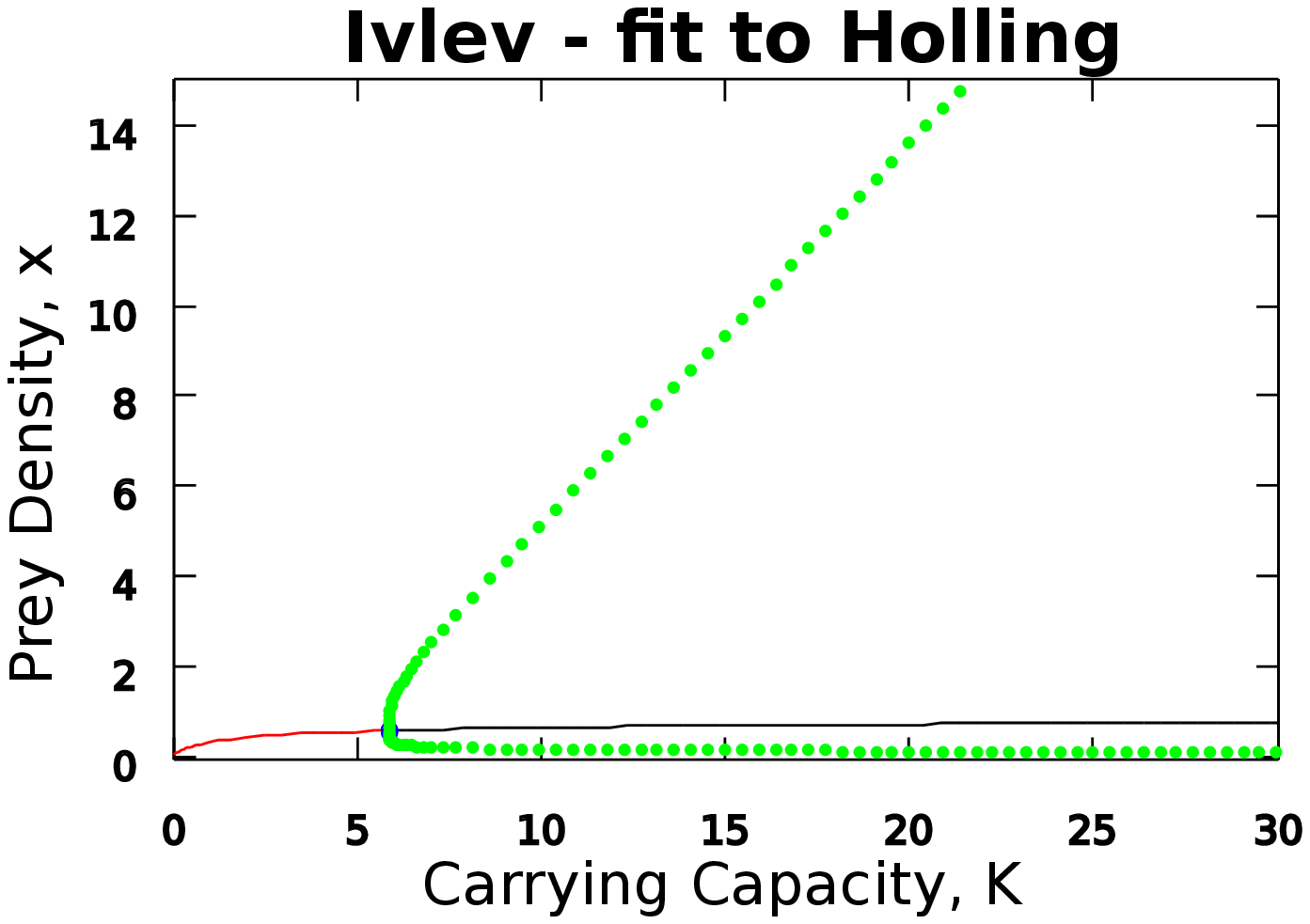}
        \includegraphics[width=0.32\linewidth]{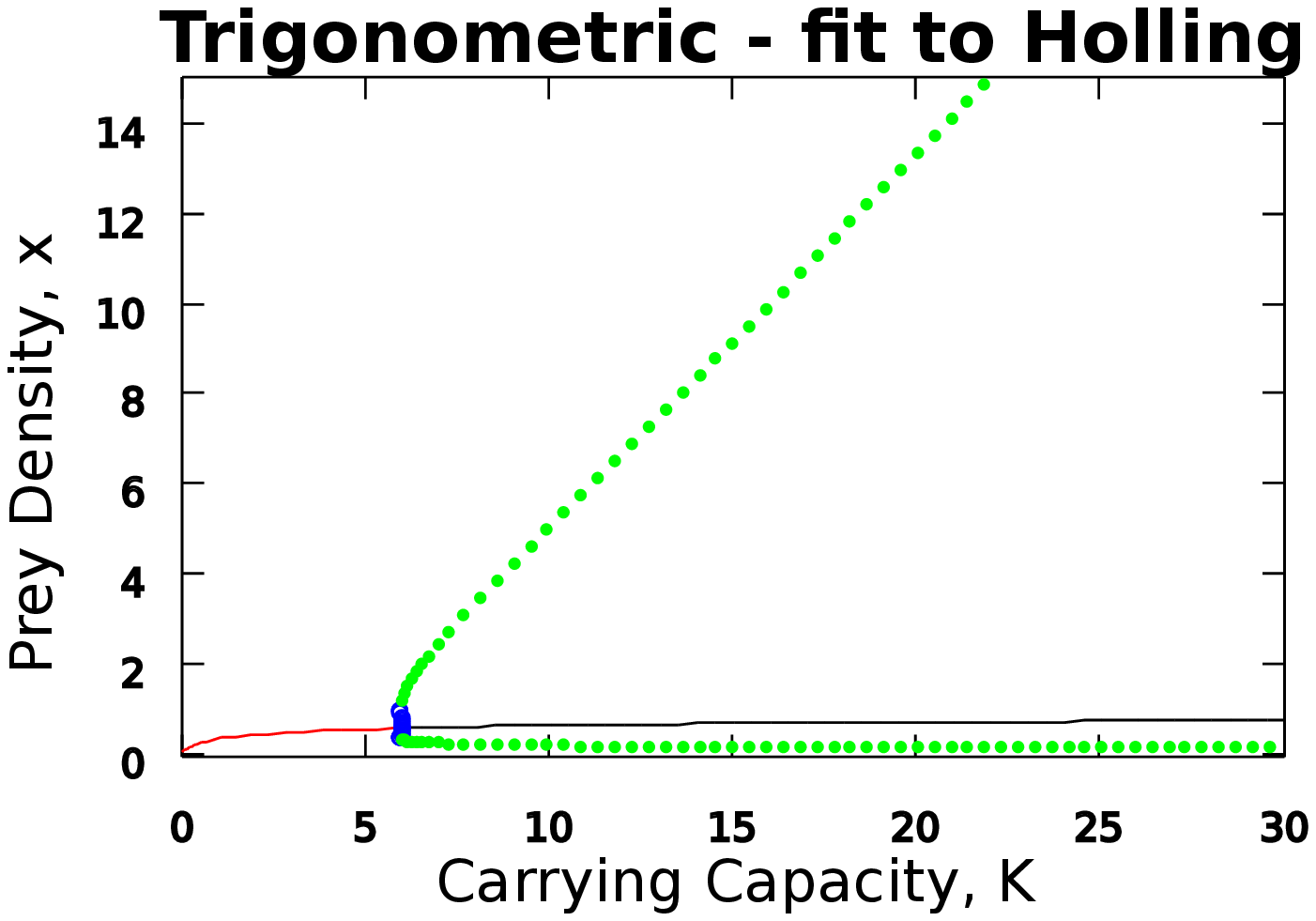}
        \subcaption{Fitted Functional Response}
        \label{fig:CoexSSBifurMay}
    \end{subfigure}
    \\
    \begin{subfigure}{\linewidth}
        \includegraphics[width=0.32\linewidth]{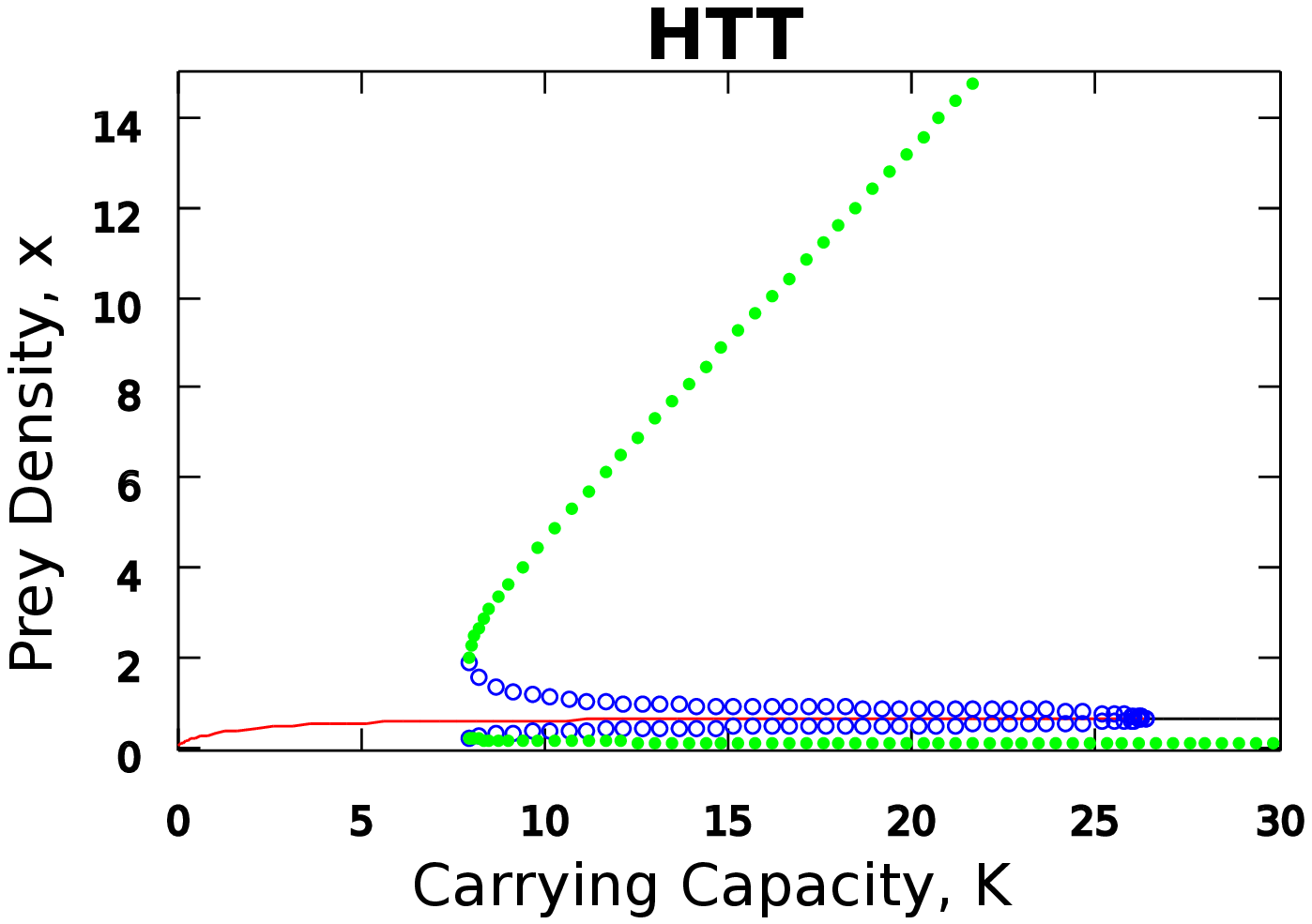}
        \includegraphics[width=0.32\linewidth]{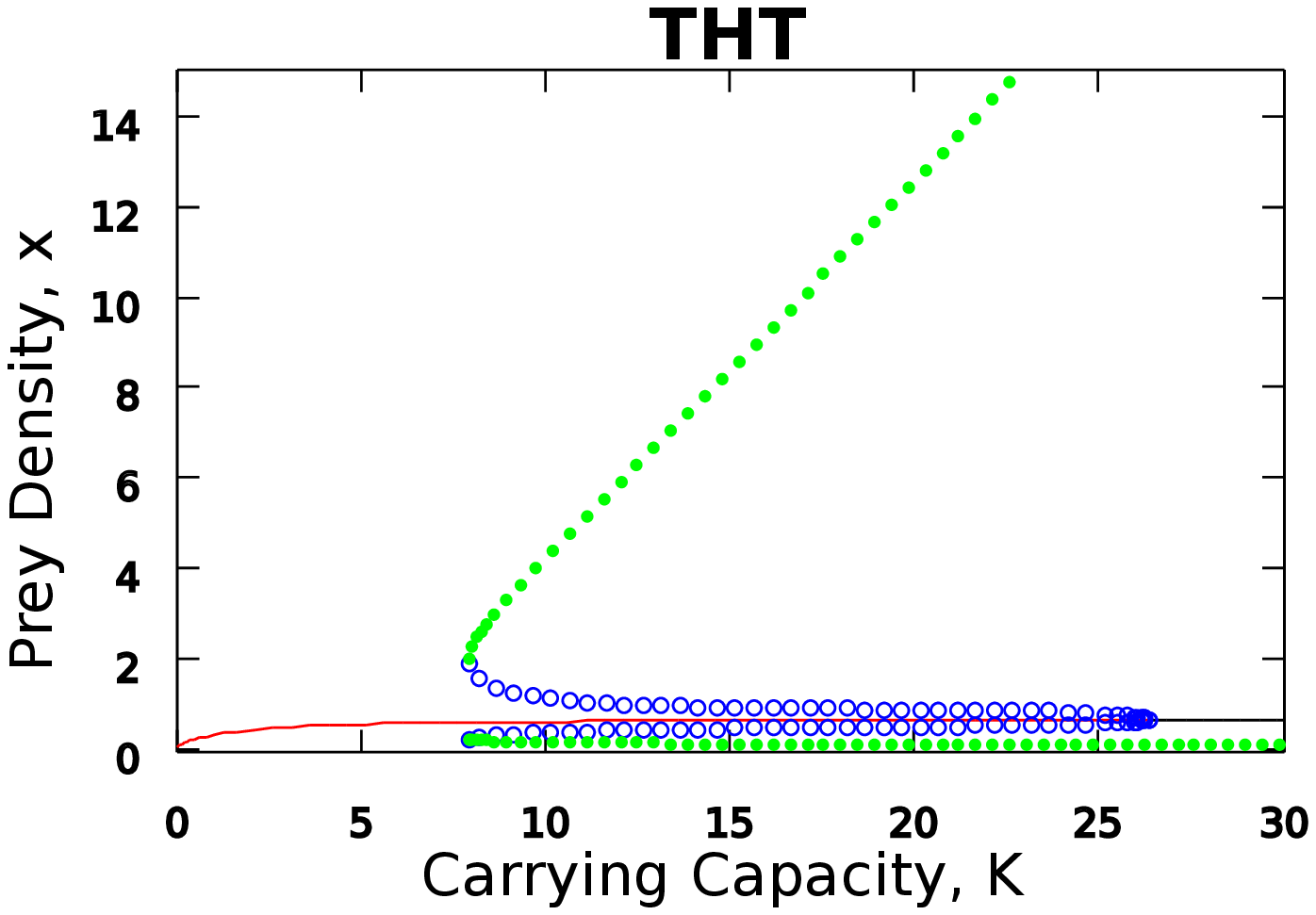}
        \includegraphics[width=0.32\linewidth]{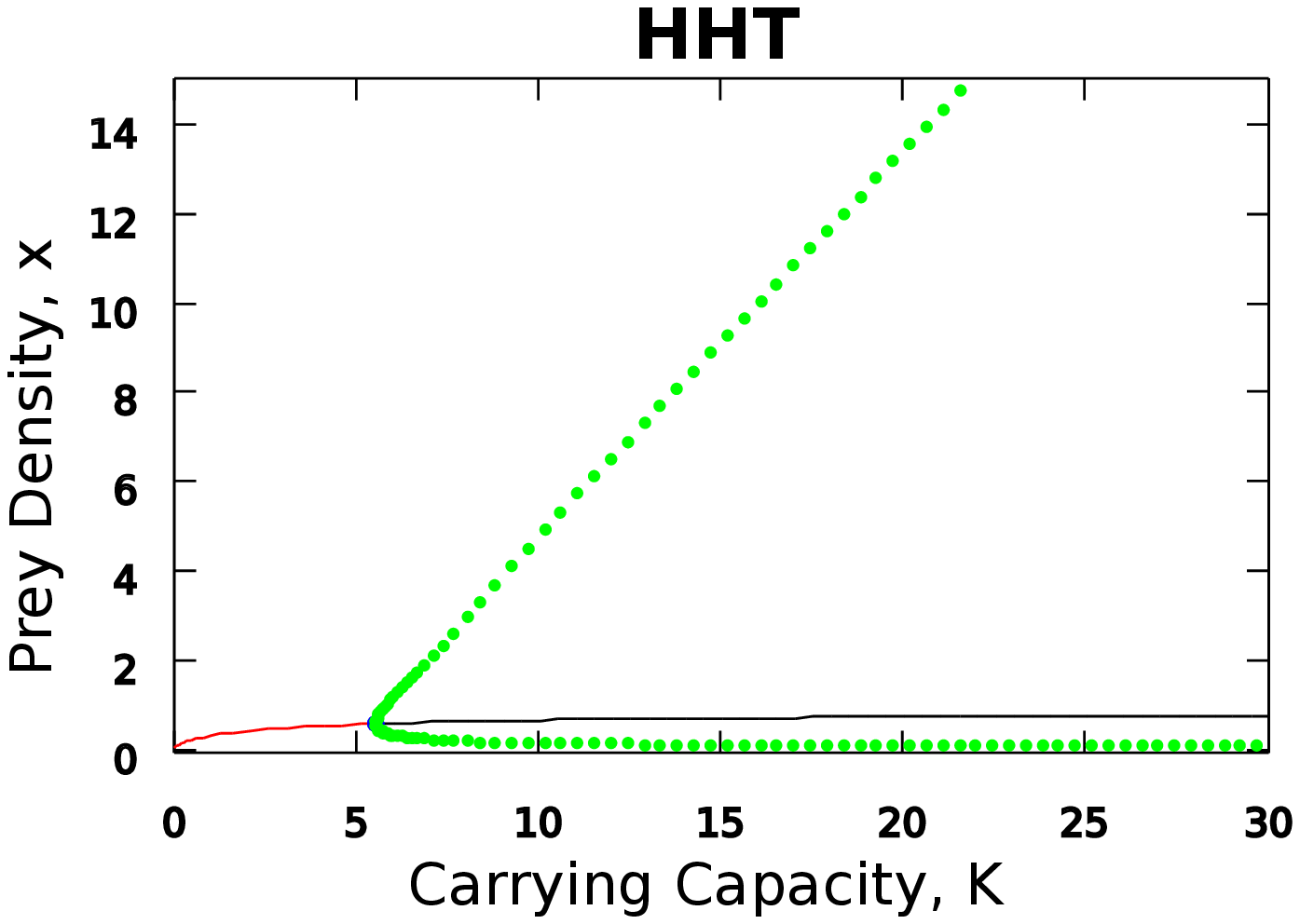}
        \subcaption{Piecewise Identical Functional Response}
        \label{fig:SwitchingBifurMay}
    \end{subfigure}
    \\
    \includegraphics[width=\linewidth]{Pictures/FIGLeg2.eps}
    \caption{Bifurcation diagrams for the Leslie-Gower-May model \eqref{May}: a) the original model; b) our stochastic Leslie-Gower-May model with $\sigma = 50$ where both the deterministic and stochastic models were obtained by running time series and discarding the transient solutions, so these diagrams only show the quasi-steady state behaviour; c) sub-models using least squares applied to the functional responses over the prey density interval $[0.3,1]$ where the Holling functional response is fixed; and d) sub-models using piecewise identical functional responses. The corresponding functional responses are given in Figure \ref{fig:MayCurves} and the bifurcation values are given in Table \ref{table:BifurValues}.}
\end{figure}

The bifurcation diagrams (quasi-steady state behaviour) for the stochastic Leslie-Gower-May model are given in Figure \ref{fig:StochBifur}. This figure shows better similarity in the bifurcation structures with the addition of stochasticity. The bifurcation structures are obtained by plotting the maximum and minimum of the stochastic time series (after removal of transients). Thus, these windows show the smallest value of the carrying capacity for which the coexistence steady states may destabilize. Limit cycles in our stochastic model appear to have greater amplitude than those of the deterministic model. We find that the stochastic bifurcation structures show better similarity than the deterministic model's bifurcation structures. This similarity suggests that there is a decrease in structural sensitivity of the model.

%The bifurcation diagrams (quasi-steady state behaviour) for the stochastic Leslie-Gower-May model are given in Figure \ref{fig:StochBifur}. This figure shows that the deterministic structure remains and is robust to the addition of stochasticity. Limit cycles in our stochastic model appear to have greater amplitude than those of the deterministic model. This result is an artefact of our plotting procedure: the bifurcation plots are obtained by plotting the maximum and minimum of the stochastic time series (after removal of transients). We find that the average value of the stochastic trajectories converge to the deterministic model's behaviour. This similarity suggests that there is no decrease in structural sensitivity of the model.

Plots of time series (Figure S2) indicate that bistability remains a property of the stochastic Ivlev and trigonometric models. In the bistable regions, the solution trajectories of the Ivlev and trigonometric models alternate randomly between the coexistence steady state and the stable limit cycle, and also spend some time at the unstable limit cycle. This ability to alternate between stable and unstable steady states is unique to the stochastic model, and has been observed in other work \citep{abbott:2015, abbott:2017}. These observations occur for small to intermediate stochastic amplitudes, larger amplitudes lead to rapid extinction, and more generally, to the stochastic behaviour washing out the deterministic structure. 

\begin{figure}[ht]
    \centering
    \begin{subfigure}{\linewidth}
        \includegraphics[width=0.32\linewidth]{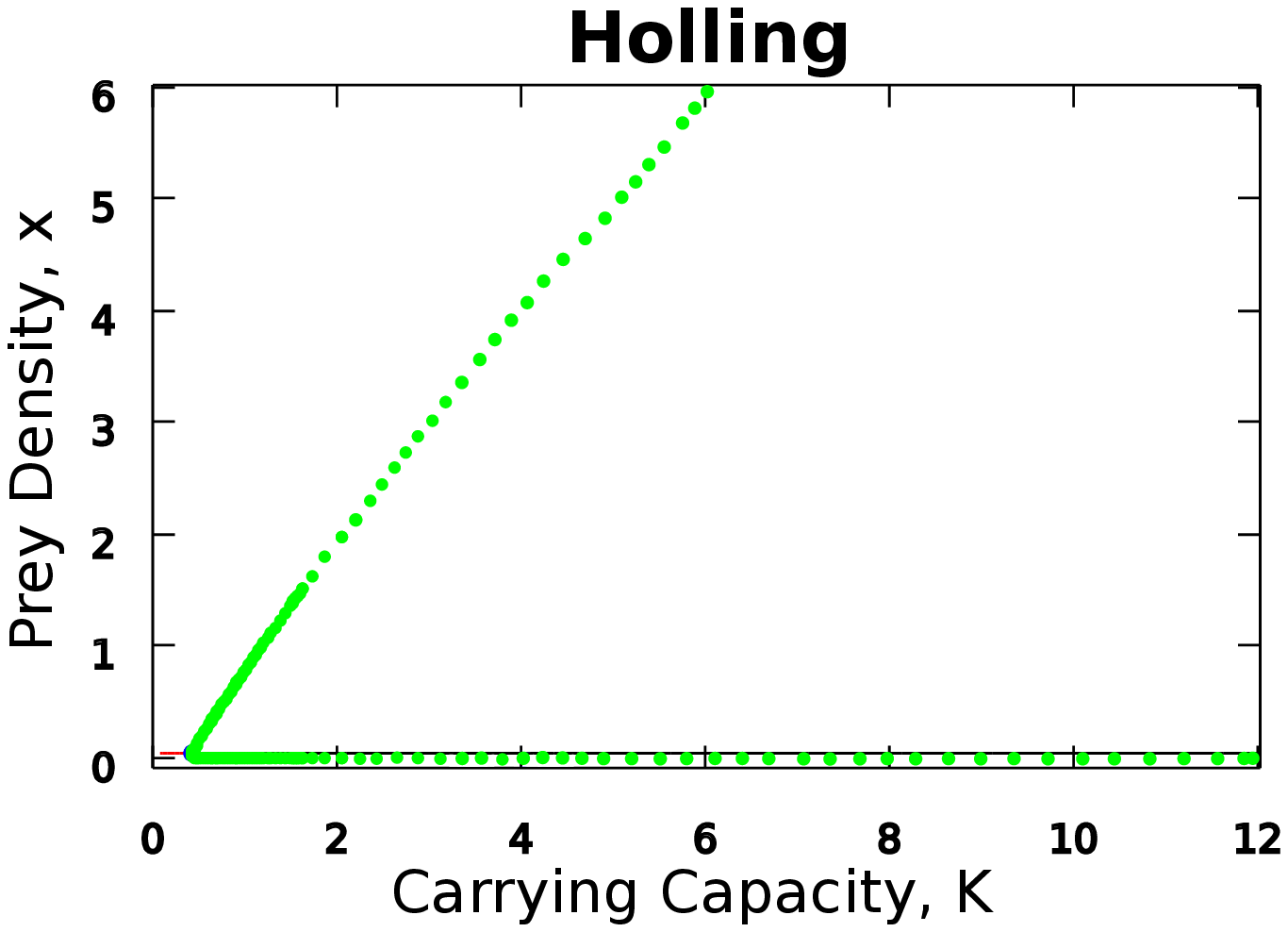}
        \includegraphics[width=0.32\linewidth]{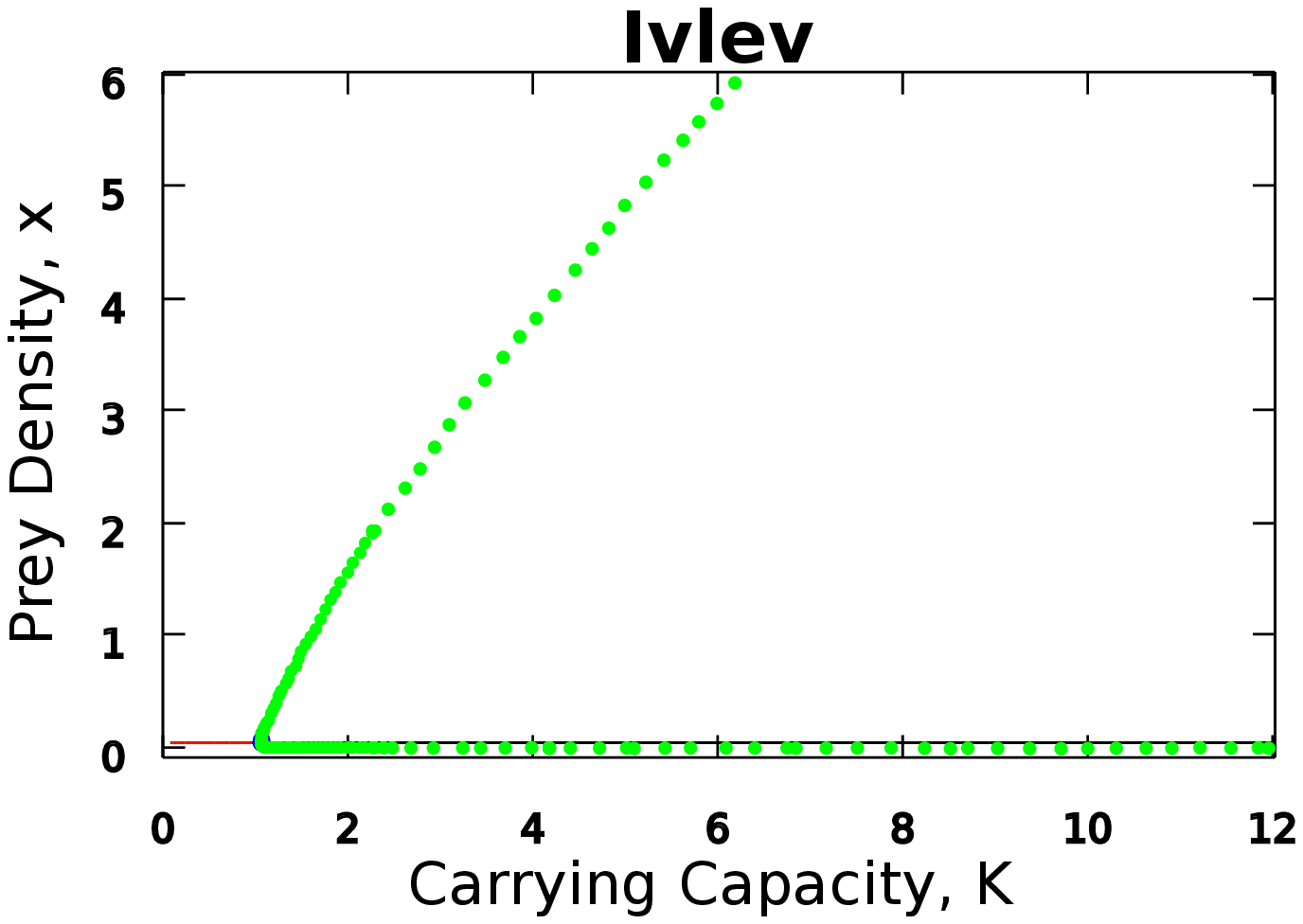}
        \includegraphics[width=0.32\linewidth]{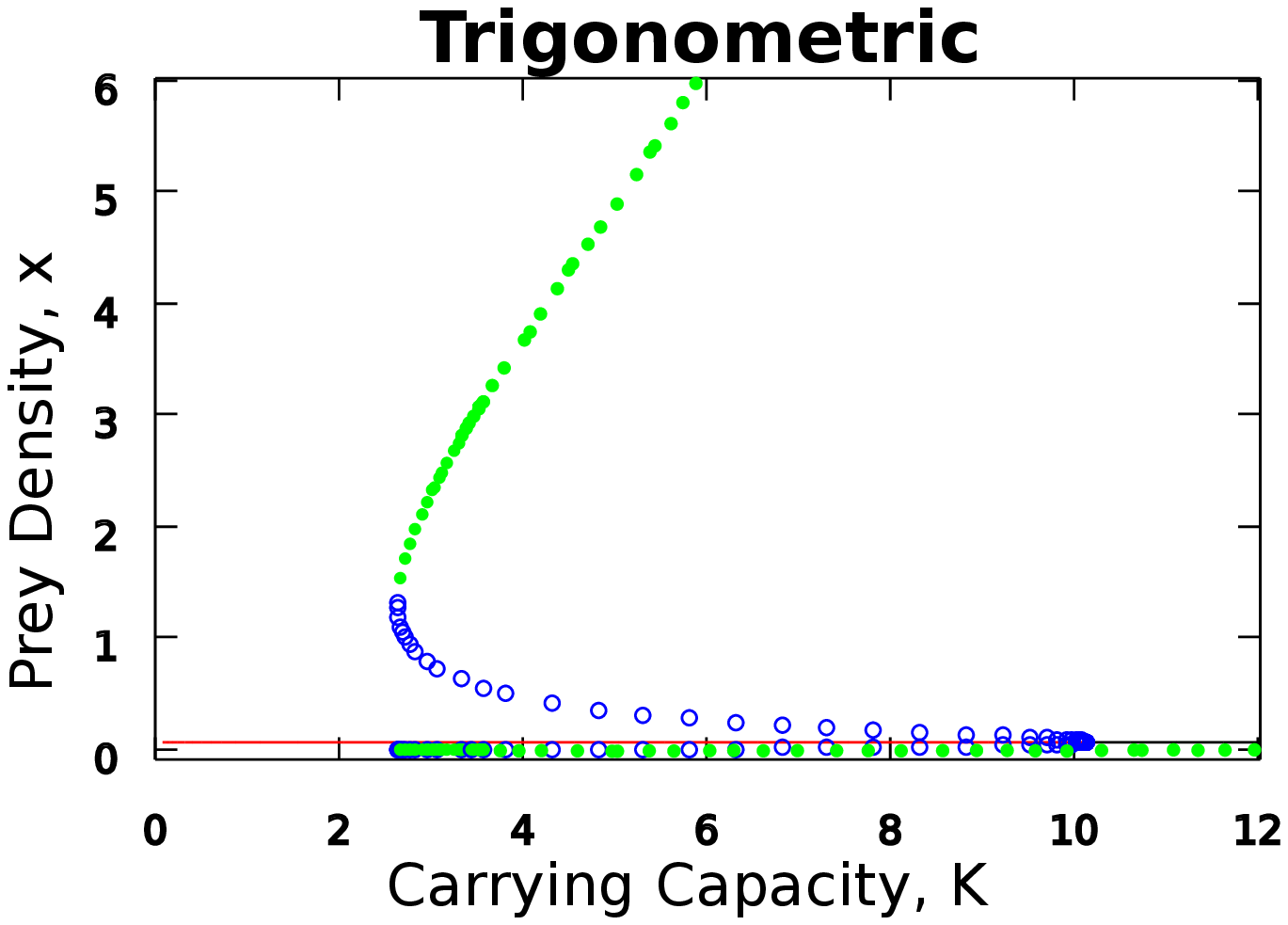} 
        \subcaption{Original Functional Response}
        \label{fig:BifurRM}
    \end{subfigure}
    \\
    \begin{subfigure}{\linewidth}
        \includegraphics[width=0.32\linewidth]{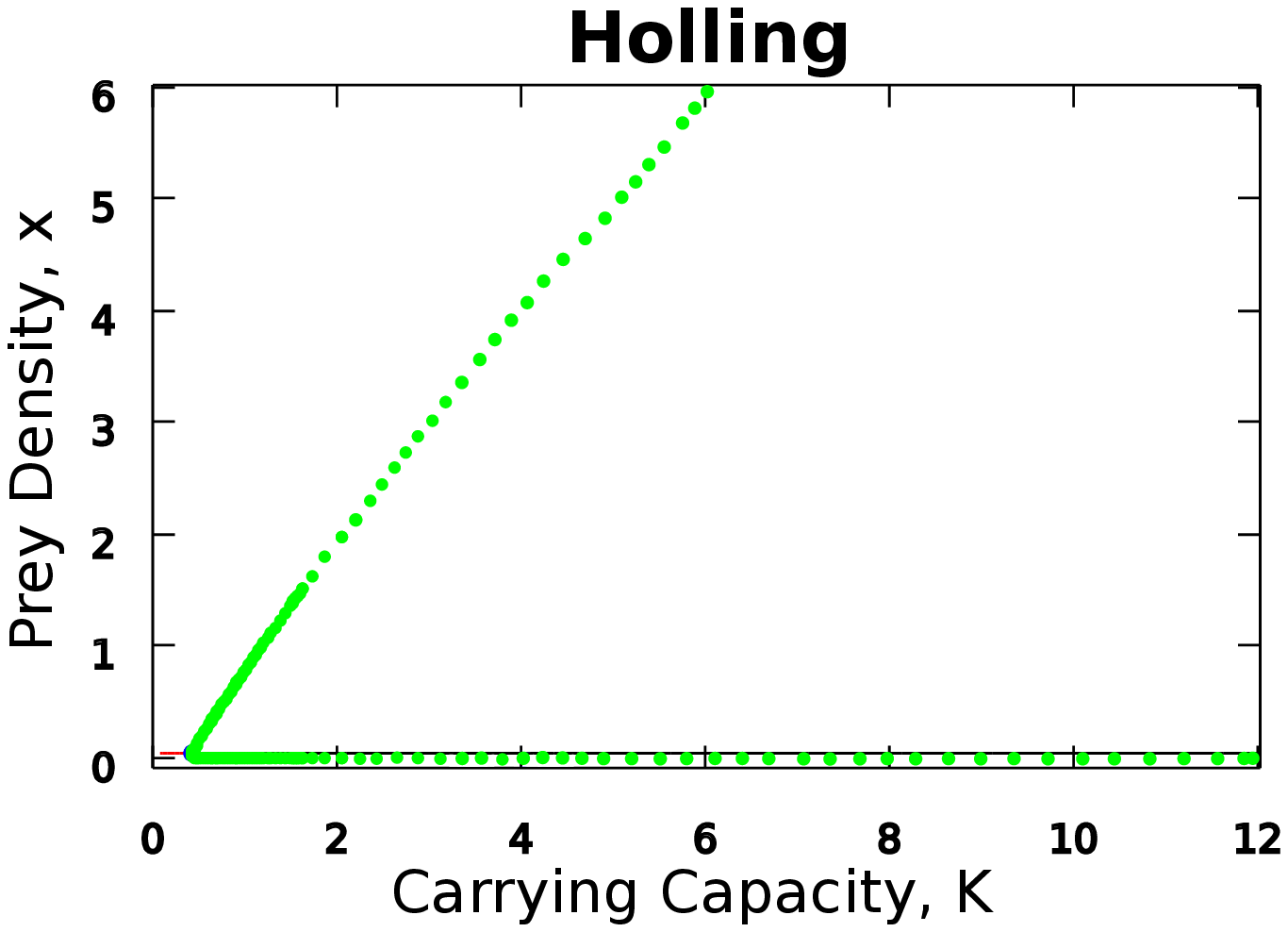}
        \includegraphics[width=0.32\linewidth]{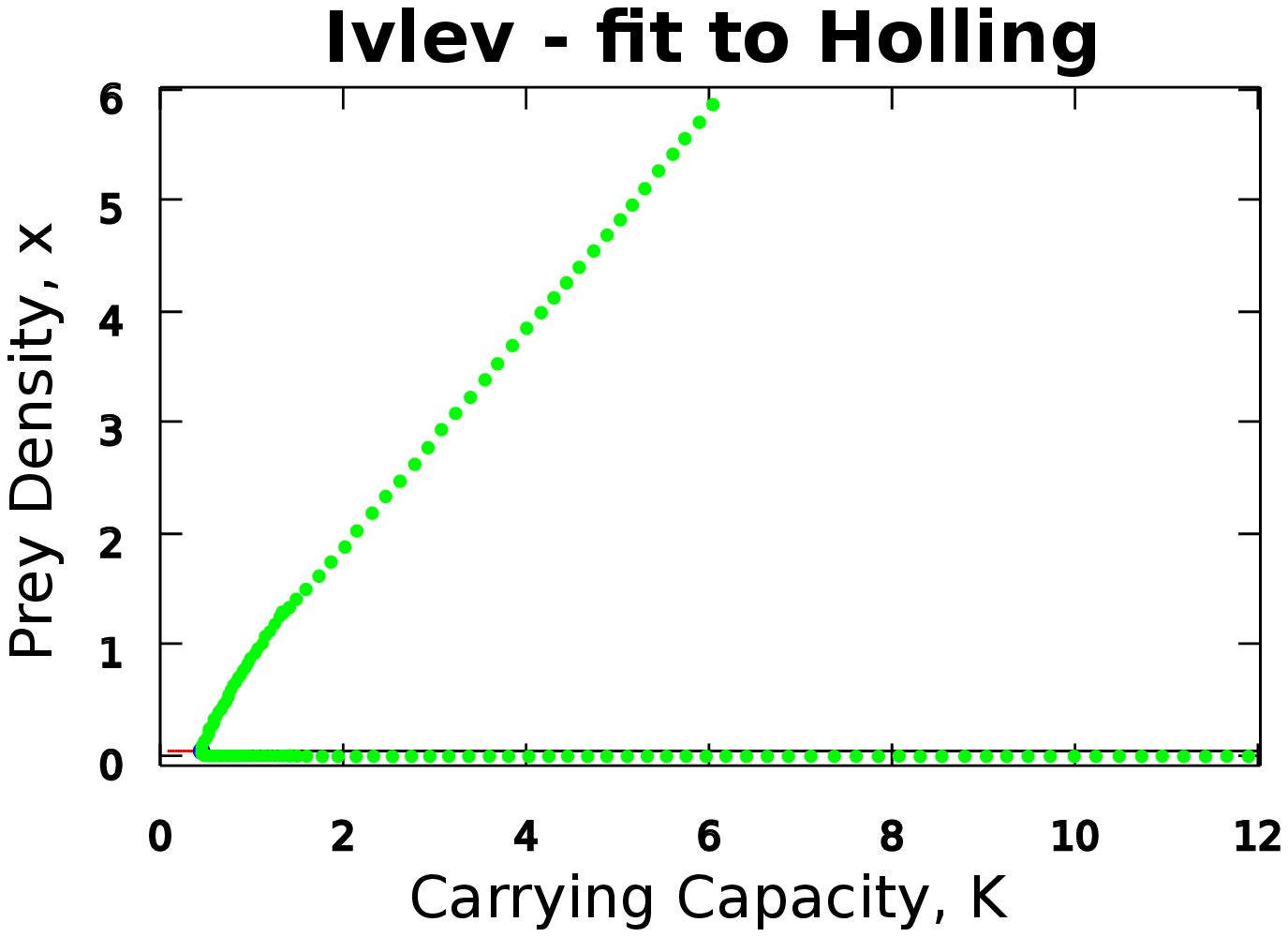}
        \includegraphics[width=0.32\linewidth]{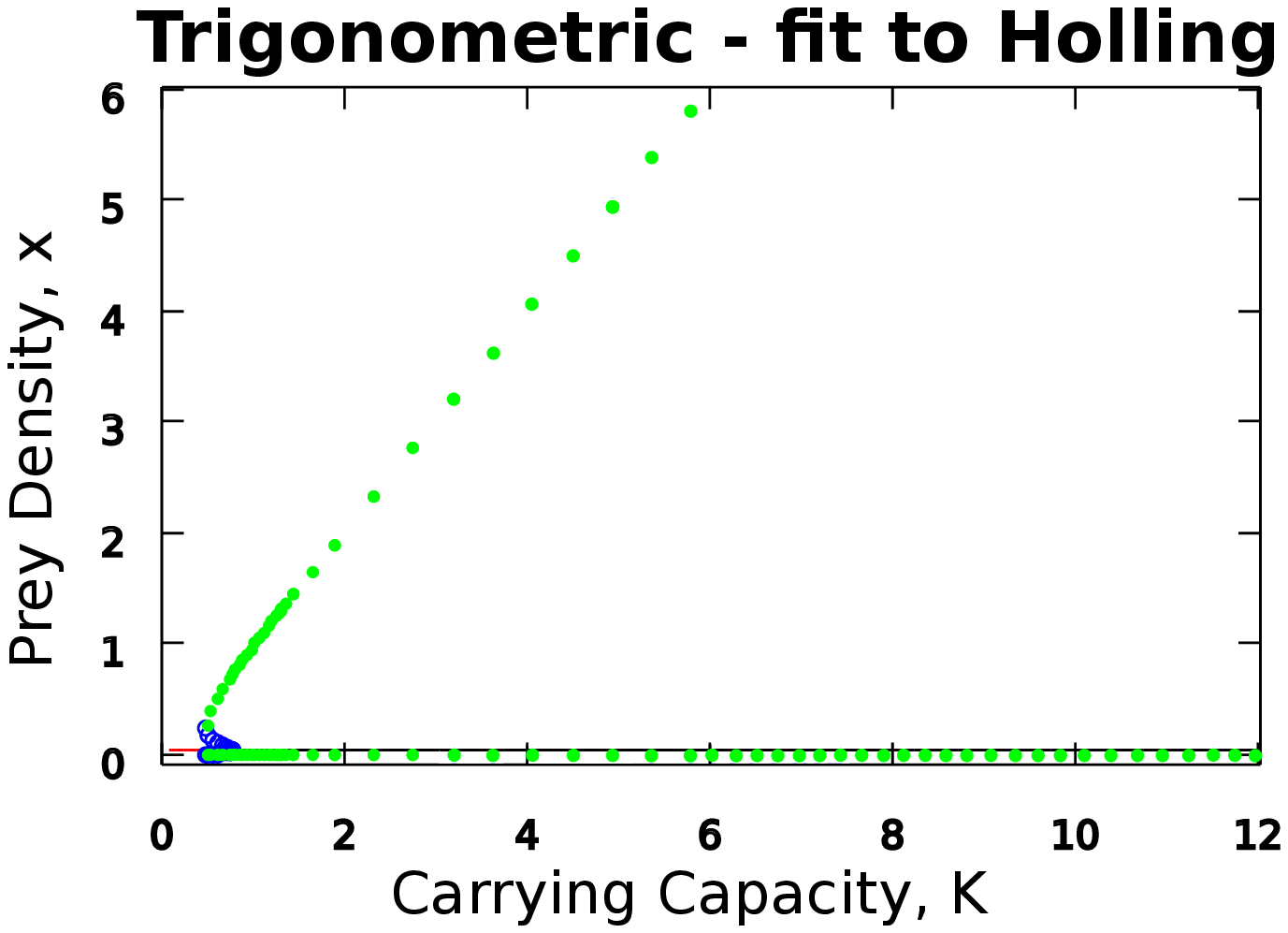}
        \subcaption{Fitted Functional Response}
        \label{fig:CoexSSBifurRM}
    \end{subfigure}
    \\
    \begin{subfigure}{\linewidth}
        \includegraphics[width=0.32\linewidth]{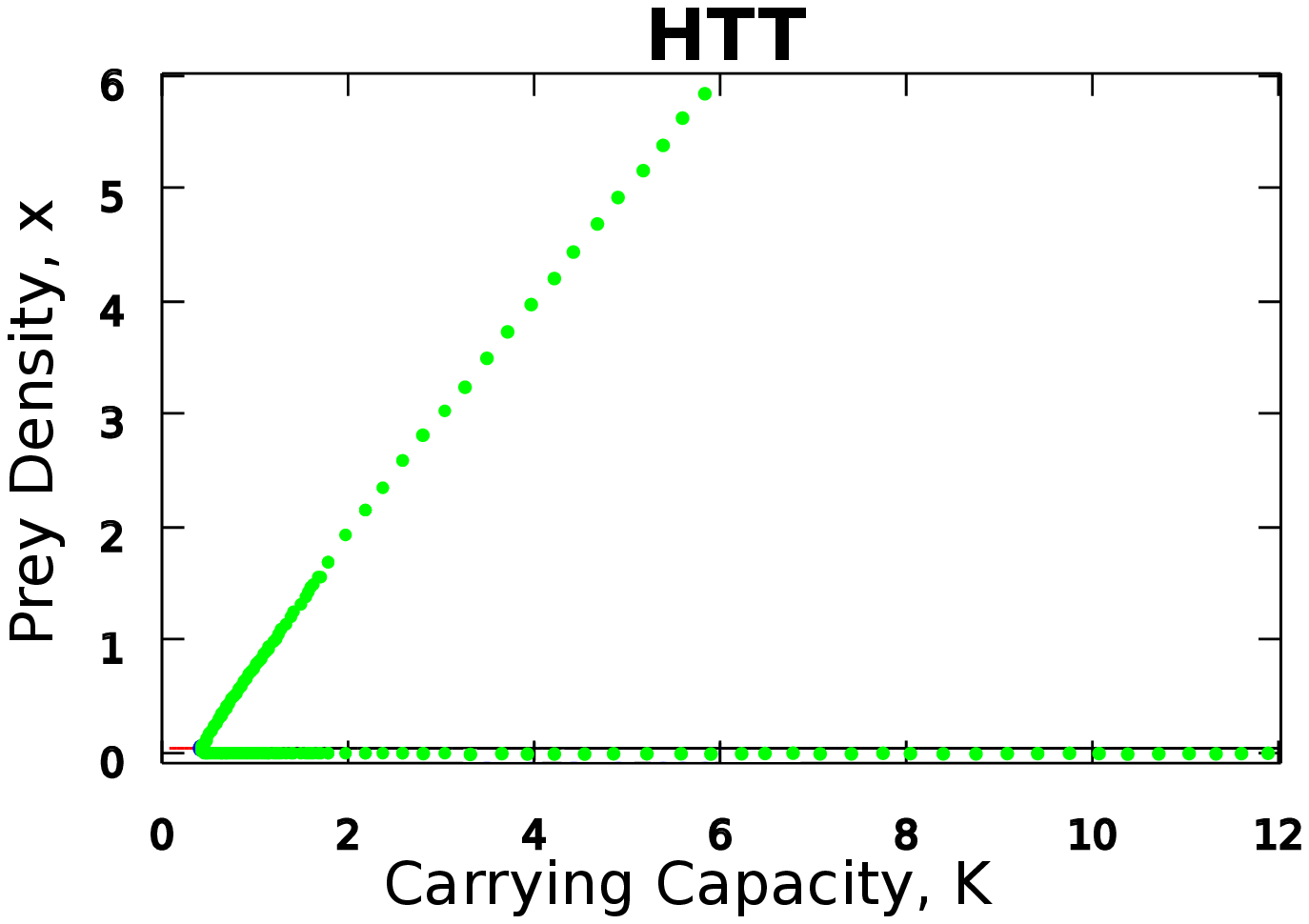}
        \includegraphics[width=0.32\linewidth]{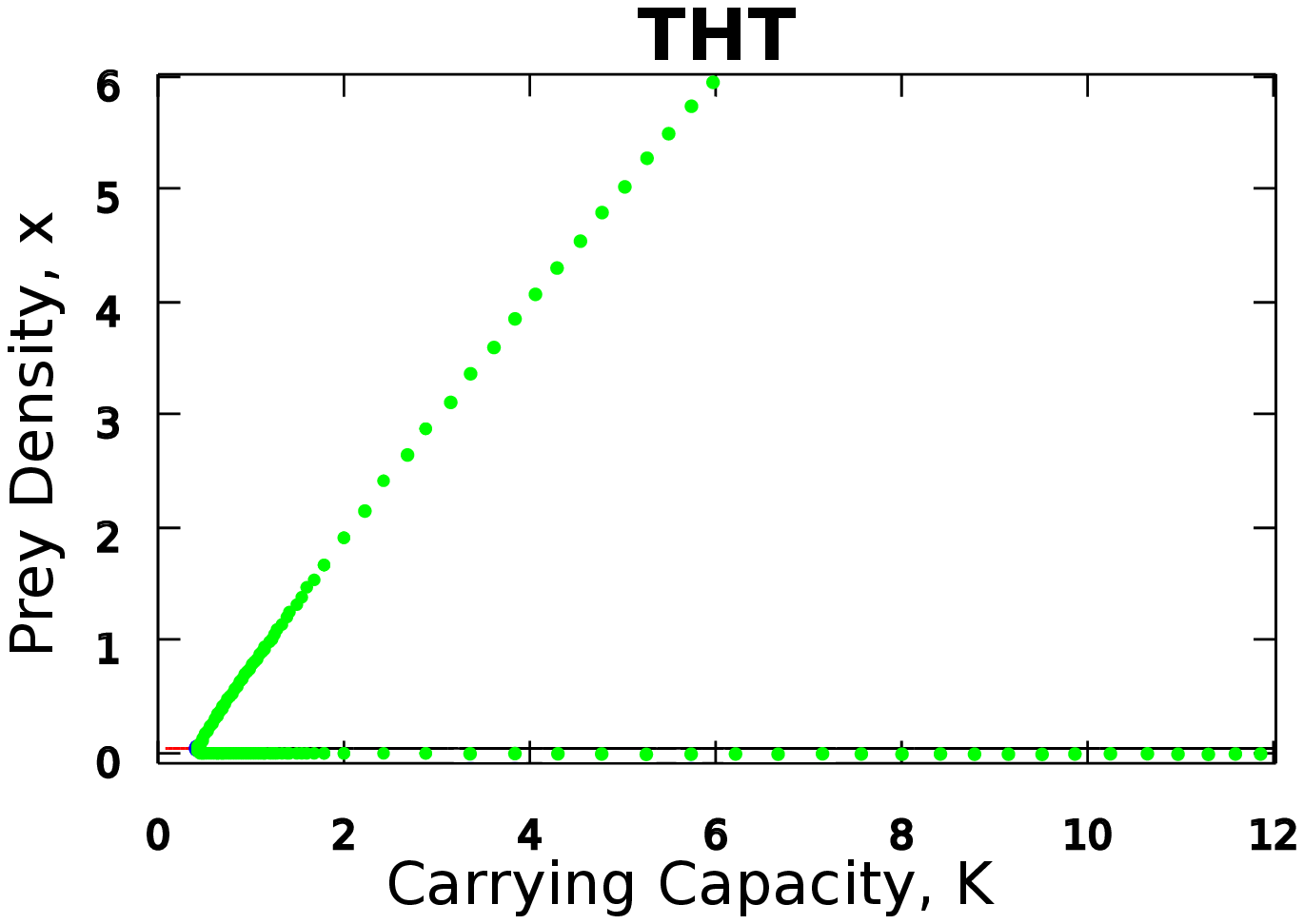}
        \includegraphics[width=0.32\linewidth]{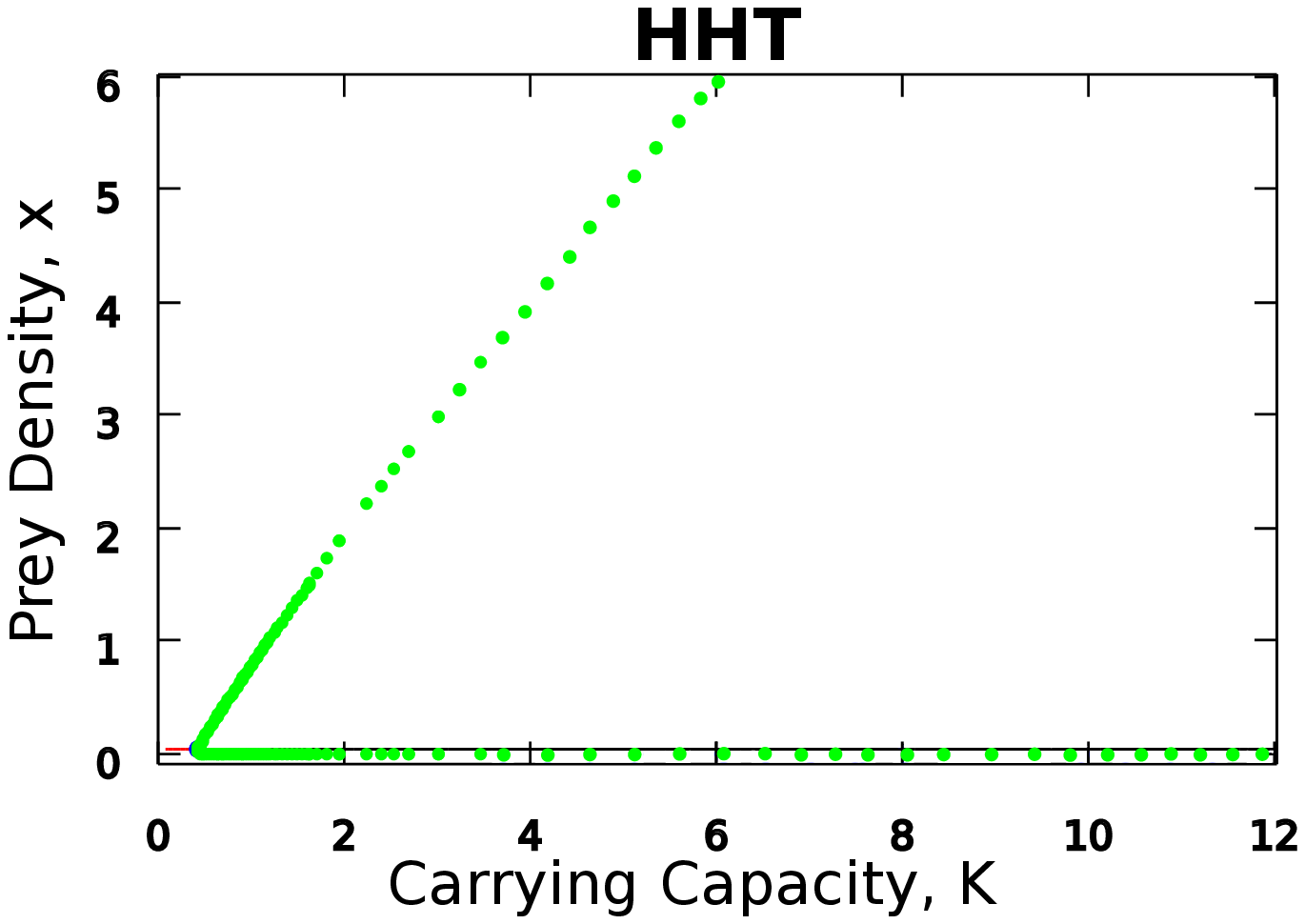}
    \\
        \includegraphics[width=0.32\linewidth]{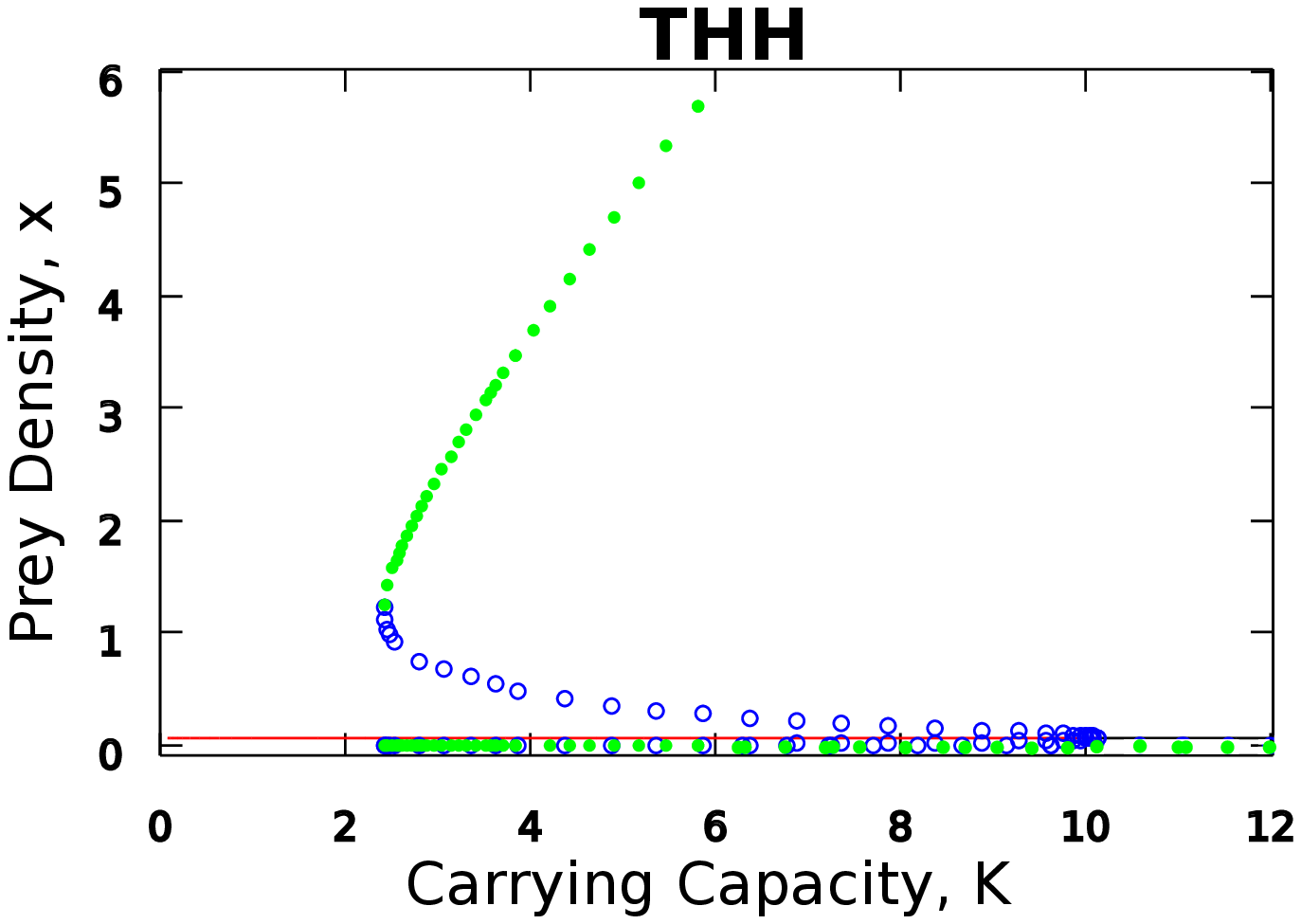}
        \includegraphics[width=0.32\linewidth]{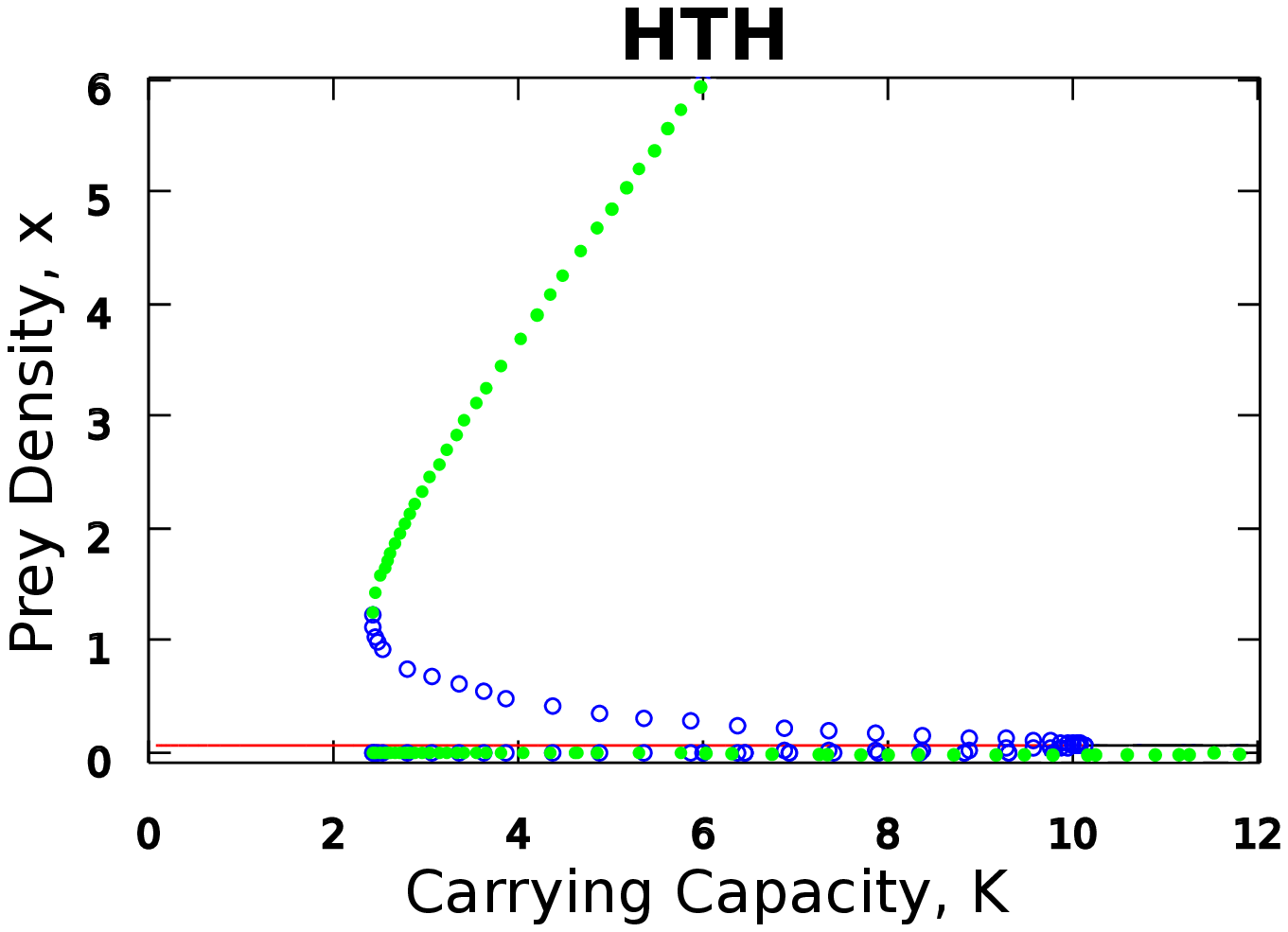}
        \includegraphics[width=0.32\linewidth]{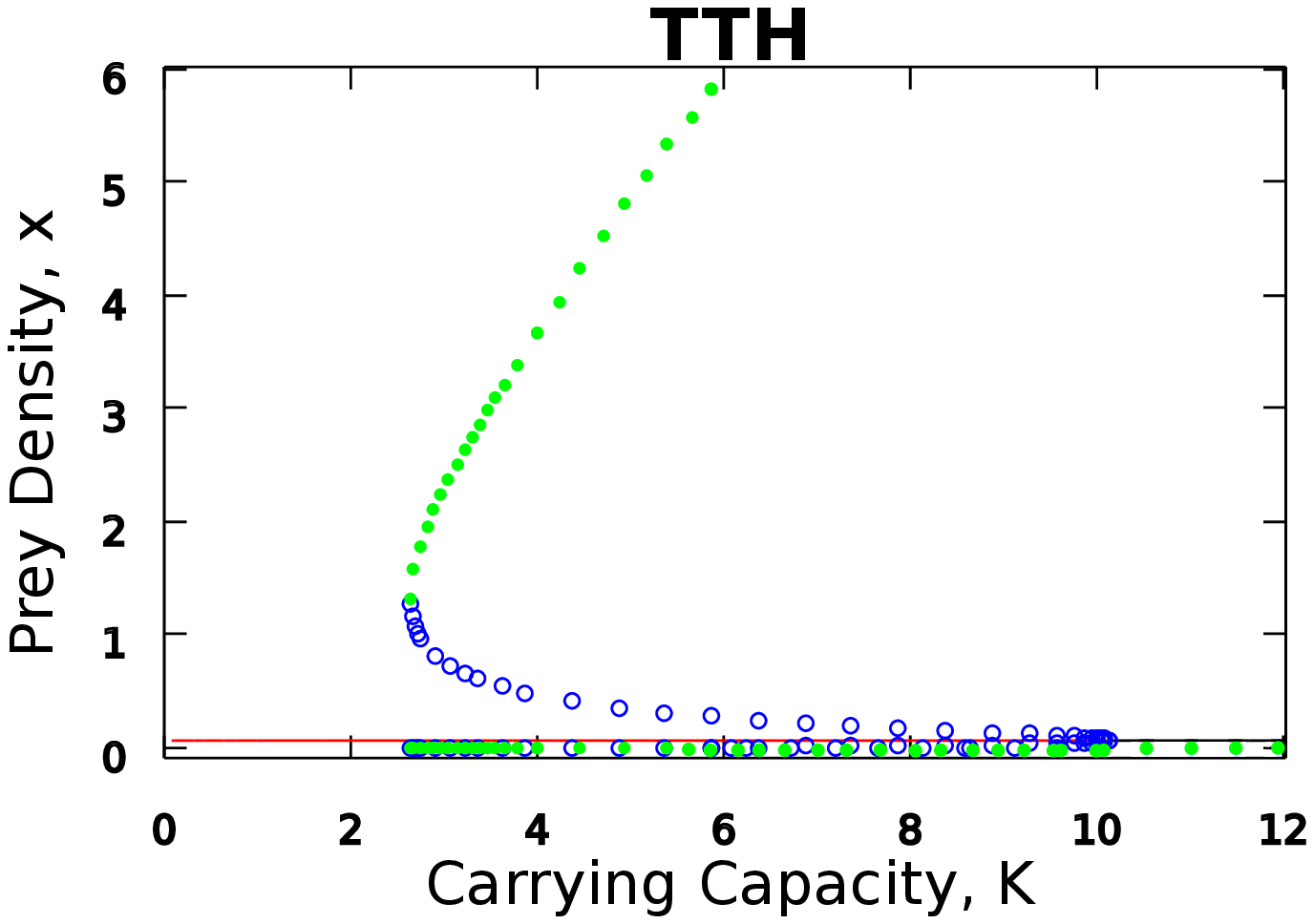}
        \subcaption{Piecewise Identical Functional Response}
        \label{fig:SwitchingBifurRM}
    \end{subfigure}
    \\
    \includegraphics[width=\linewidth]{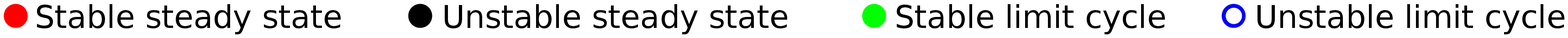}
    \caption{Bifurcation diagrams for the Rosenzweig-MacArthur model \eqref{RM}: a) the original model; b) sub-models using least squares applied to the functional responses over the prey density interval $[0,0.1]$ where the Holling functional response is fixed; and c/d) sub-models using piecewise identical functional responses. The corresponding functional responses are given in Figure S1 and the bifurcation values are given in Table \ref{table:BifurValues}. As described in Section \ref{Sect:StochResults}, the addition of stochasticity causes extinction in the Rosenzweig-MacArthur model and thus, no bifurcation diagrams could be produced.}
\end{figure}

\subsection{Fitted Functional Response Results}
%------------------------------------------------------------------------------------
\label{sect:LSResults}
We only consider the fixed Holling case here as the fixed Ivlev and fixed trigonometric cases are similar and given in the Supplementary Material (Figures S4 and S5).

The Rosenzweig-MacArthur fitted functional response sub-models yield the bifurcation plots shown in Figure~\ref{fig:CoexSSBifurRM}. The bifurcation values for these sub-models are given in Table \ref{table:BifurValues}. These fits result in better similarity of the bifurcation structures and less structural sensitivity as compared to the original model (Figure \ref{fig:BifurRM}). For example, when the Ivlev functional response is fit to the fixed Holling functional response, the Ivlev fit to Holling sub-model shows behaviour similar to the Holling only sub-model as there is a leftward shift in the destabilizing Hopf bifurcation. This result holds across each of the fits. Hence, considering the prey densities in the chosen fitting region appears to be a better metric to describe functional responses in the Rosenzweig-MacArthur model; in the sense that model behaviour is less sensitive to small changes in the functional response using this method.

The Leslie-Gower-May fitted functional response sub-models yield the bifurcation diagrams given in Figure \ref{fig:CoexSSBifurMay} and the value of the bifurcation parameter for each bifurcation is given in Table \ref{table:BifurValues}. The resulting sub-models show better similarity in bifurcation structure when compared to the fixed functional response. For example, Ivlev fit to Holling shows a leftward shift of the Hopf bifurcation as compared to the Ivlev only bifurcation structure and the disappearance of the saddle bifurcation of limit cycles. So, the Ivlev fit to Holling sub-model shows behaviour similar to the Holling only sub-model. The same result holds for each of the other cases, and there is better similarity in model behaviour when compared to the fixed functional response. Therefore, as in the Rosenzweig-MacArthur model, fitting the functional responses over the chosen fitting region appears to decrease structural sensitivity in the Leslie-Gower-May model.

We would like to note that these theoretical results are subjective to which functional response is chosen to be fixed. That is, the model behaviour will always become more similar to the fixed model. In the case of the results shown in this section, the model behaviour always becomes more similar to the Holling sub-model behaviour. Our method when applied to a model that uses data to fit the functional responses, however, is not subjective. We will discuss this application of our method to data further in Section \ref{Sect:Discuss:Data}.

\begin{table}[ht]
\centering
\begin{tabular}{|c|c|c|c|c|c|c|}
    \hline 
   & \multicolumn{3}{|c|}{Rosenzweig-MacArthur} & \multicolumn{3}{|c|}{Leslie-Gower-May} \\
    \hline
    
   \textbf{\large Original Model (\ref{fig:BifurRM}, \ref{fig:BifurMay})} & \textbf{Holling} & \textbf{Ivlev} & \textbf{Trig} & \textbf{Holling} & \textbf{Ivlev} & \textbf{Trig} \\
    \hline
    Hopf & 0.4452 & 1.071 & 10.12 & 5.554 & 12.09 & 26.36 \\
    Limit Cycle Saddle & - & - & 2.644 & - & 10.51 & 12.94\\
    \hline
    
    \textbf{\large Holling Fixed (\ref{fig:CoexSSBifurRM}, \ref{fig:CoexSSBifurMay})} & \textbf{Holling} & \textbf{Ivlev} & \textbf{Trig} & \textbf{Holling} & \textbf{Ivlev} & \textbf{Trig} \\
    \hline
    Hopf & 0.4452 & 0.4632 & 0.7729 & 5.554 & 5.857 & 6.052 \\
    Limit Cycle Saddle & - & - & 0.5057 & - & - & 5.964\\
    \hline

    \textbf{\large Pw-I Holling (\ref{fig:SwitchingBifurRM},  \ref{fig:SwitchingBifurMay})} & \textbf{HTT} & \textbf{THT} & \textbf{HHT} & \textbf{HTT} & \textbf{THT} & \textbf{HHT} \\
    \hline
    Hopf & 0.4451 & 0.4451 & 0.4451 & 26.36 & 26.36 & 5.554 \\
    Limit Cycle Saddle & - & - & - & 7.95 & 7.95 & - \\
    \hline
    
    \textbf{\large Pw-I Trigonometric (\ref{fig:SwitchingBifurRM}, S8)} & \textbf{THH} & \textbf{HTH} & \textbf{TTH} & \textbf{THH} & \textbf{HTH} & \textbf{TTH} \\
    \hline
    Hopf & 10.12 & 10.12 & 10.12 & 5.554 & 5.554 & 26.36 \\
    Limit Cycle Saddle & 2.431 & 2.431 & 2.644 & - & - & 13.13 \\
    \hline
\end{tabular}
    \caption{Values of the bifurcation parameter, the carrying capacity, for each bifurcation shown in the respective figures. The values in parentheses give the figure number of the corresponding bifurcation diagrams. Pw-I stands for piecewise identical, H stands for Holling, and T stands for trigonometric.}
    \label{table:BifurValues}
\end{table}

\subsection{Piecewise Identical Functional Response Results}
%------------------------------------------------------------------------------------
\label{Sect:SensRegionsResults}

Suppose that medium prey densities determine the bifurcation structure of a model. Then, according to our hypothesis, we expect YXY and ZXZ to show model behaviours more similar to the X model behaviour (be that H, I, or T). Based on the results in Section \ref{sect:LSResults}, we might actually expect low prey densities to be the most important, in which case we expect XYY and XZZ to have better similarity to the X model behaviour. However, our results show that none of the three regions we consider separately determine the model behaviour. In this section, we describe only the piecewise identical functional responses made up of the Holling and trigonometric functional responses. The bifurcation diagrams for these cases of the Rosenzweig-MacArthur and Leslie-Gower-May sub-models are given in Figures \ref{fig:SwitchingBifurRM} and \ref{fig:SwitchingBifurMay}. The rest of the cases are similar and are given in the Supplementary Material (Figures S6-S8).

In the case of the Rosenzweig-MacArthur model, there is not better similarity of the Holling to trigonometric sub-models (i.e., sub-models of the form HTT, THT, or HHT) behaviour towards the trigonometric behaviour. This result suggests that both the first and second regions determine the resulting bifurcation behaviour of the Rosenzweig-MacArthur model and that these regions are the most structurally sensitive.

However, the Rosenzweig-MacArthur trigonometric to Holling piecewise identical sub-models (i.e., sub-models of the form THH, HTH, or TTH) are slightly different. In the THH plot, there is a small leftward shift of the saddle bifurcation of limit cycles from the original value (Table \ref{table:BifurValues}) at $K=2.644$ to $K=2.431$. There is, however, no change in the location of the Hopf bifurcation. In the HTH plot, the same leftward shift occurs. Lastly, the TTH plot has a resulting bifurcation structure that matches that of the trigonometric only functional response model behaviour. Overall, the trigonometric functional response has a strong effect on the bifurcation structure. Further, when the Holling functional response takes up two of the three regions (in the THH and HTH cases), the Holling functional response has an effect, though limited, on the bifurcation structure. This effect causes the aforementioned shift in the saddle bifurcation of limit cycles. The piecewise identical functional responses seem to affect the resulting model behaviour without following any clear rules. Hence, we are unable to make any claims on which prey density regions determine the Rosenzweig-MacArthur model behaviour.

In the case of the Leslie-Gower-May model, the HTT and THT sub-models show better similarity when compared to the trigonometric bifurcation structure while the resulting bifurcation structure for the HHT sub-model matches that of the Holling only functional response. We found that in the Leslie-Gower-May model, the functional response that is used across more prey density regions is the one that determines where the Hopf bifurcation occurs and whether there exists a saddle bifurcation of limit cycles. However, the functional response used across only one prey density region causes a shift in the saddle bifurcation of limit cycles (if it exists) for every case, except surprisingly the IIH sub-model where there is no shift from the Ivlev behaviour (Figure S8). This result means that the sub-models with the Holling functional response describing two prey density regions result in only Holling model behaviour, but the sub-models with either the Ivlev or trigonometric functional response describing two prey density regions result in mixed dynamics.

\section{Discussion}
%------------------------------------------------------------------------------------
\label{Discussion}
\subsection{Reducing Structural Sensitivity}
%------------------------------------------------------------------------------------
Structural sensitivity is an intrinsic feature of a wide range of ecological models \citep{adamson:2013, adamson:2014}. In our work, we focus on predator-prey models and investigate the structural sensitivity that occurs with respect to the predation function. Our work has attempted to answer how to fit data in such a way to gain confidence in the results, regardless of which curve is used. We found that only the sub-models that use the functional responses fit around the coexistence steady state consistently decrease structural sensitivity across both the Rosenzweig-MacArthur and Leslie-Gower-May models. Using this result, we also give suggestions as to which intervals of prey density are most relevant when fitting a functional response to data and which intervals of prey density are most relevant for ecologists to collect data.

\subsubsection{Stochasticity}
%------------------------------------------------------------------------------------
\label{Sect:Discuss:Stoch}
One of the differences between the functional response curves in \citet{fussmann:2005} and real data is that the latter include stochasticity. A model's steady states and deterministic structure only provide a partial understanding of the stochastic structure and so the stochastic system may behave in unanticipated ways \citep{abbott:2015}. Hence, it is important to consider a stochastic model when analyzing biological systems. We examine the effect of adding stochasticity to the functional responses of the Rosenzweig-MacArthur and Leslie-Gower-May models and, specifically, whether or not this addition can remove or decrease structural sensitivity.

In our work, we observe that the destabilizing of the coexistence steady state can occur earlier in the stochastic system than the deterministic system. Most notably, there is greater similarity in the stochastic bifurcation structures from the three functional responses implying decreased structural sensitivity. We find that the maximum and minimum values of the deterministic system lie within the interval determined by the maximum and minimum values of the stochastic system. This result supports the observation that trajectories in the stochastic system may approach the steady states of the deterministic system \citep{aguirre:2013, predescu:2007}. How stochasticity alters a model's behaviour from the deterministic system depends on the amplitude and distribution of the stochastic random variable. For example, white noise weakens bistable behaviour while red noise can amplify the bistability observed in a system \citep{abbott:2015}. Using a different distribution can change the shape and depth of potential wells, and thus produce different stability properties \citep{abbott:2015}. We do not, however, study other stochastic distributions here and leave this investigation for future work.

%In our work, we observe that the bifurcation structure of the deterministic system is robust to stochasticity. This result supports the observation that trajectories in the stochastic system approach the steady states  of the deterministic system \citep{aguirre:2013, predescu:2007}. How stochasticity alters a model's behaviour from the deterministic system depends on the amplitude and distribution of the stochastic random variable. For example, white noise weakens bistable behaviour while red noise can amplify the bistability observed in a system \citep{abbott:2015}. Using a different distribution can change the shape and depth of potential wells, and thus produce different stability properties \citep{abbott:2015}. We do not, however, observe any substantial changes to the Leslie-Gower-May model's stability properties for any of the distributions we employ. Altogether, adding stochasticity to the functional responses does not decrease structural sensitivity significantly enough to be considered a better way to describe the functional responses in either of the Rosenzweig-MacArthur or Leslie-Gower-May models.

White noise is widely used in studies of stochastic models \citep{vanDerBolt:2018}. It is characterized by having the same variance at all frequencies. However, in recent years, there has been interest in varying the colour, or auto-correlation, of applied noise. Red noise, which has more auto-correlation and is characterized by low frequencies, has been found to better describe the stochasticity found in ecological systems \citep{vasseur:2004, greenman:2003, grove:2020, jonsson:2019}. In this paper, we only consider white noise. An interesting future step is to consider the effect of red noise on the bifurcation structure of the stochastic system.

\subsubsection{Fitted Functional Responses}
%------------------------------------------------------------------------------------
\label{Sect:Discuss:Metrics}
It is sometimes possible to improve mathematical model predictions through describing data in a particular way. In order to help ecologists make accurate decisions regarding conservation efforts, it is important to obtain the best model predictions possible. We discuss the use of data in these model predictions further in Section \ref{Sect:Discuss:Data} below.

In our models, the slopes of the functional response curves are greatest near the coexistence steady state, meaning that this region is also where the predator response to a change in prey density is strongest \citep{aldebert:2016}. This property suggests that having a higher density of data points in the region of the functional response domain near the coexistence steady state prey density value, or giving these data points more weight, is important. Further, since the initial behaviour of all three models is similar and all are stable at the coexistence steady state before any destabilizing occurs, we hypothesized that if the corresponding part of the functional response (region around the stable or unstable coexistence steady state prey density value) held more weight in the nonlinear least squares fitting, then the bifurcation structures would be more similar. That is, if the fitted functional responses method (Section \ref{Sect:CoexSSMethods}) is a better metric to decrease structural sensitivity, we expect to observe better similarity of the resulting model's bifurcation structure. We found that this hypothesis is supported in both the Rosenzweig-MacArthur and the Leslie-Gower-May models.

\subsubsection{Piecewise Identical Functional Responses}
%-----------------------------------------------------------------------
We found that more similar bifurcation structures are obtained for functional responses that fit very closely at low and intermediate prey density values, but also disagree wildly at large prey density values near saturation. So, the piecewise identical functional response results suggest that low and medium prey densities are the most important regions to fit to decrease structural sensitivity.

Our results, however, appear to be in contrast from previous findings, indicating that it is most important to have data at saturation or high prey densities \citep{aldebert:2018-2}. However, this contradiction is only apparent. In our work, we attempt to decrease structural sensitivity and describe the functional responses such that any of the three can be used to provide equally accurate predictions. On the other hand, \citet{aldebert:2018-2} aim to choose one functional response that fits the model best. To do so, they suggest collecting data around saturation such that the functional responses result in models with vastly different bifurcation structures. Since structural sensitivity in the resulting models is increased, it becomes easier to determine which functional response and model fit the data best. While the authors are able to choose a best functional response to fit data, the model remains structurally sensitive.

\subsection{Using Data in Models}
%------------------------------------------------------------------------------------
\label{Sect:Discuss:Data}

The main application of our work is in models that use data to fit the functional response. As we noted in Section \ref{sect:LSResults}, the choice of the fitted functional response is subjective if our method is used merely theoretically. However, in the case of applied models, our method is no longer subjective since we fit the functional responses to data instead of a subjectively chosen fixed functional response. Once all three functional responses are fit using our method to make their bifurcations structures more similar, the resulting bifurcation structure should show the behaviour of the data rather than a single chosen functional response.

Ecological data, however, is typically sparse as it is often difficult to collect \citep{SparseData}. Furthermore, the data is also typically very noisy. As a result, it is important that ecologists know what data is most important to collect to build accurate models that can be used with confidence to predict future population dynamics. In Section \ref{sect:LSResults}, we found that model predictions are most consistent if the functional responses match over the coexistence steady state prey density region. The implication is that better predictions can be obtained if there is additional data near the coexistence steady state. We believe that this region of prey densities is most important in fitting the functional responses because it corresponds to the steady state behaviour for values of the carrying capacity less than destabilizing bifurcation. The idea being that if the model behaviour is more similar before the destabilizing bifurcation, then the destabilizing bifurcation may occur for more similar values of the carrying capacity. This result opposes previous work that states saturation is the most important part to match of the functional responses \citep{aldebert:2018-2}.

To apply our result to ecological data and to make predictions in predator-prey models, we propose the following procedure:
\begin{enumerate}
    \item Collect data as per normal methods,
    \item Fit the three functional responses to this data and determine the parameter values for each $a_i$ and $b_i$,
    \item Use the functional responses to choose a fitting region around the stable or unstable coexistence steady state prey density values across all three models,
    \item Do one or both of the following:
    \begin{enumerate}
        \item Collect more data for prey densities in the coexistence steady state prey density region,
        \item Give more weight to the data in the coexistence steady state prey density region,
    \end{enumerate}
    \item Fit the three functional responses again and use these functional responses to make predictions regarding the system.
\end{enumerate}
This algorithm takes into account our result that collecting data for prey densities around the coexistence steady state appears to be most important while also considering the difficulty of collecting data (step 4b). Note that we have only considered structural sensitivity with respect to the functional response in this work. Models that are structurally sensitive with respect to other terms, such as prey growth, may require a different algorithm.

\section{Conclusion}
%------------------------------------------------------------------------------------
\label{Conclusion}
Structural sensitivity in the functional response of a model can cause similar functional responses to produce very different bifurcation structures. In this paper, we investigated three potential methods to achieve better similarity of resulting bifurcation structures. We found that 1) stochasticity results in better similarity of the bifurcation structures, 2) functional responses fit over the coexistence steady state prey density region seem to decrease structural sensitivity, and 3) piecewise identical functional responses do not consistently decrease structural sensitivity. Overall, adding stochasticity and fitting functional responses over the coexistence steady state prey density region appear to reduce structural sensitivity. Further investigation into the fitted functional response method is necessary to determine if there is an optimal interval size that decreases structural sensitivity the most. 

Investigating other types of stochasticity, specifically red noise or noise with stronger auto-correlation, remains an interesting question for further investigation. Working with ecological data to test our proposed procedure and make predictions regarding the population dynamics of a predator-prey system is a very useful application of this work. It is also of interest to investigate more predator-prey models and more functional response curves to determine if our results generalize. It is also possible to have structural sensitivity occur with respect to the prey growth function \citep{adamson:2013}. This type of structural sensitivity remains an area for future investigation.

As part of this work, we also investigated the derivatives of the functional responses. That is, we fit the parameters $a_i$ and $b_i$ such that the derivatives of the functional responses were as similar as possible. This method of fitting did not show promise and requires further investigation to determine if it can decrease structural sensitivity.

\bibliographystyle{spmpscinat.bst}
\bibliography{Sources.bib}

\end{document}

\maketitle
\linenumbers

\section{Rosenzweig-MacArthur Functional Responses}
\begin{figure}[H]
    \centering
    \begin{subfigure}{0.4\linewidth}
        \includegraphics[width=\linewidth]{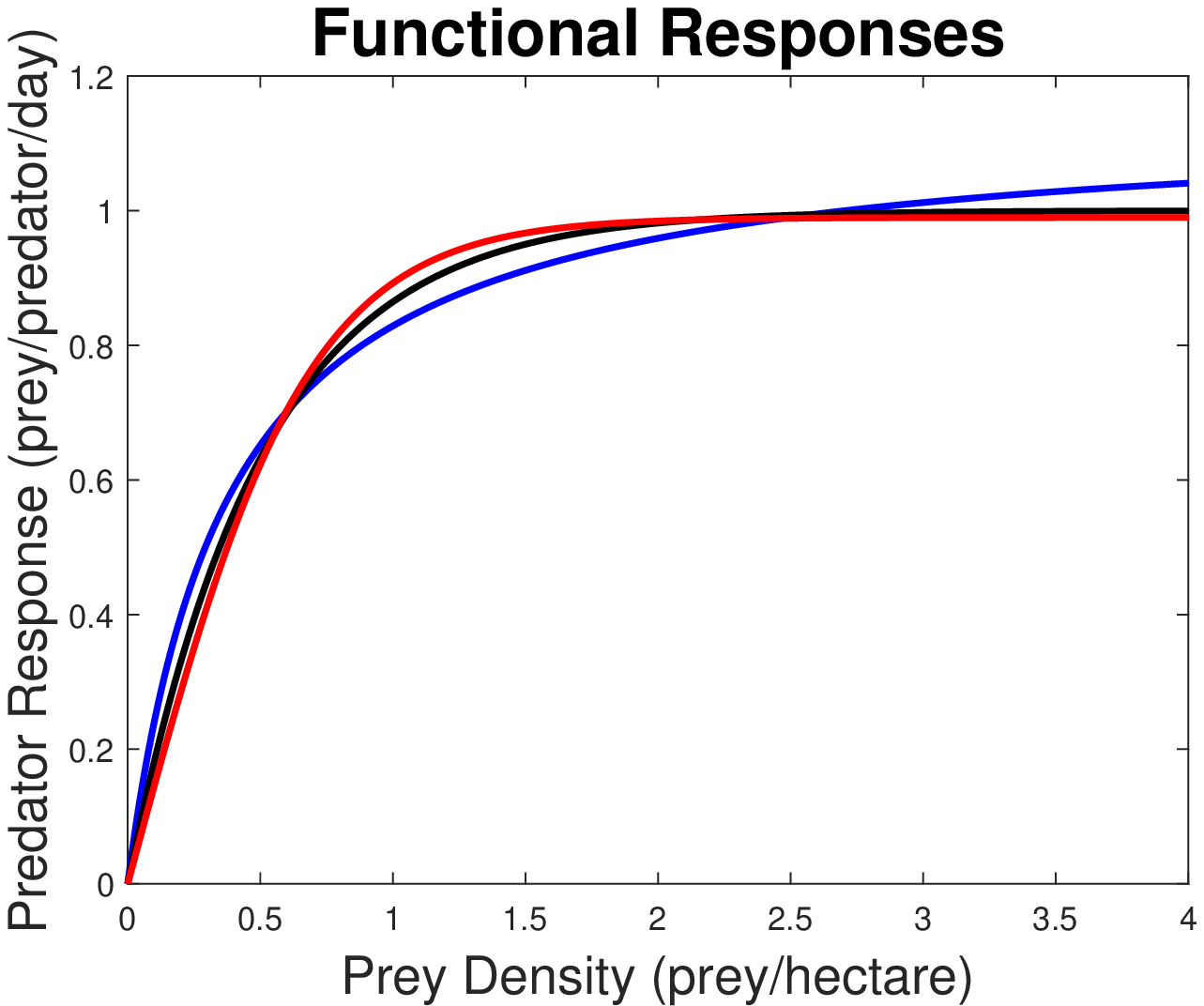}
        \subcaption{Original Functional Responses}
        \label{supfig:CurvesOrigRM}
    \end{subfigure}
    \begin{subfigure}{0.4\linewidth}
        \includegraphics[width=\linewidth]{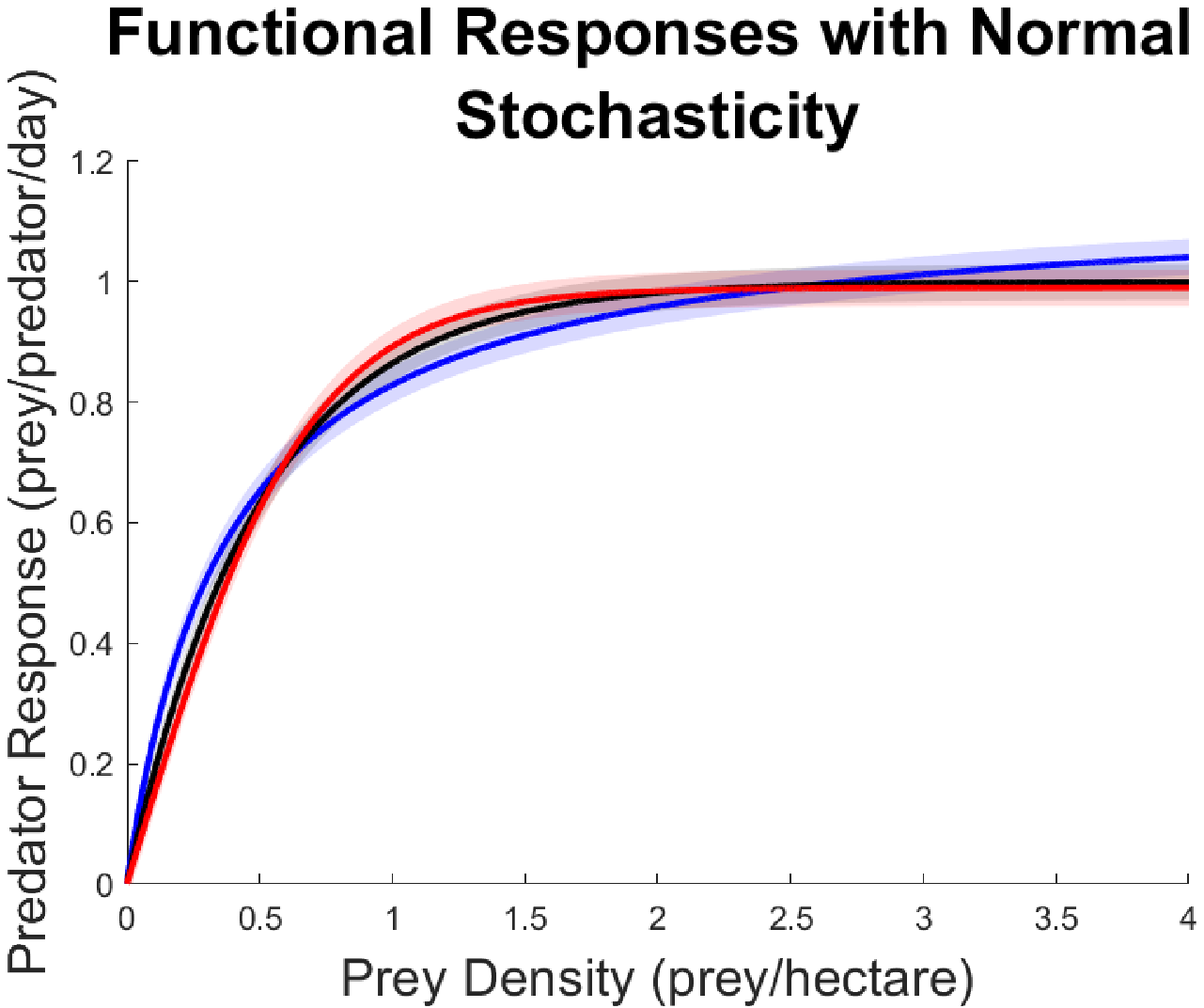}
        \subcaption{Stochastic Functional Responses}
        \label{supfig:StochCurvesRM}
    \end{subfigure}
    \\
    \begin{subfigure}{0.4\linewidth}
        \includegraphics[width=\linewidth]{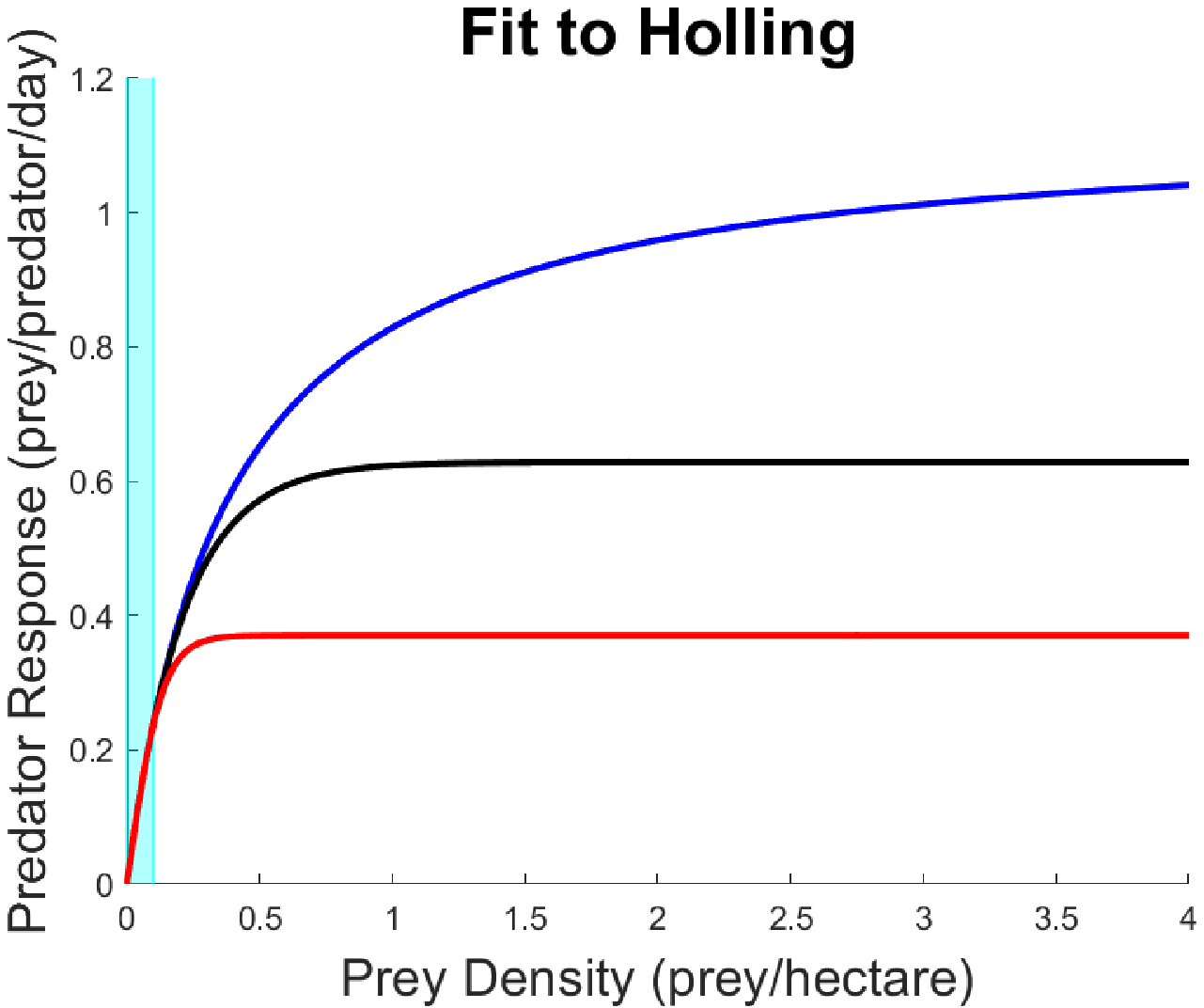}
        \subcaption{Fitted Functional Responses}
        \label{supfig:SSCurvesRM}
    \end{subfigure}
    \begin{subfigure}{0.4\linewidth}
        \includegraphics[width=\linewidth]{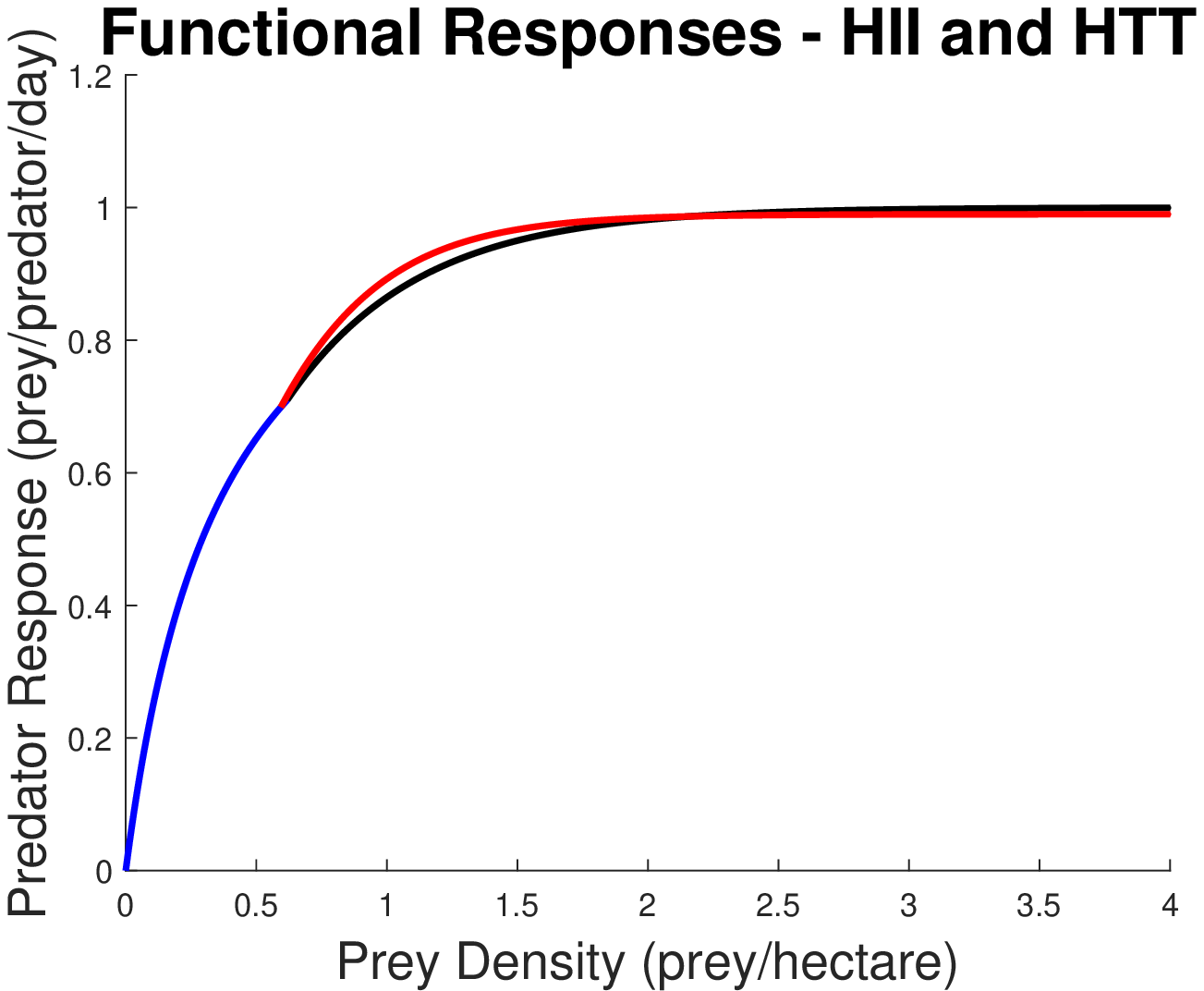}
        \subcaption{Piecewise Identical Functional Responses}
        \label{supfig:SwitchingCurvesRM}
    \end{subfigure}
    \includegraphics[width=0.7\linewidth]{Pictures/FIGLeg1.eps}
    \caption{Functional responses used in the Rosenzweig-MacArthur models: a) the original model; b) stochastic functional responses with $\sigma = 0.03$ where the shaded region represents one standard deviation; c) the functional responses with nonlinear least squares applied over prey densities in [0,0.1]; and d) sample piecewise identical functional responses. Parameter values for all of the models and intersection values for the piecewise identical responses are given in the manuscript and Table \ref{suptable:ParametersLS}.}
    \label{supfig:RMCurves}
\end{figure}

%-------------------------------------------------------------
In Figure \ref{supfig:RMCurves}, we give the functional responses for the original Rosenzweig-MacArthur model and the three methods to decrease structural sensitivity that we consider. Note that the stochastic functional response has been changed to $\sigma = 0.03$ since $\sigma = 0.01$ as used in the manuscript is not visible in the figure.

\section{Stochasticity: sample time series other amplitudes}
\begin{figure}[H]
    \centering
    \includegraphics[width=0.325\linewidth]{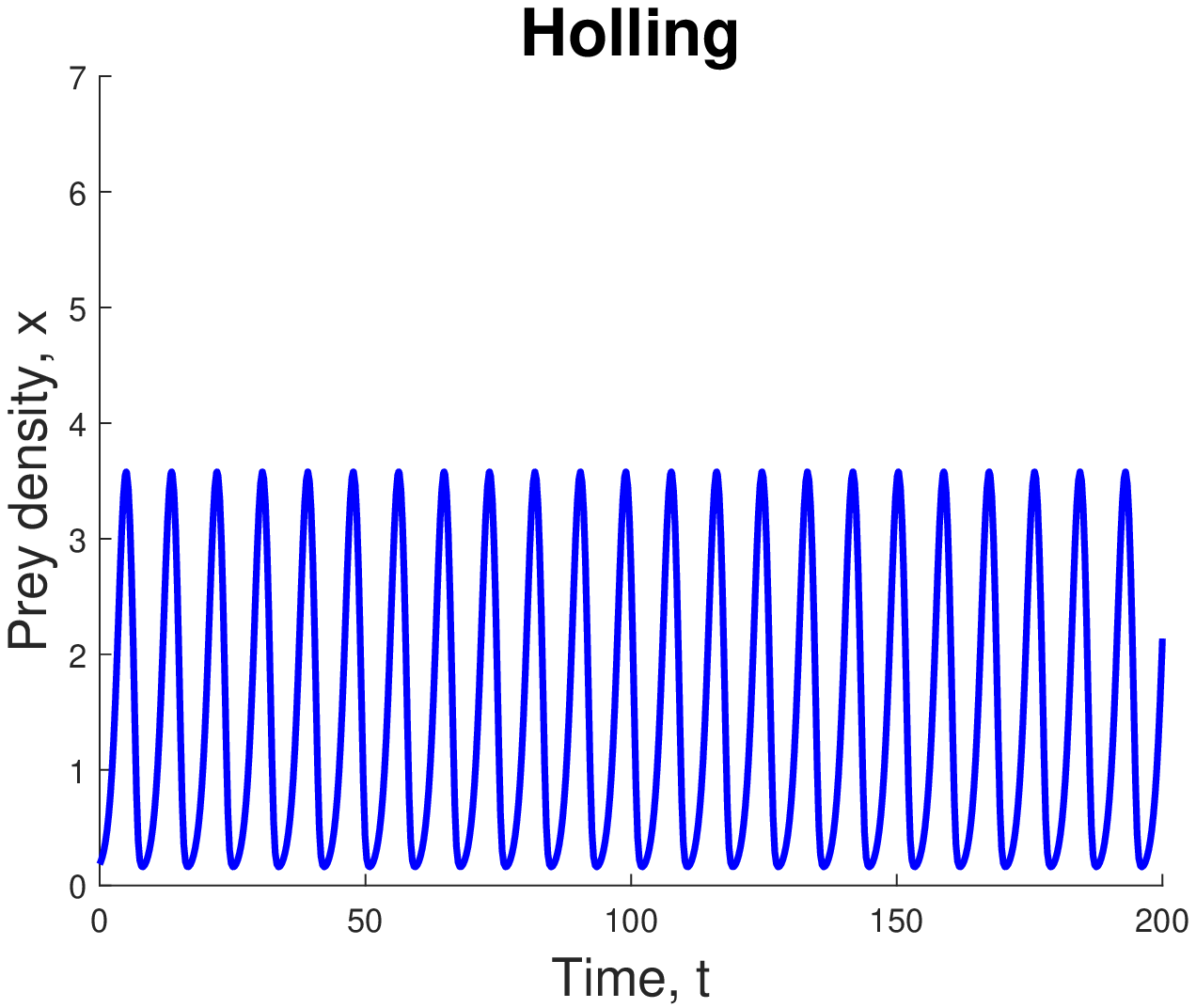}
    \includegraphics[width=0.325\linewidth]{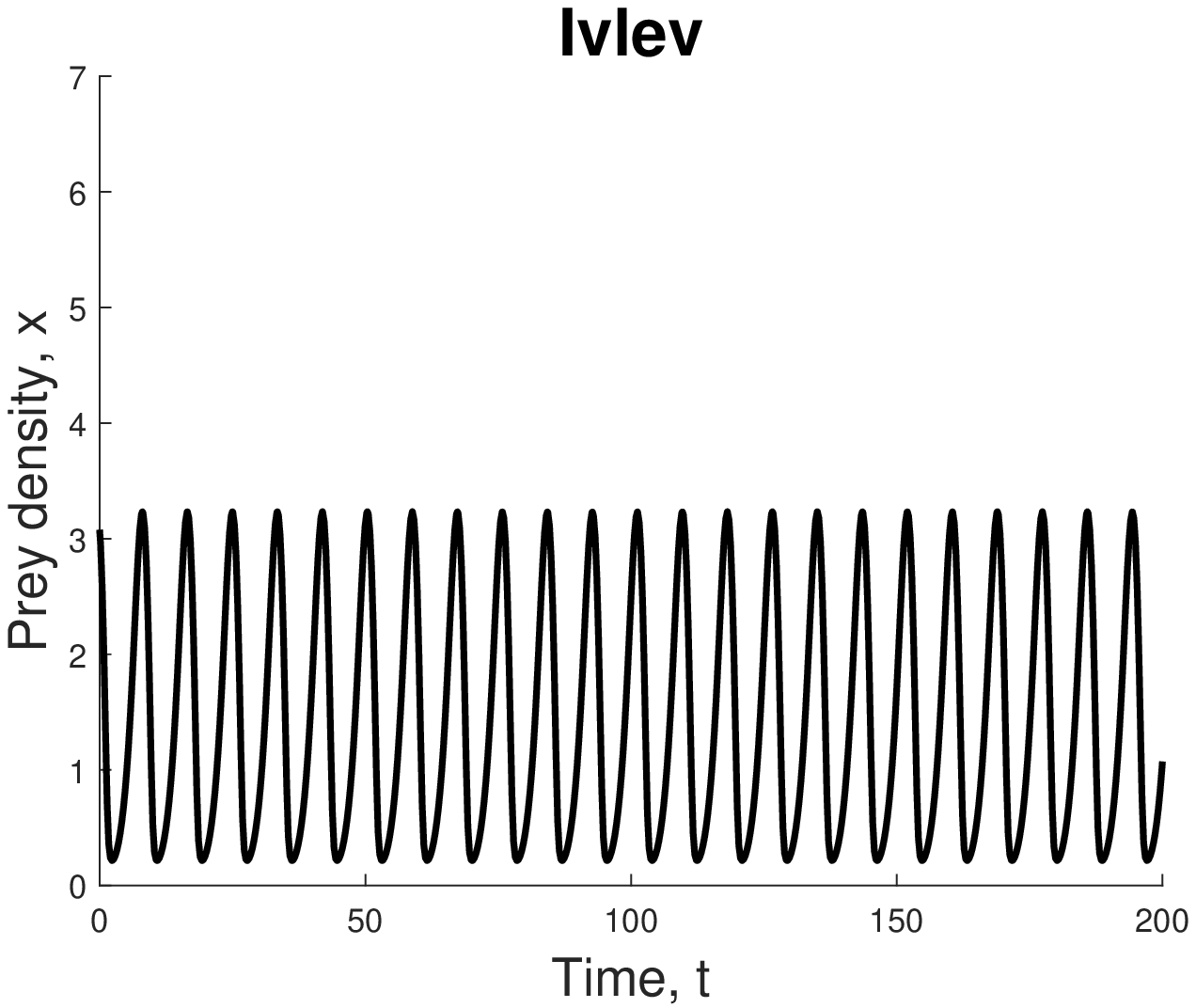}
    \includegraphics[width=0.325\linewidth]{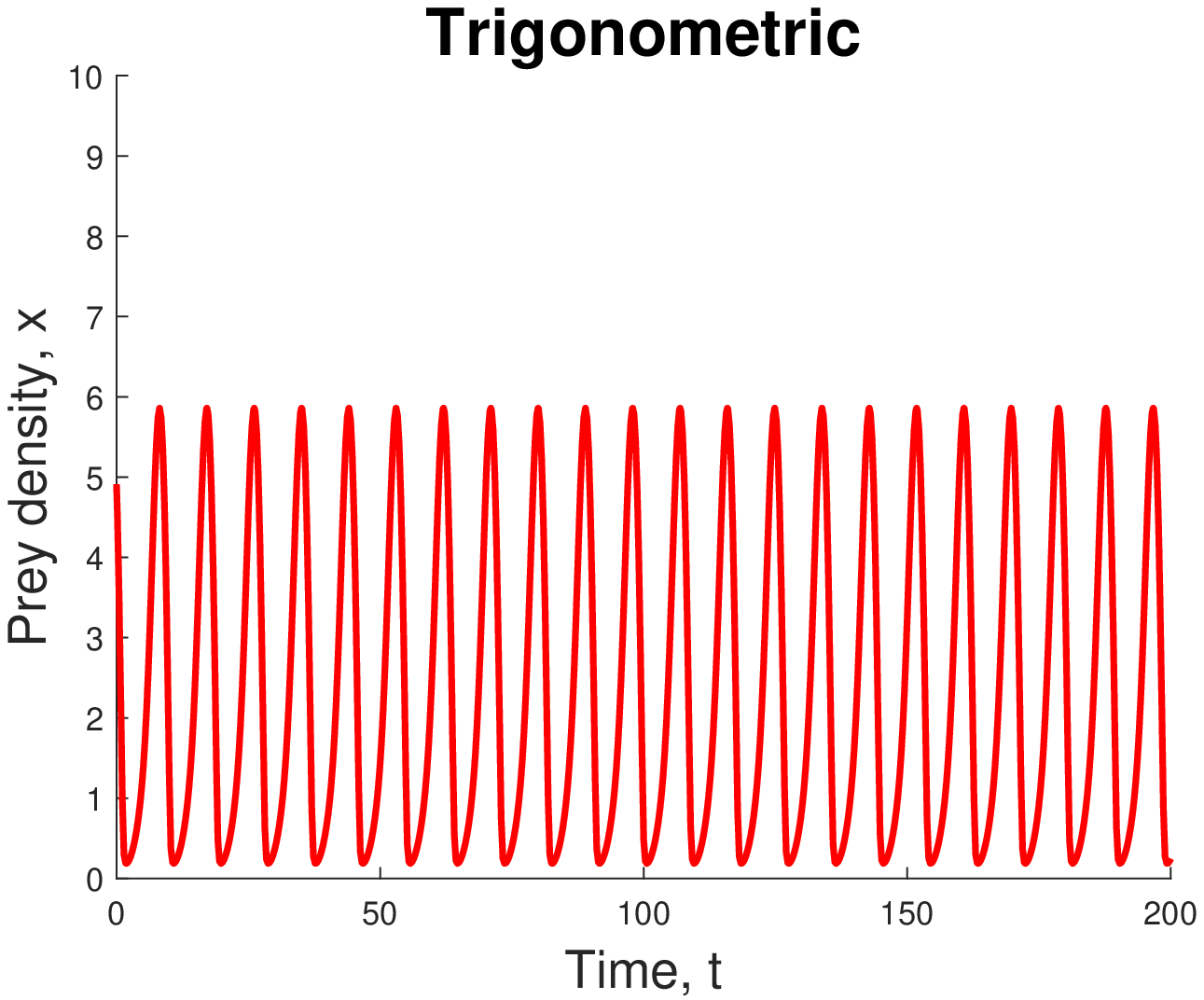}
    \\
    \includegraphics[width=0.325\linewidth]{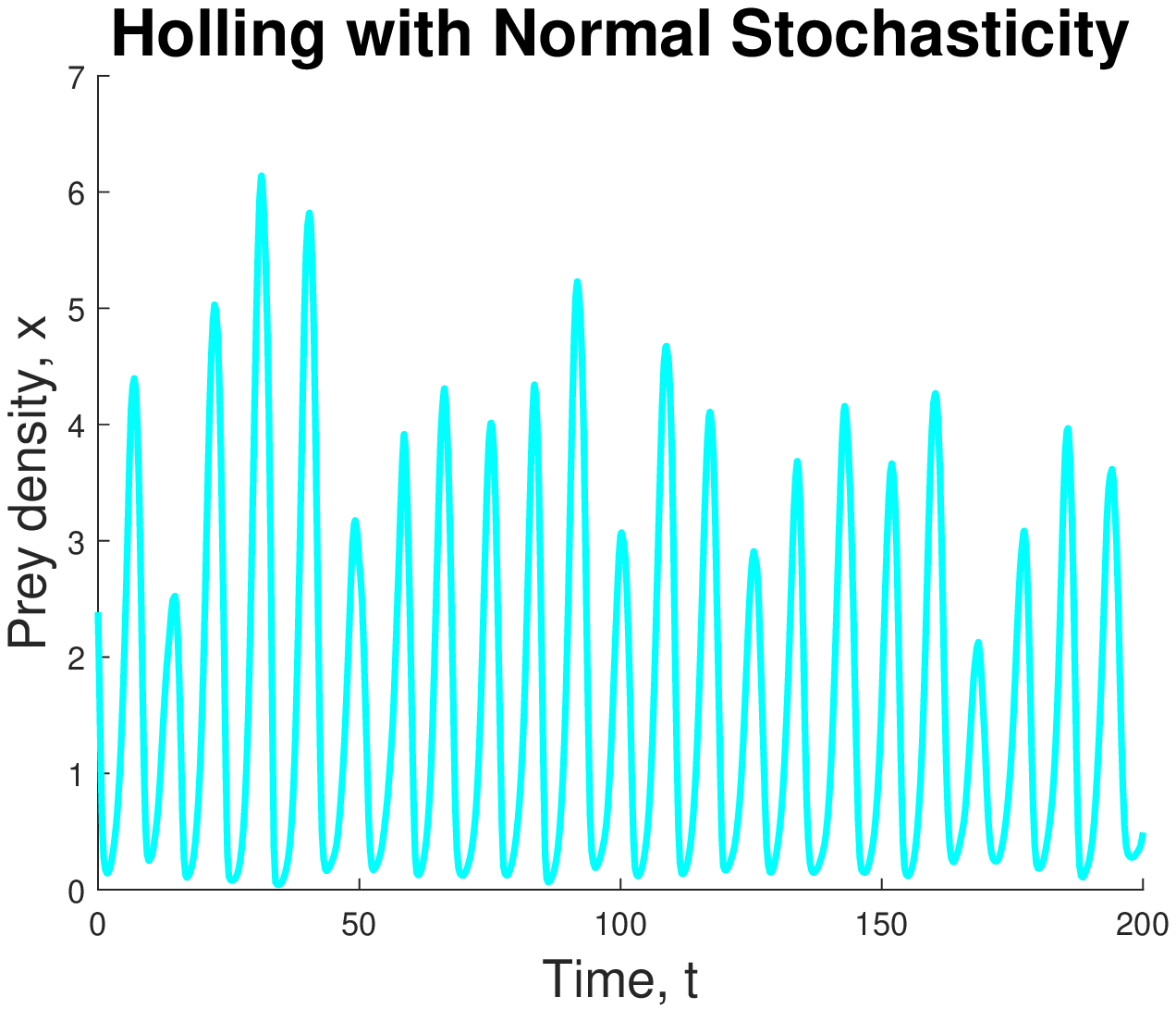}
    \includegraphics[width=0.325\linewidth]{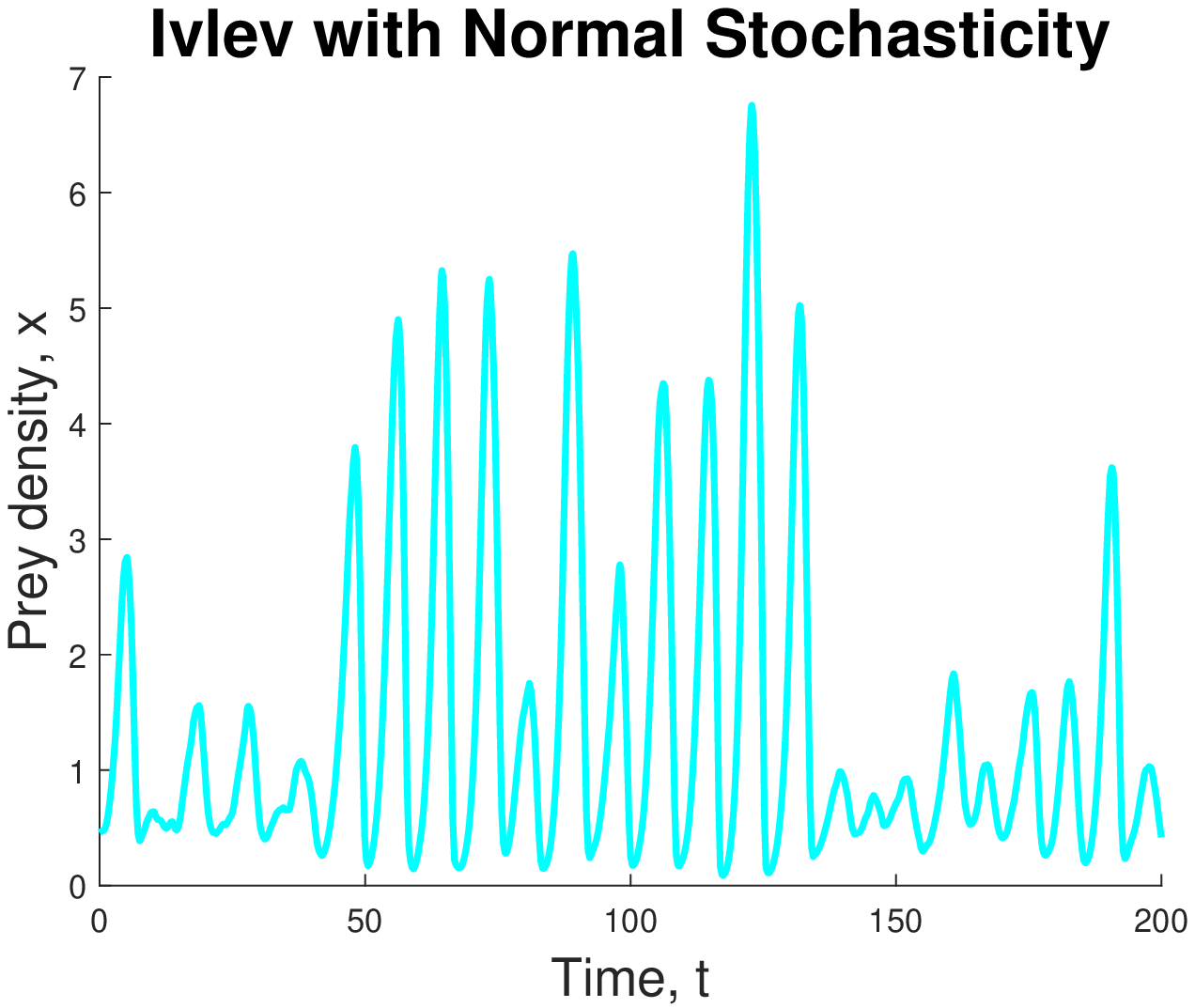}
    \includegraphics[width=0.325\linewidth]{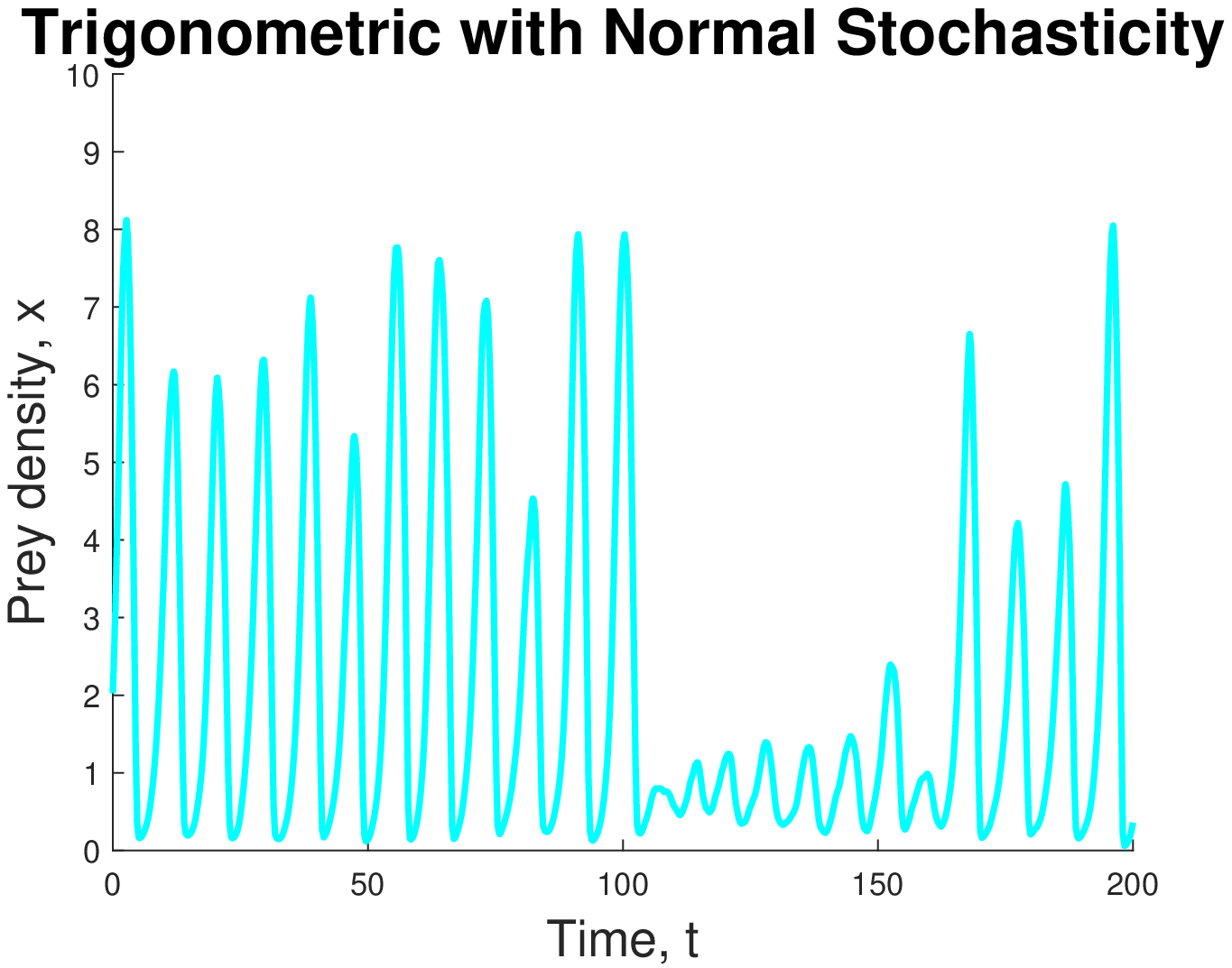}
    \\
    \includegraphics[width=\linewidth]{Pictures/FIGLeg3.eps}
    \caption{Sample time series of the Leslie-Gower-May models. Top: Deterministic limit cycles. Bottom: Stochastic behaviour with $\sigma = 50$. The time series are run using carrying capacity $K=9$ (Holling), $K=11$ (Ivlev), and $K=15$ (trigonometric). Note that the axes vary across the models.}
    \label{supfig:Bistability}
\end{figure}

%-------------------------------------------------------------
In Figure \ref{supfig:Bistability}, we provide sample times series for the Leslie-Gower-May models. The Ivlev and trigonometric time series show that bistability remains a property of the stochastic models.

\begin{figure}[H]
    \centering
    \includegraphics[width=0.325\linewidth]{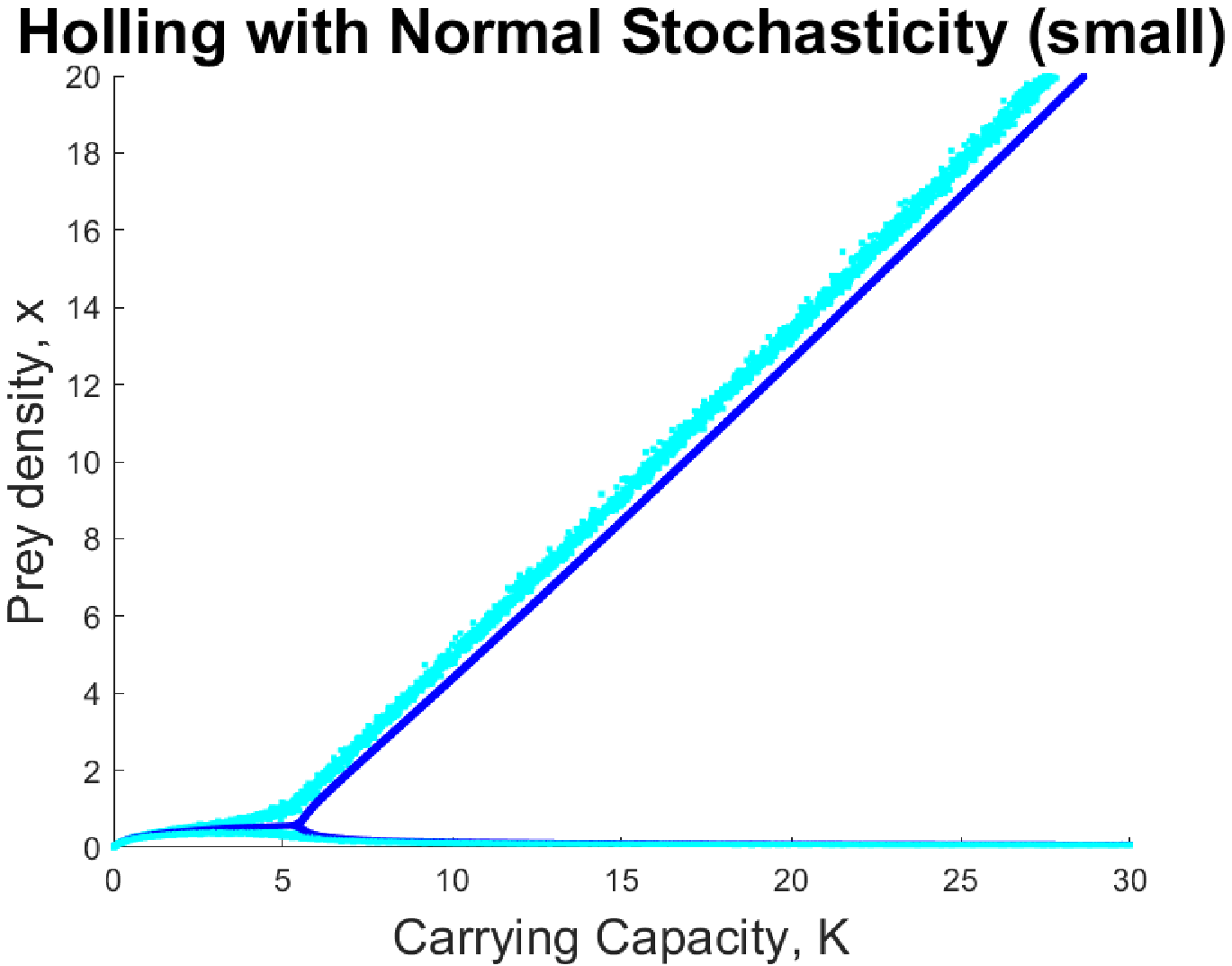}
    \includegraphics[width=0.325\linewidth]{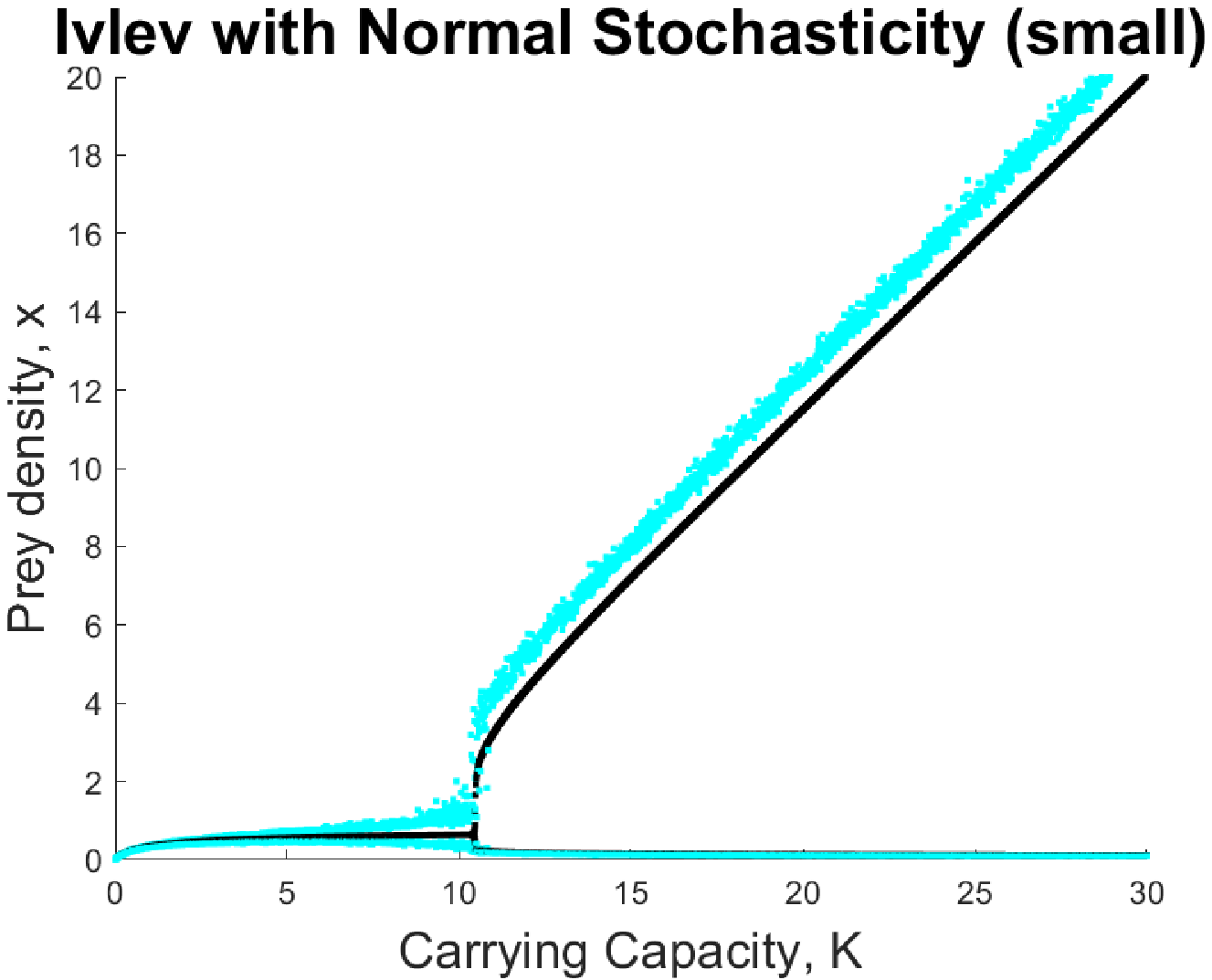}
    \includegraphics[width=0.325\linewidth]{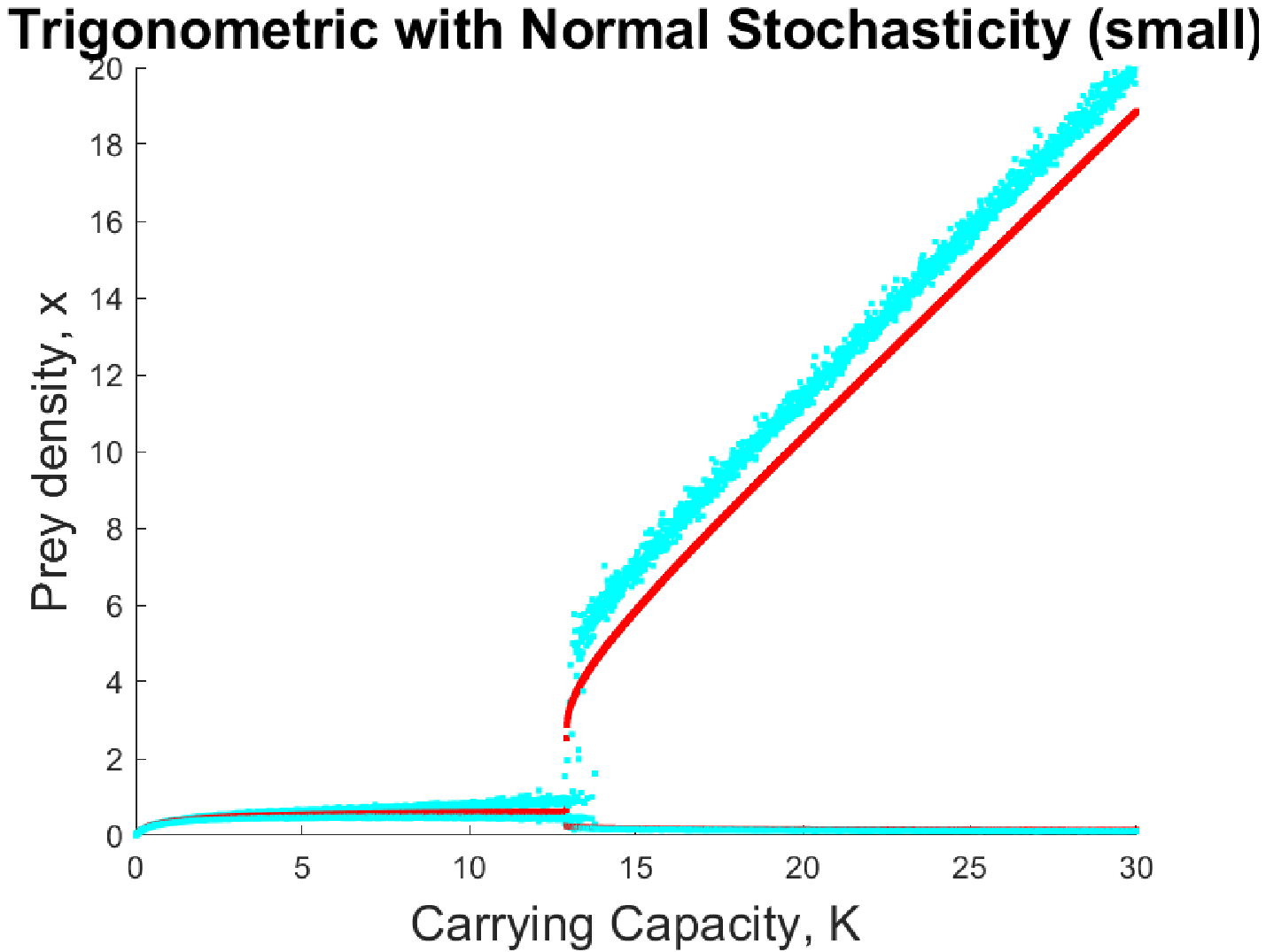}
    \\
    \includegraphics[width=0.325\linewidth]{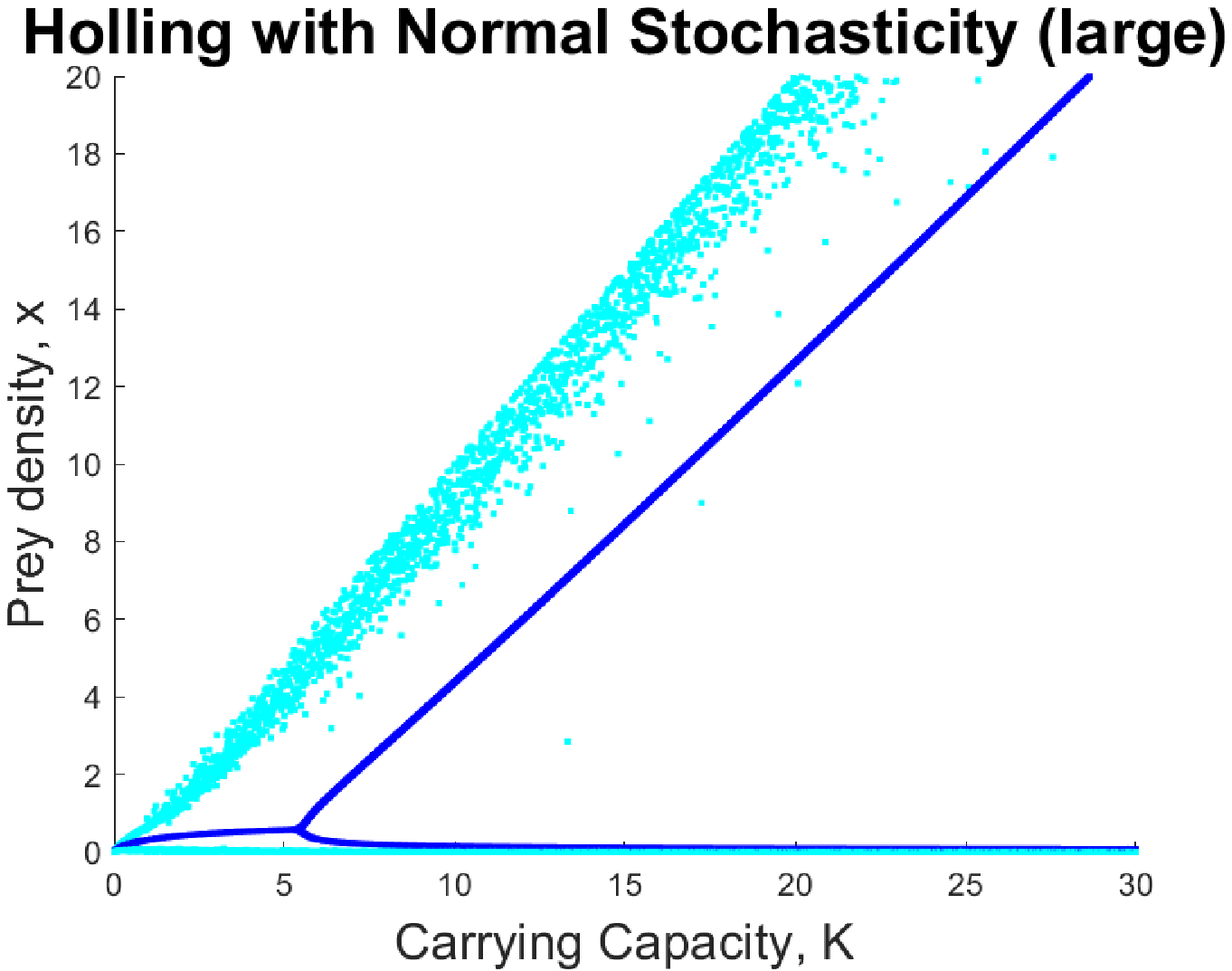}
    \includegraphics[width=0.325\linewidth]{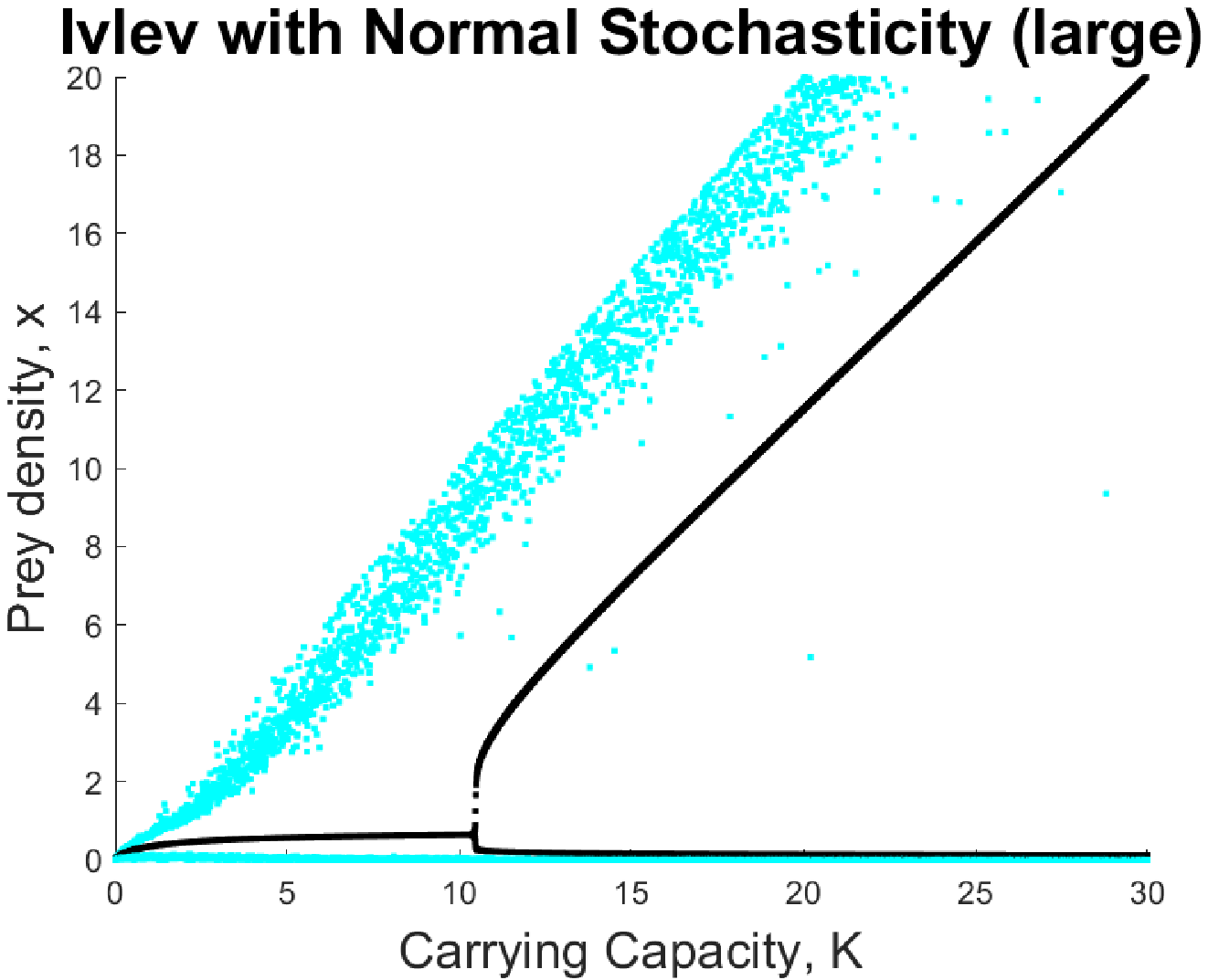}
    \includegraphics[width=0.325\linewidth]{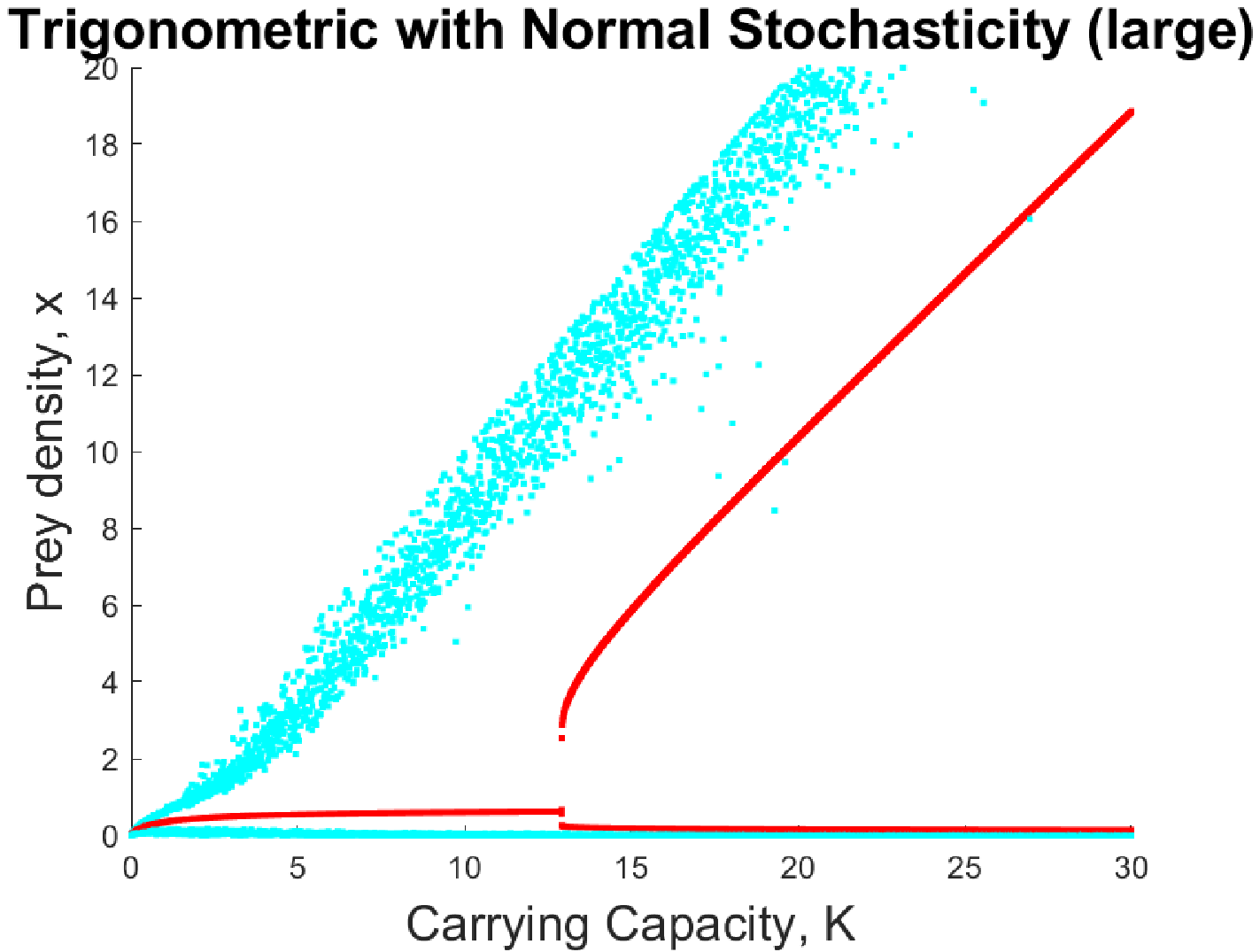}
    \\
    \includegraphics[width=\linewidth]{Pictures/FIGLeg3.eps}
    \caption{Top: Bifurcation diagrams for the stochastic Leslie-Gower-May model with $\sigma = 10$. Bottom: Bifurcation diagrams for the stochastic Leslie-Gower-May model with $\sigma = 100$. Both the deterministic and stochastic models were obtained by running time series and discarding the transient solutions, so these diagrams only show the stable behaviour.}
    \label{fig:StochBifurLargeSmall}
\end{figure}

In Figure \ref{fig:StochBifurLargeSmall} (top), we add normally distributed stochasticity to our models, with a smaller amplitude than we used in the manuscript. The stochastic bifurcation diagrams show no significant change from the deterministic diagrams. In Figure \ref{fig:StochBifurLargeSmall} (bottom), we add stochasticity with amplitudes larger than we used in the manuscript. These models show convergence to the prey-only and extinction steady states, and only show stable coexistence for small carrying capacities.

\section{Fitted Functional Response: Other fits}
\subsection{Functional Responses}
\begin{figure}[ht]
    \centering
    \includegraphics[width=0.4\linewidth]{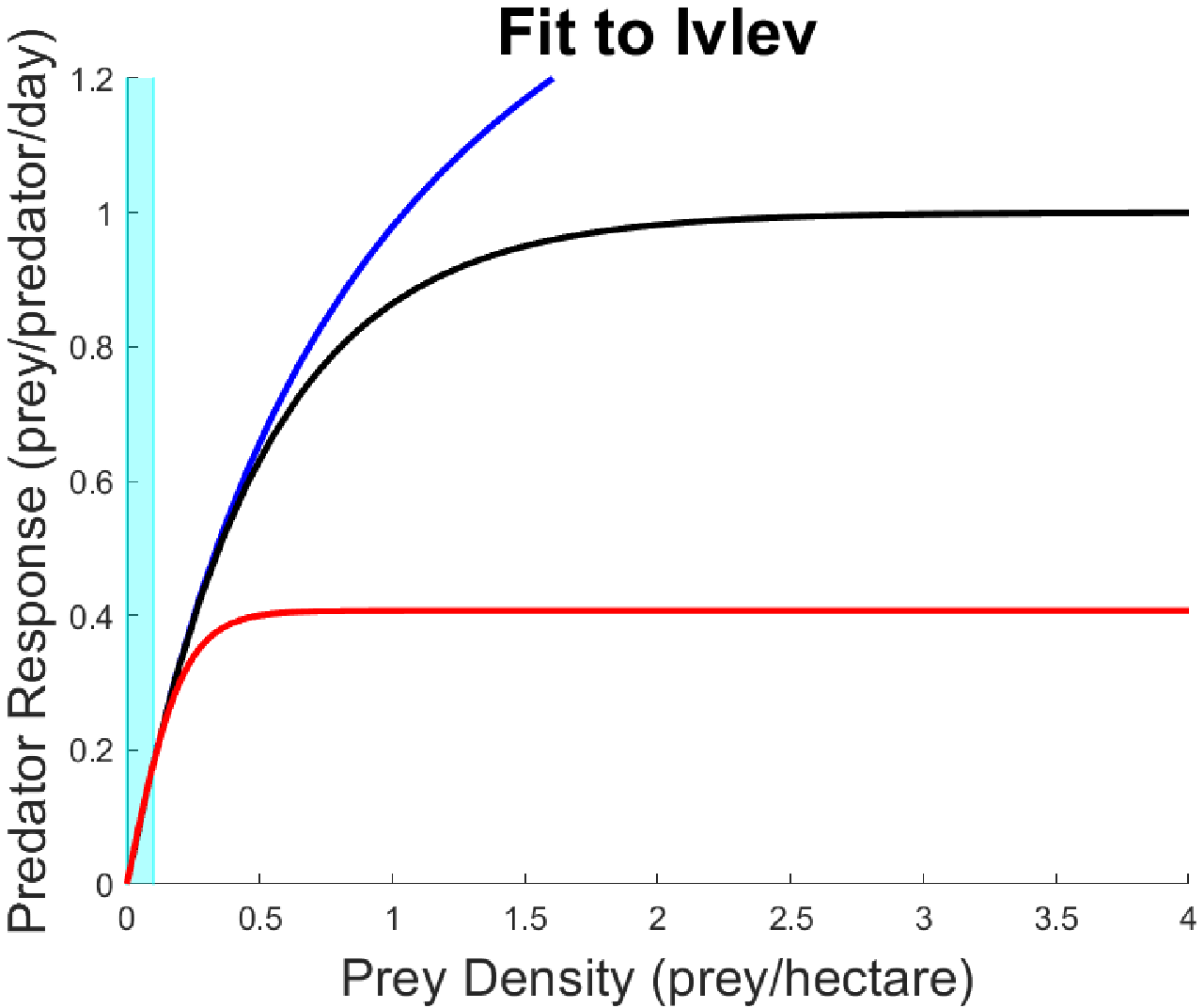}
    \includegraphics[width=0.4\linewidth]{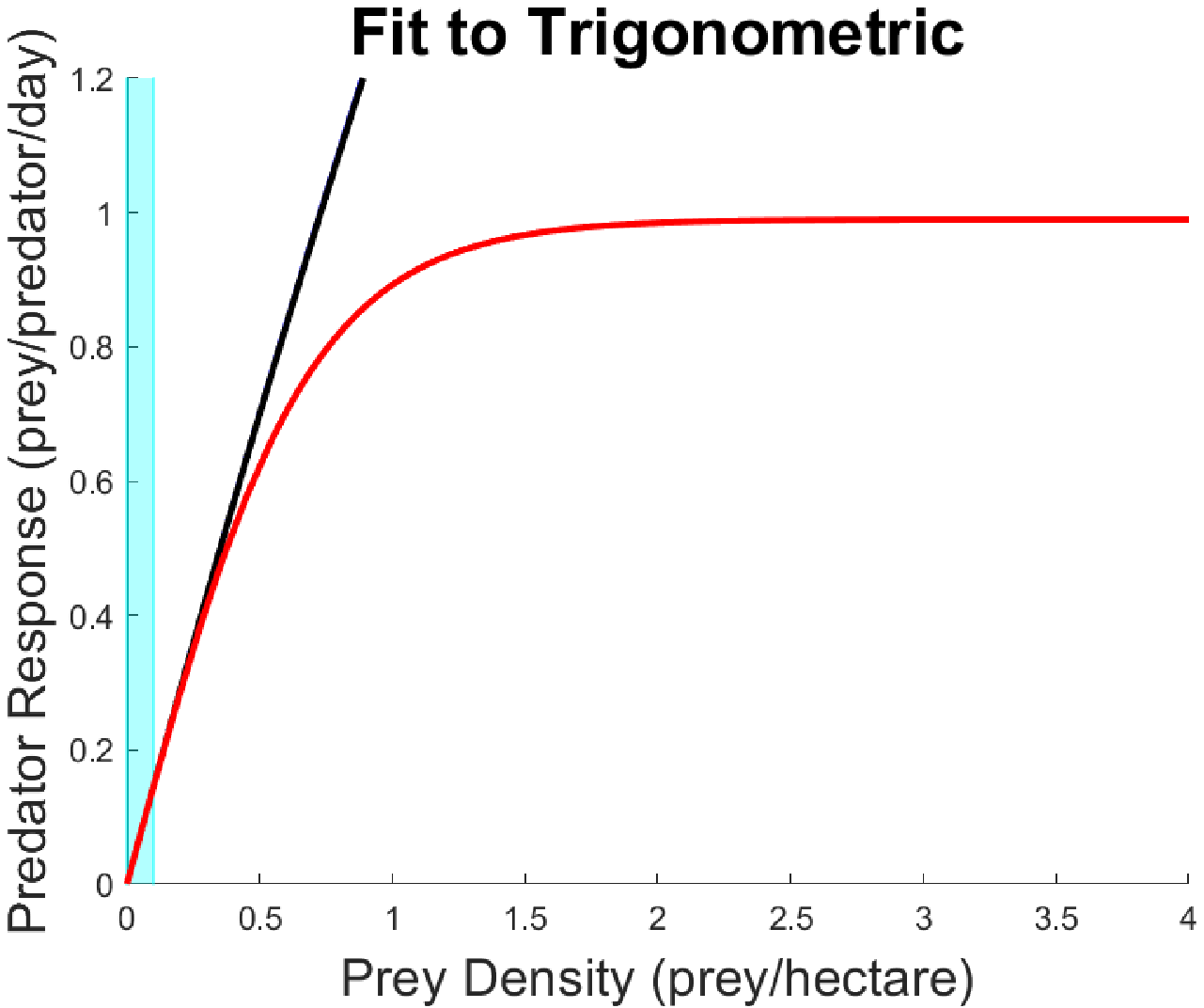}
    \\
    \includegraphics[width=0.4\linewidth]{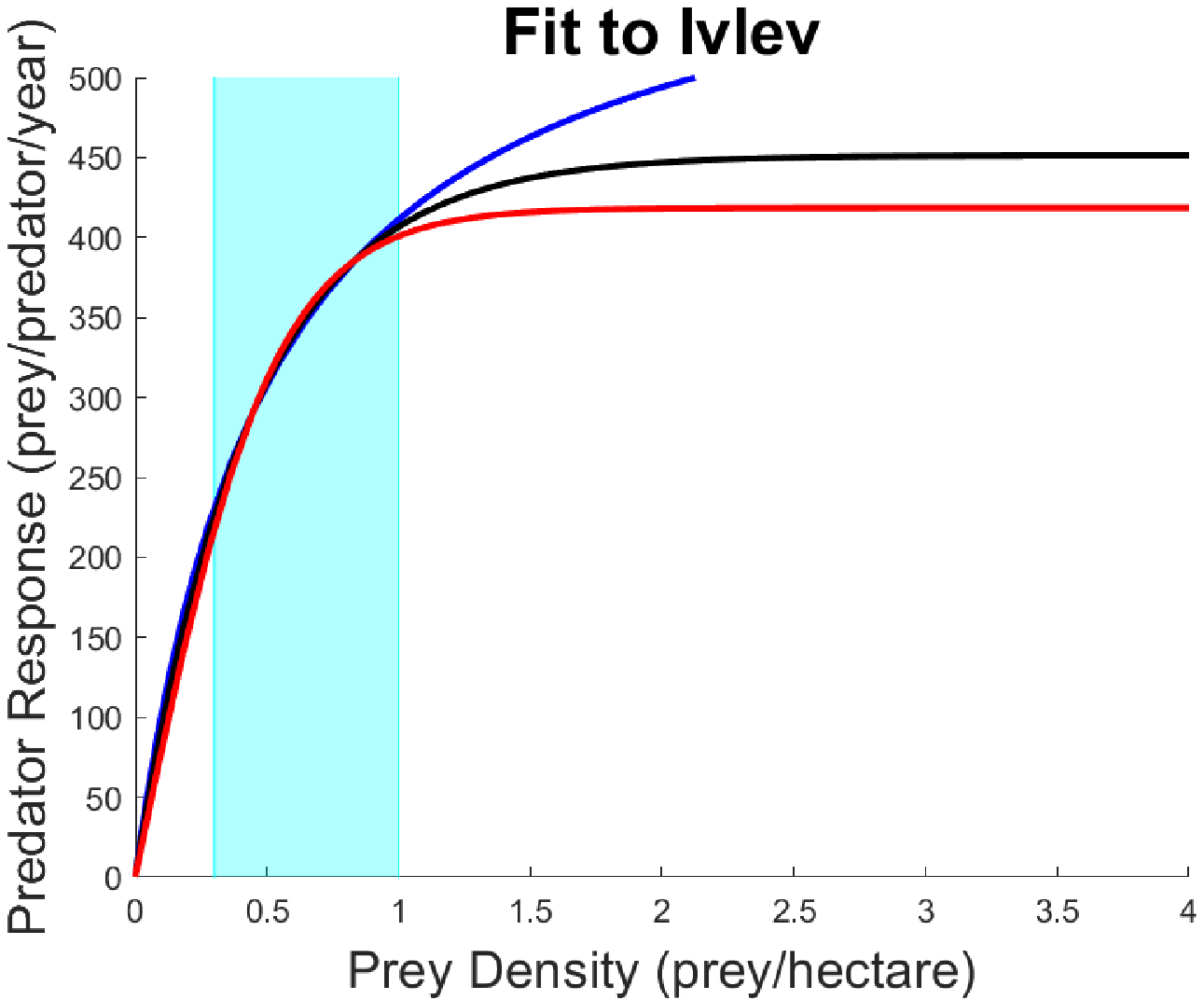}
    \includegraphics[width=0.4\linewidth]{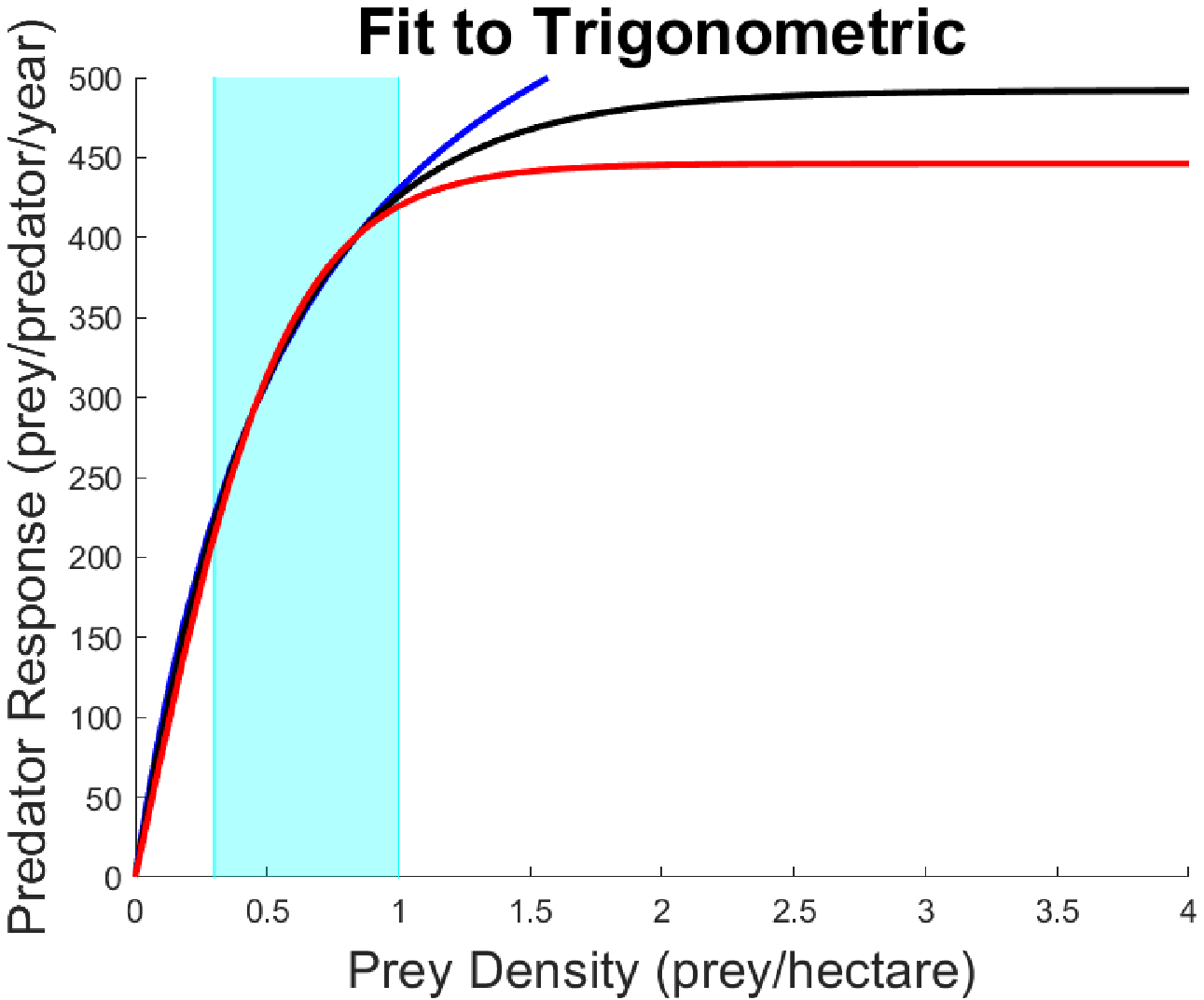}
    \includegraphics[width=0.7\linewidth]{Pictures/FIGLeg1.eps}
    \caption{Top: The Rosenzweig-MacArthur functional responses with nonlinear least squares applied over prey densities in [0,0.1]. Bottom: The Leslie-Gower-May functional responses with nonlinear least squares applied over prey densities in [0.3,1]. Both models are subject to the constraint $f_i(0)=0$ and fit to the specified functional response. The parameter values obtained from these fits are given in Table \ref{suptable:ParametersLS}. }
    \label{supfig:SSCurves}
\end{figure}

%-----------------------------------------------------------------
As in the manuscript, we fit the functional responses to a fixed functional response for prey density values in a region around the coexistence steady state. The resulting functional responses are given in Figure \ref{supfig:SSCurves} and the parameters found in the fit are given in Table \ref{suptable:ParametersLS}.

\begin{table}[ht]
\centering
\begin{tabular}{|c|c|c|c|c|c|c|}
    \hline 
   & \multicolumn{3}{|c|}{Rosenzweig-MacArthur} & \multicolumn{3}{|c|}{Leslie-Gower-May} \\
    \hline
     & Holling & Ivlev & Trigonometric & Holling & Ivlev & Trigonometric \\
    \hline
    
    $a_H$ & \textbf{3.05} & 2.0026 & 1.4695 & \textbf{1683.333} & 1228.2780 & 1114.6253\\
    $b_H$ & \textbf{2.68} & 1.0453 & 0.0978 & \textbf{3.333} & 1.9855 & 1.5905\\
    \hline
    
    $a_I$ & 0.6282 & \textbf{1} & 7.5417 & 399.8009 & \textbf{451.447} & 492.0405\\
    $b_I$ & 4.8204 & \textbf{2} & 0.1948 & 3.1656 & \textbf{2.313} & 2.0060\\
    \hline
    
    $a_T$ & 0.3707 & 0.4069 & \textbf{0.99} & 384.1465 & 418.5697 & \textbf{446.182}\\
    $b_T$ & 7.6281 & 4.7698 & \textbf{1.48} & 2.4028 & 1.9188 & \textbf{1.743}\\
    \hline
\end{tabular}
    \caption{Values of the parameters obtained from the fitted functional response where a specified functional response is fixed and the other functional responses are fit to the fixed response via the parameters $a_i$ and $b_i$. The entries in the second row give the fixed functional response for each corresponding column. The parameters, $a_i$ and $b_i$, used in the fixed sub-models are given in bold.}
    \label{suptable:ParametersLS}
\end{table}

\subsection{Bifurcation Diagrams}
\begin{figure}[ht]
    \centering
    \begin{subfigure}{\linewidth}
        \includegraphics[width=0.325\linewidth]{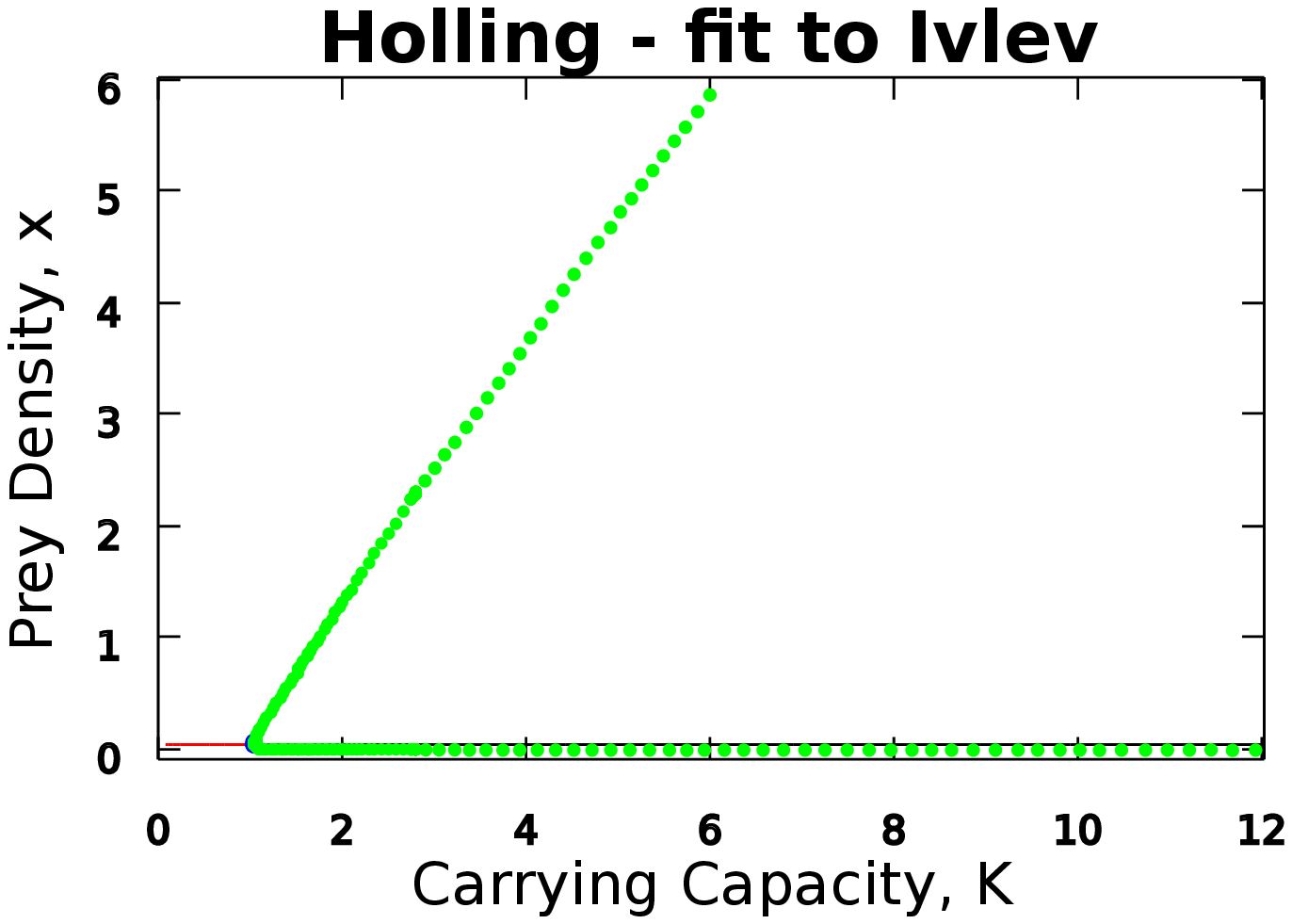}
        \includegraphics[width=0.325\linewidth]{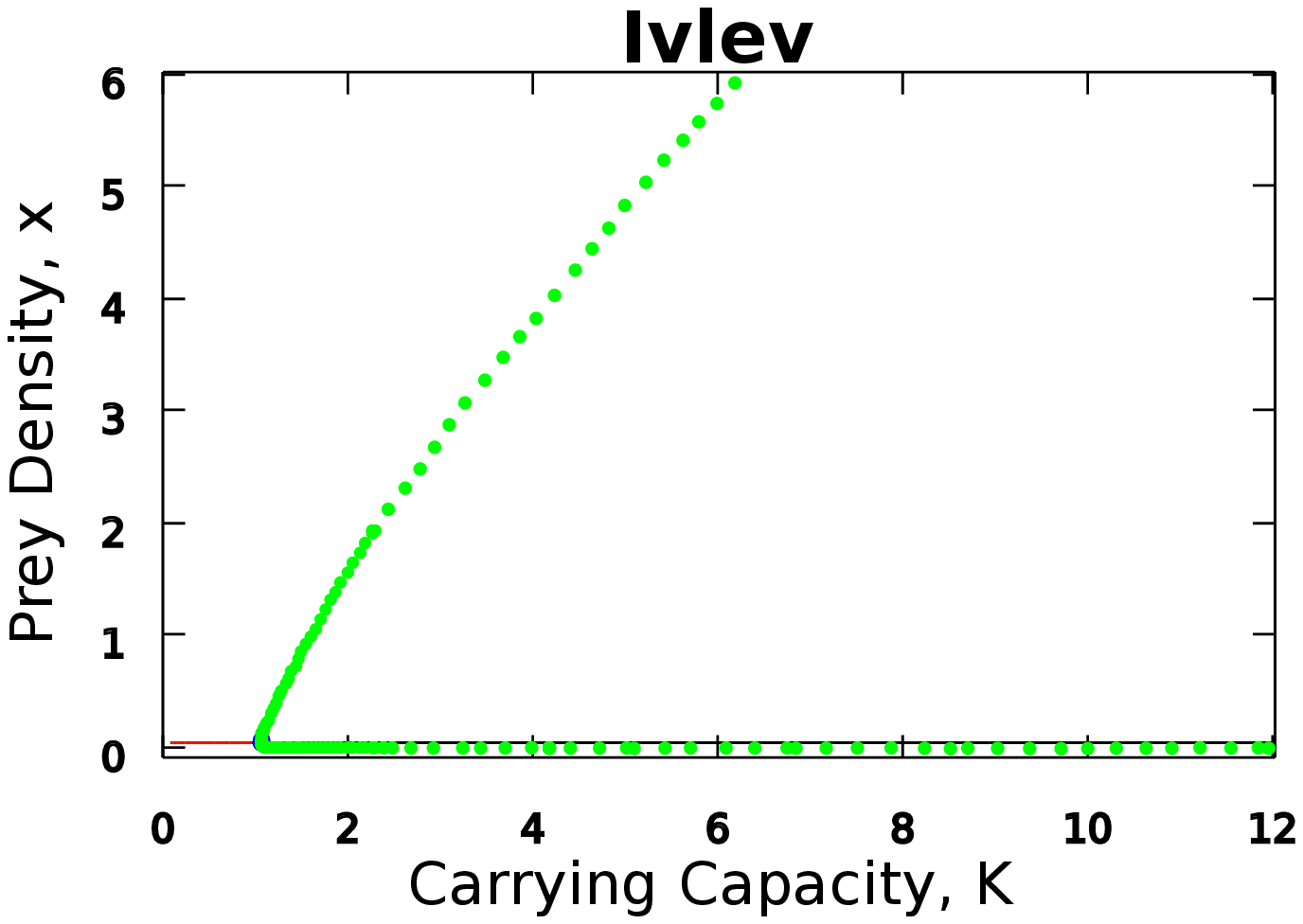}
        \includegraphics[width=0.325\linewidth]{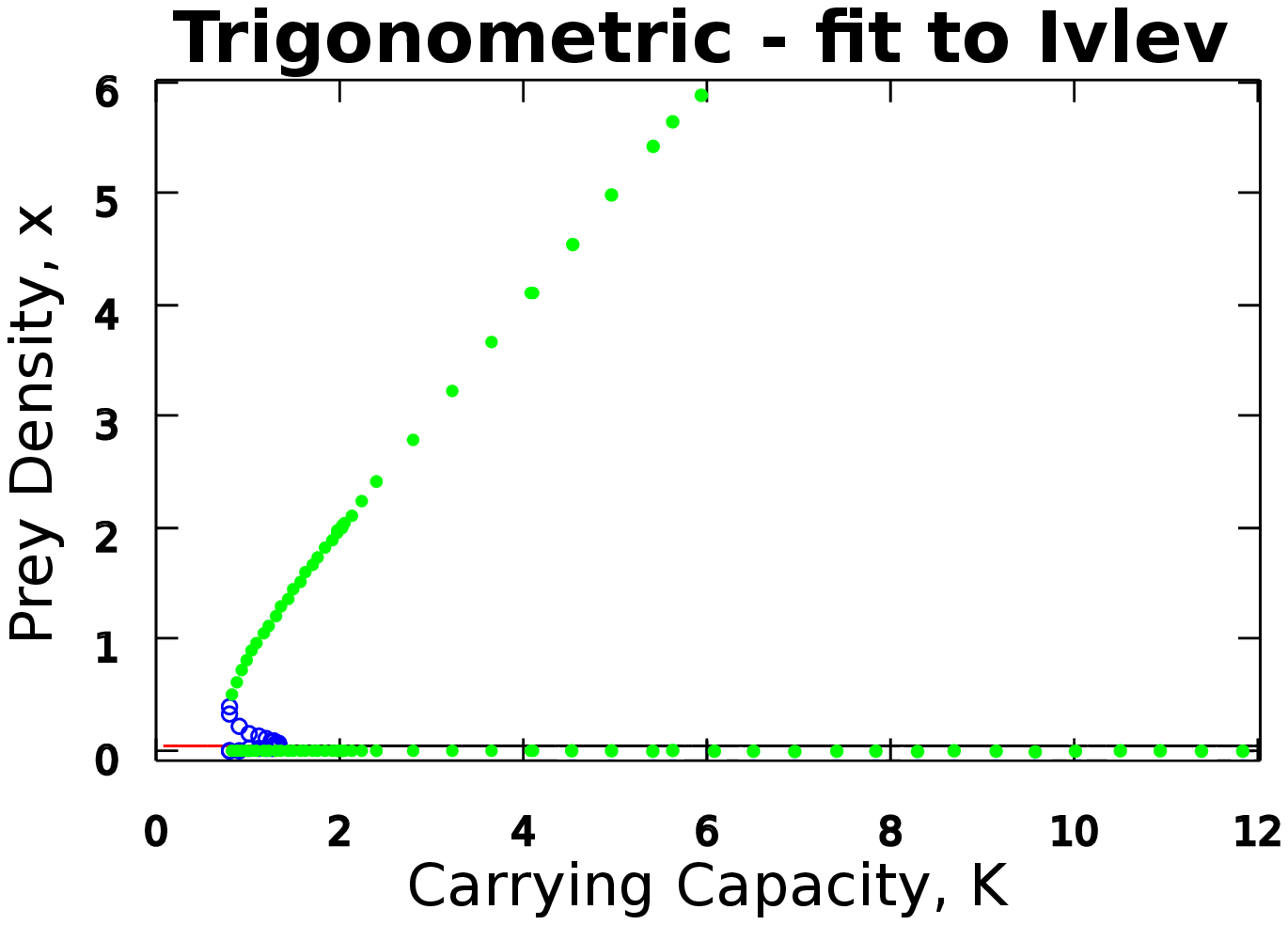}
        \subcaption{Rosenzweig-MacArthur model with Ivlev fixed}
    \end{subfigure}
    \\
    \begin{subfigure}{\linewidth}
        \includegraphics[width=0.325\linewidth]{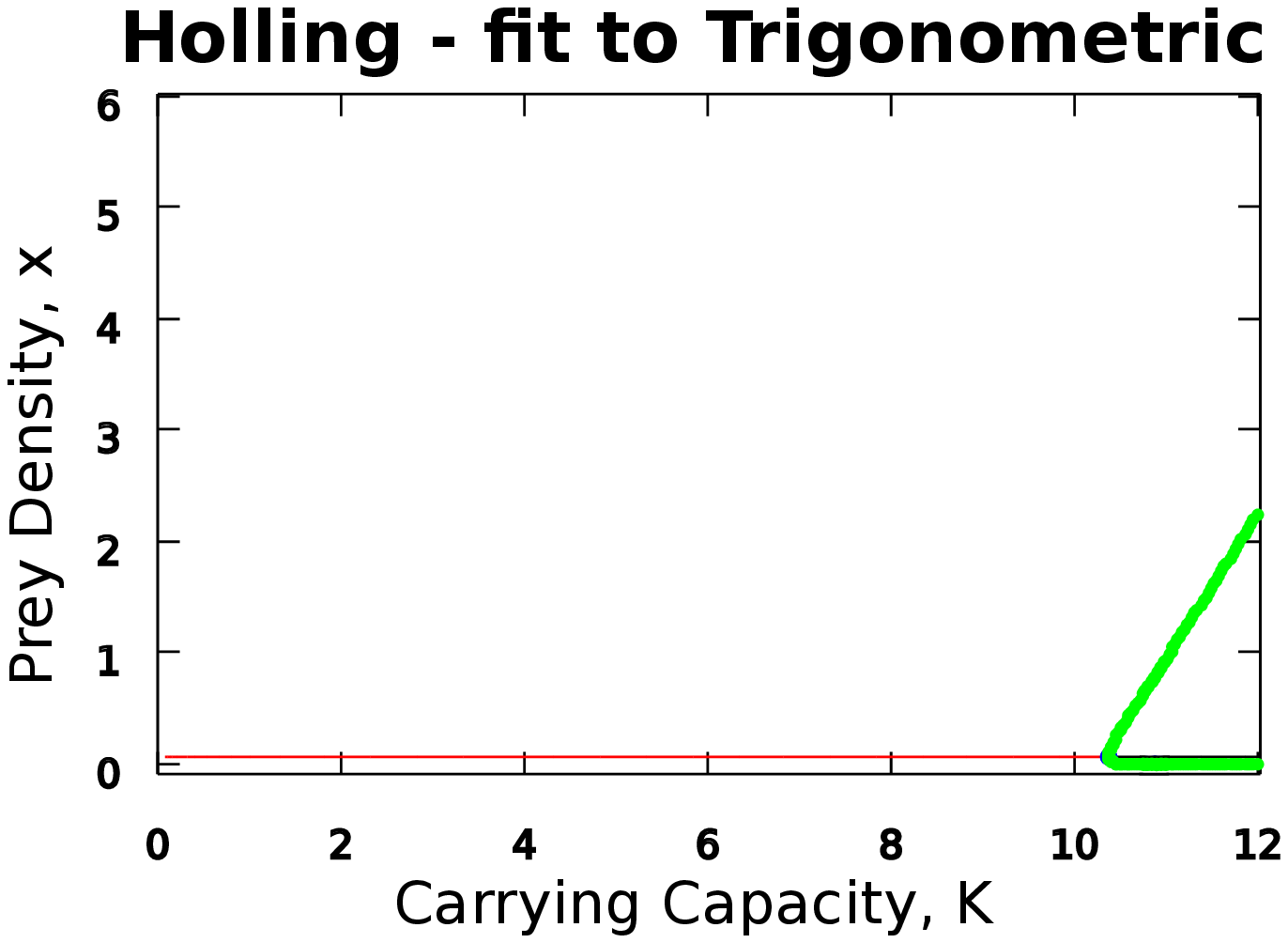}
        \includegraphics[width=0.325\linewidth]{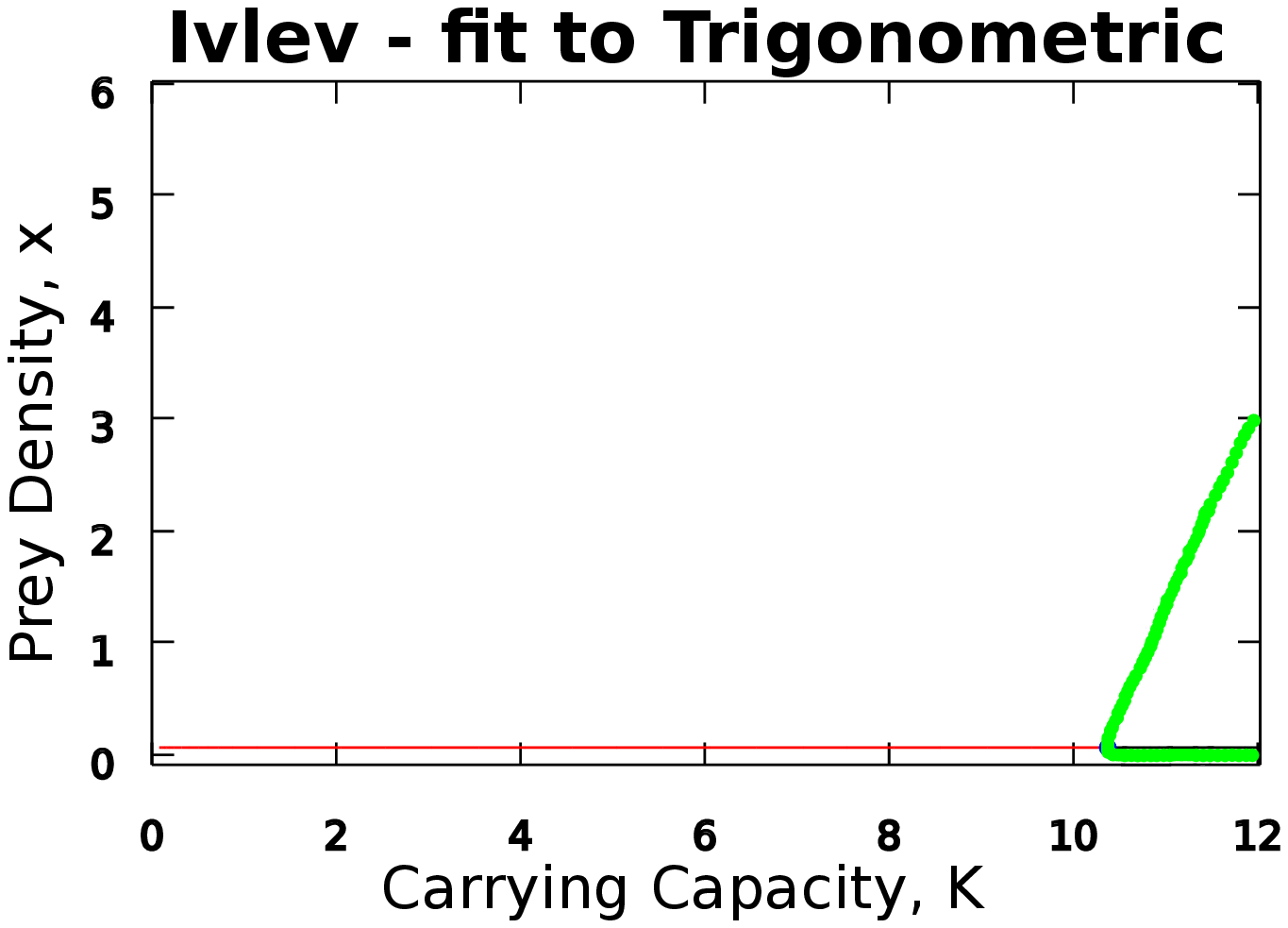}
        \includegraphics[width=0.325\linewidth]{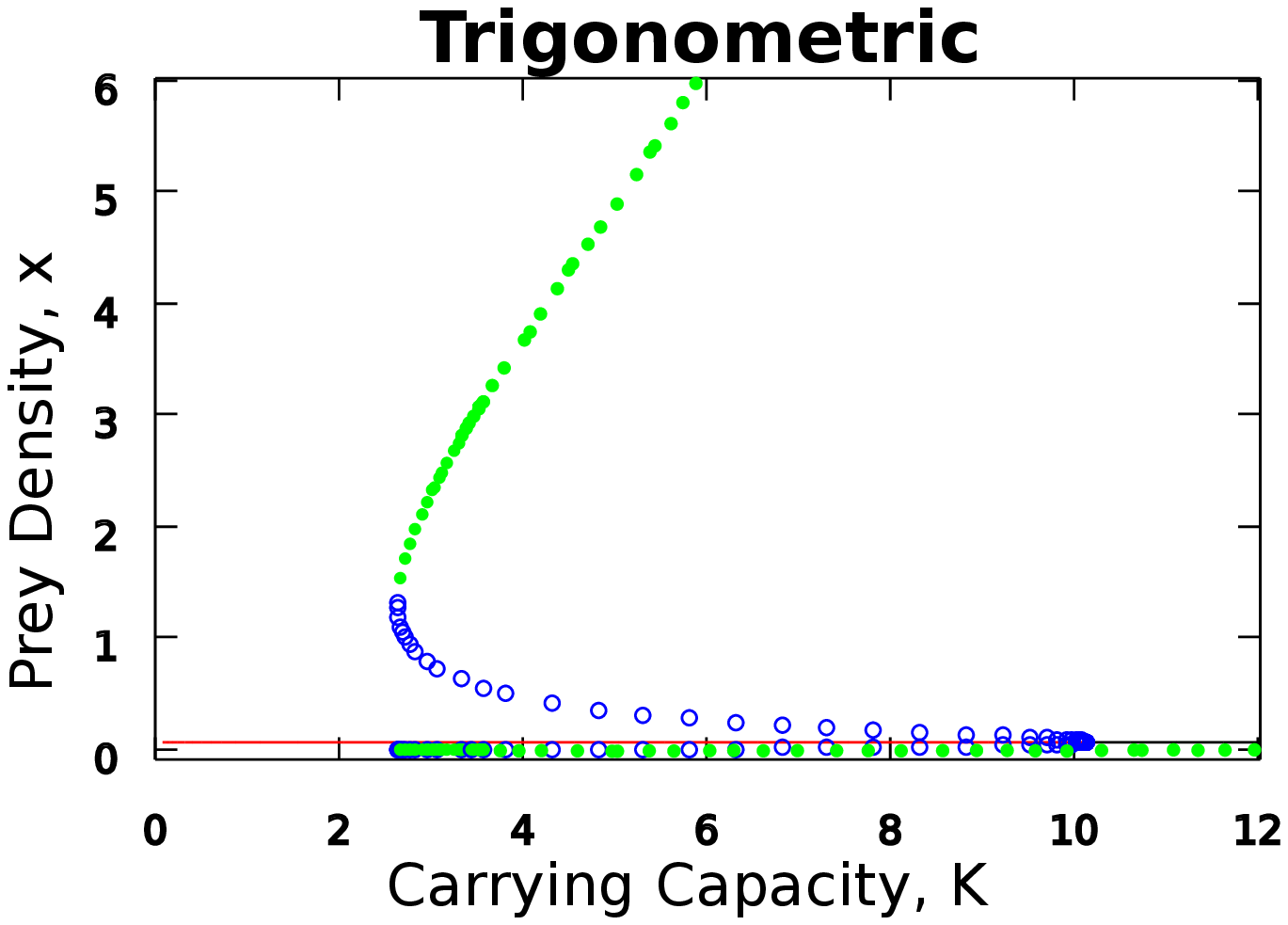}
         \subcaption{Rosenzweig-MacArthur model with trigonometric fixed}
    \end{subfigure}
    \\
    \begin{subfigure}{\linewidth}
        \includegraphics[width=0.325\linewidth]{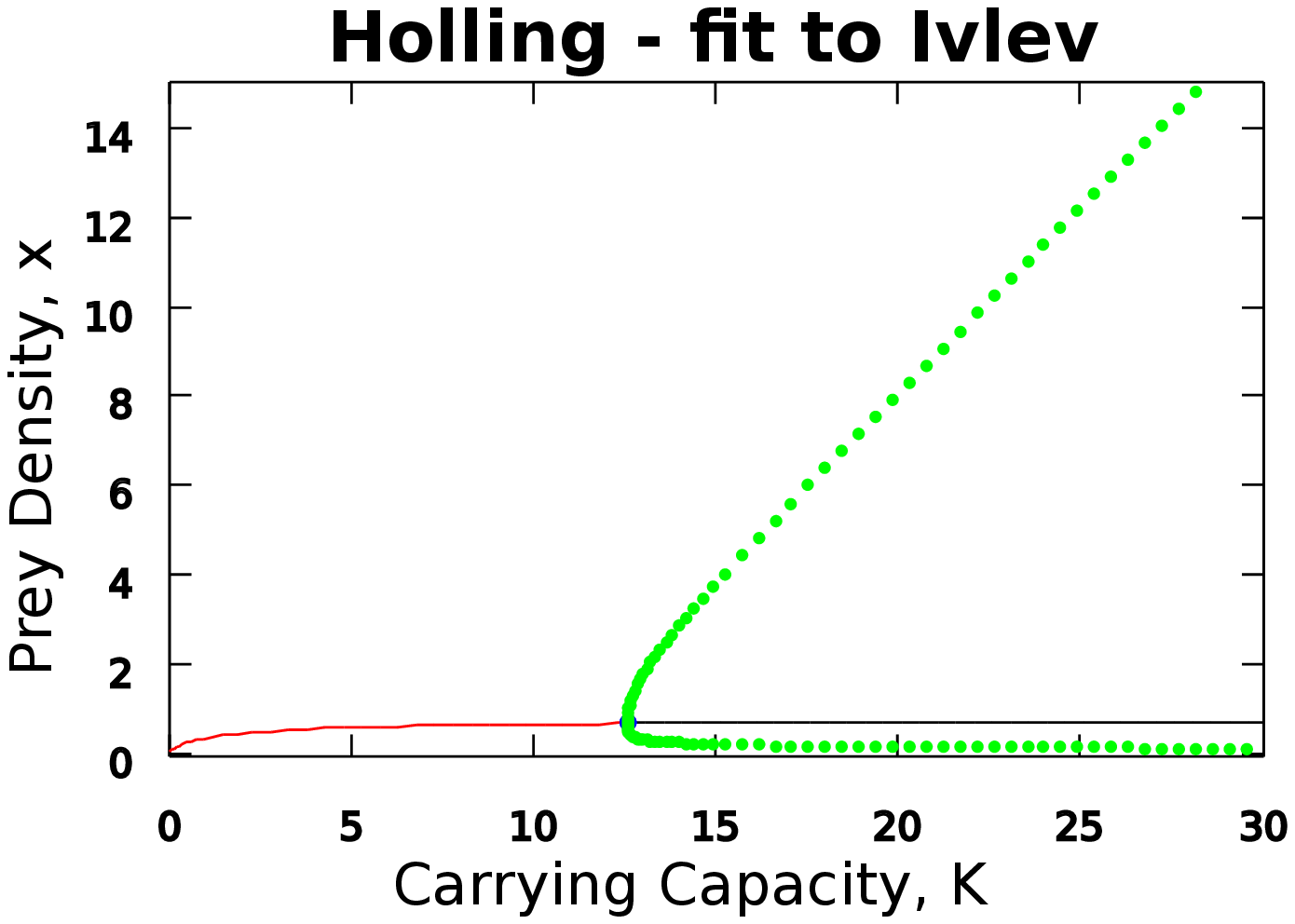}
        \includegraphics[width=0.325\linewidth]{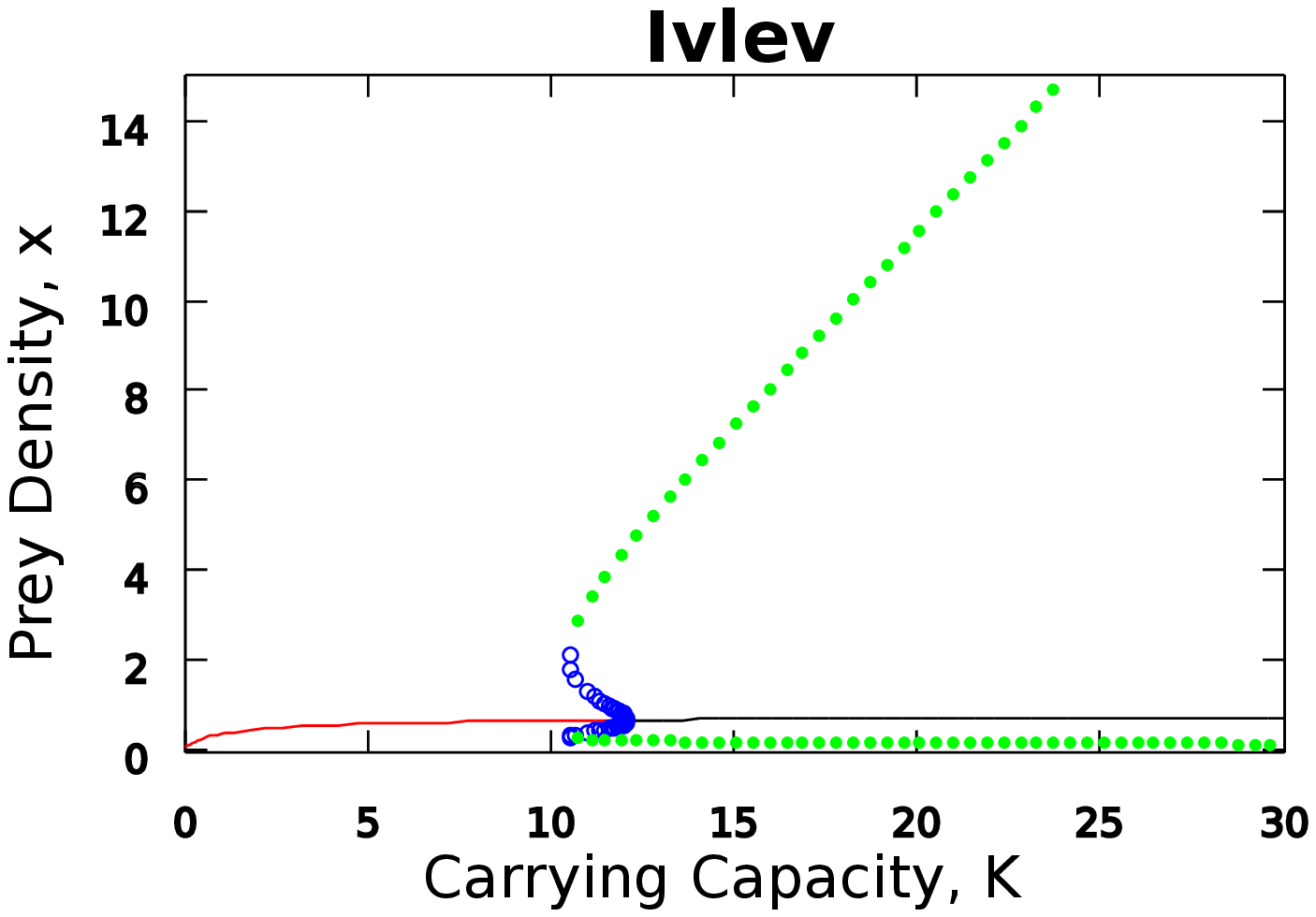}
        \includegraphics[width=0.325\linewidth]{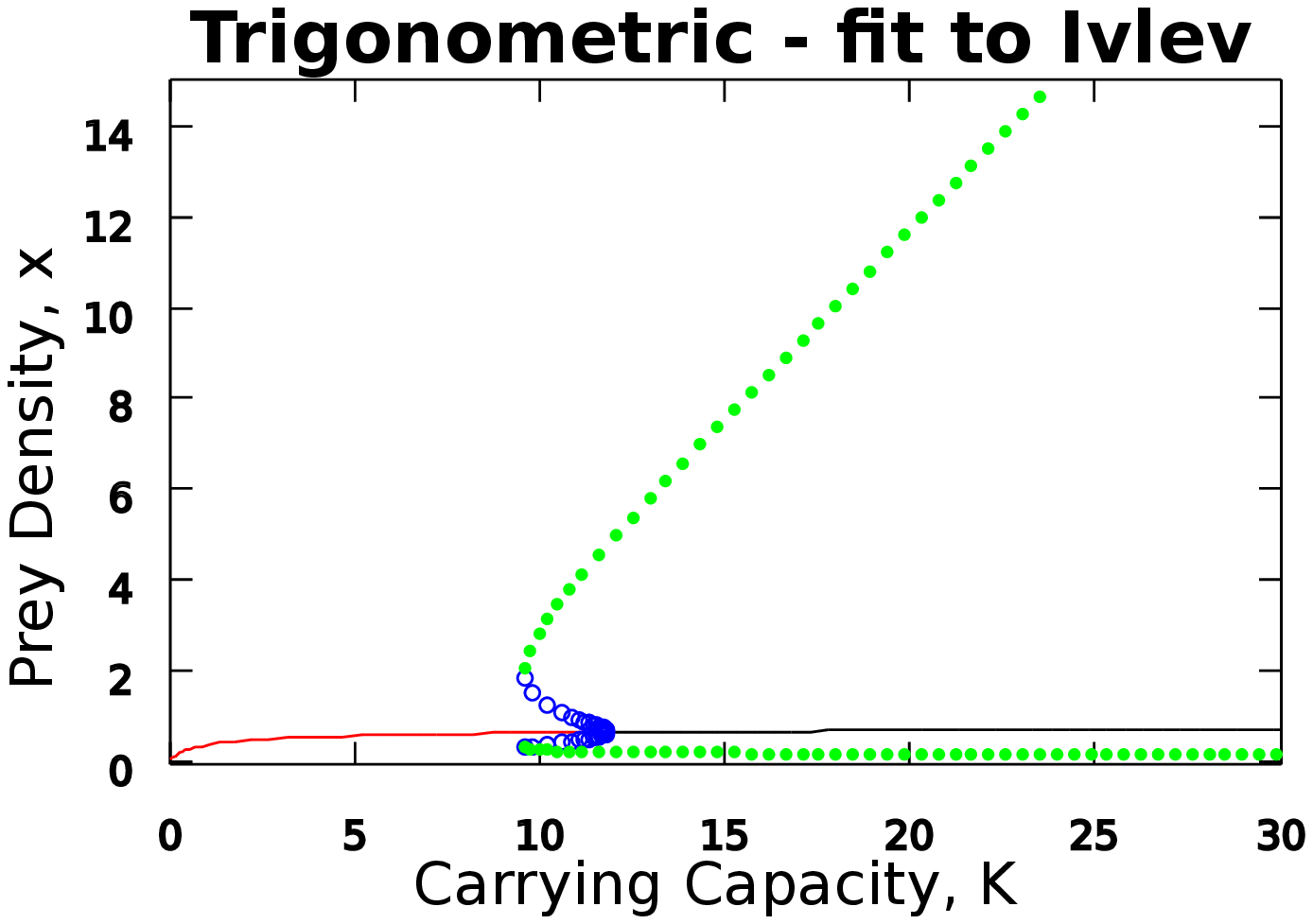}
        \subcaption{Leslie-Gower-May model with Ivlev fixed}
    \end{subfigure}
    \\
    \begin{subfigure}{\linewidth}
        \includegraphics[width=0.325\linewidth]{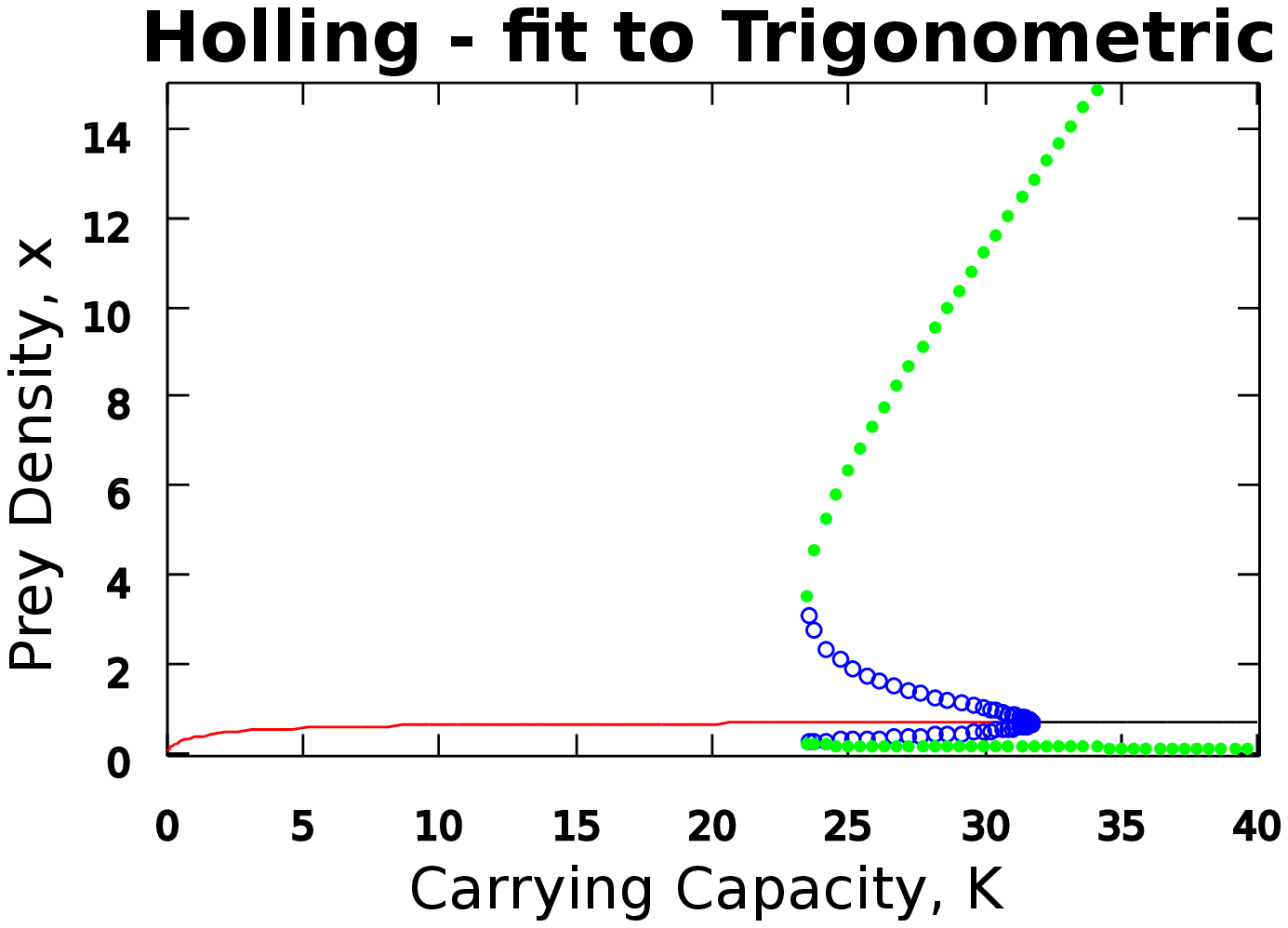}
        \includegraphics[width=0.325\linewidth]{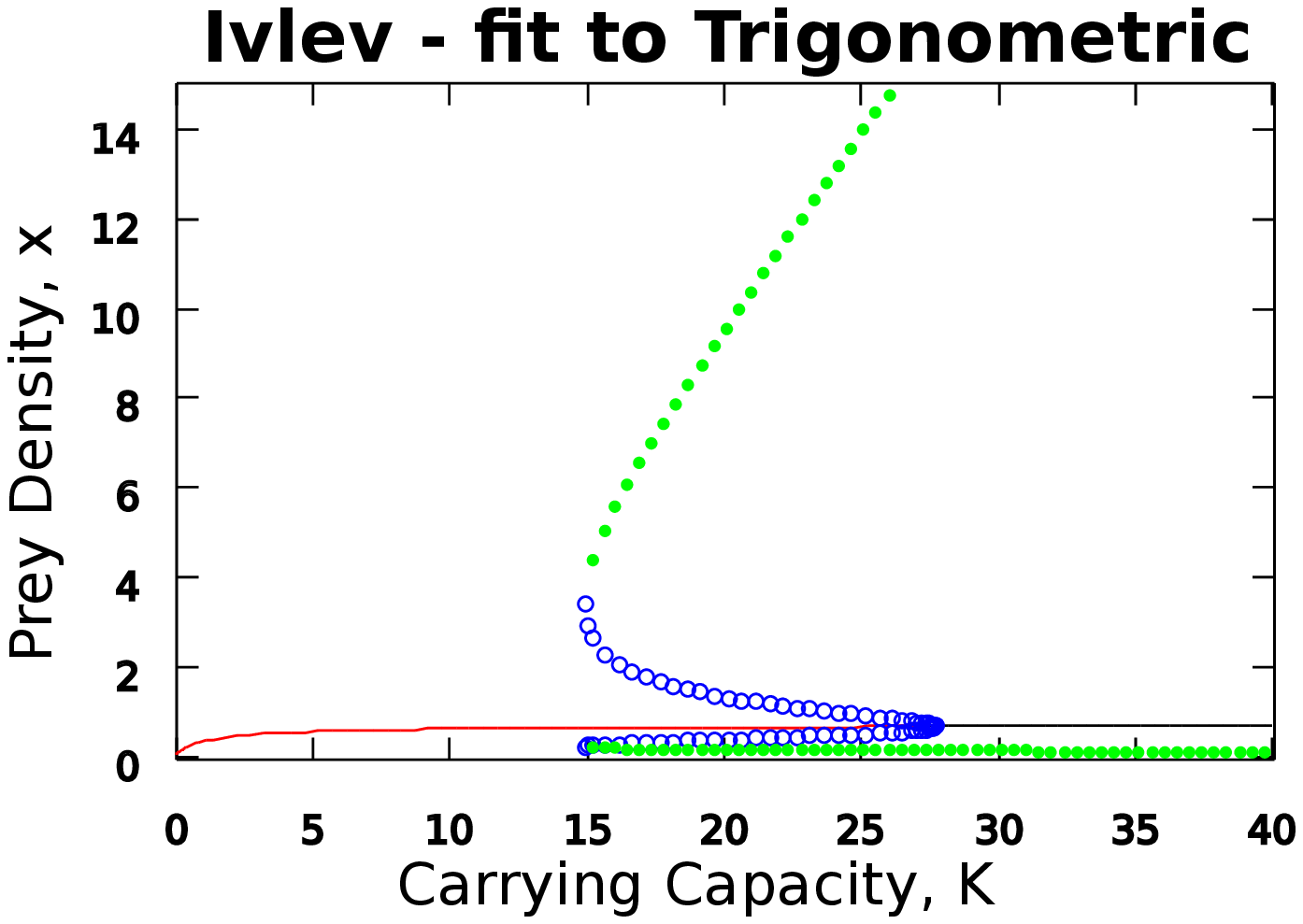}
        \includegraphics[width=0.325\linewidth]{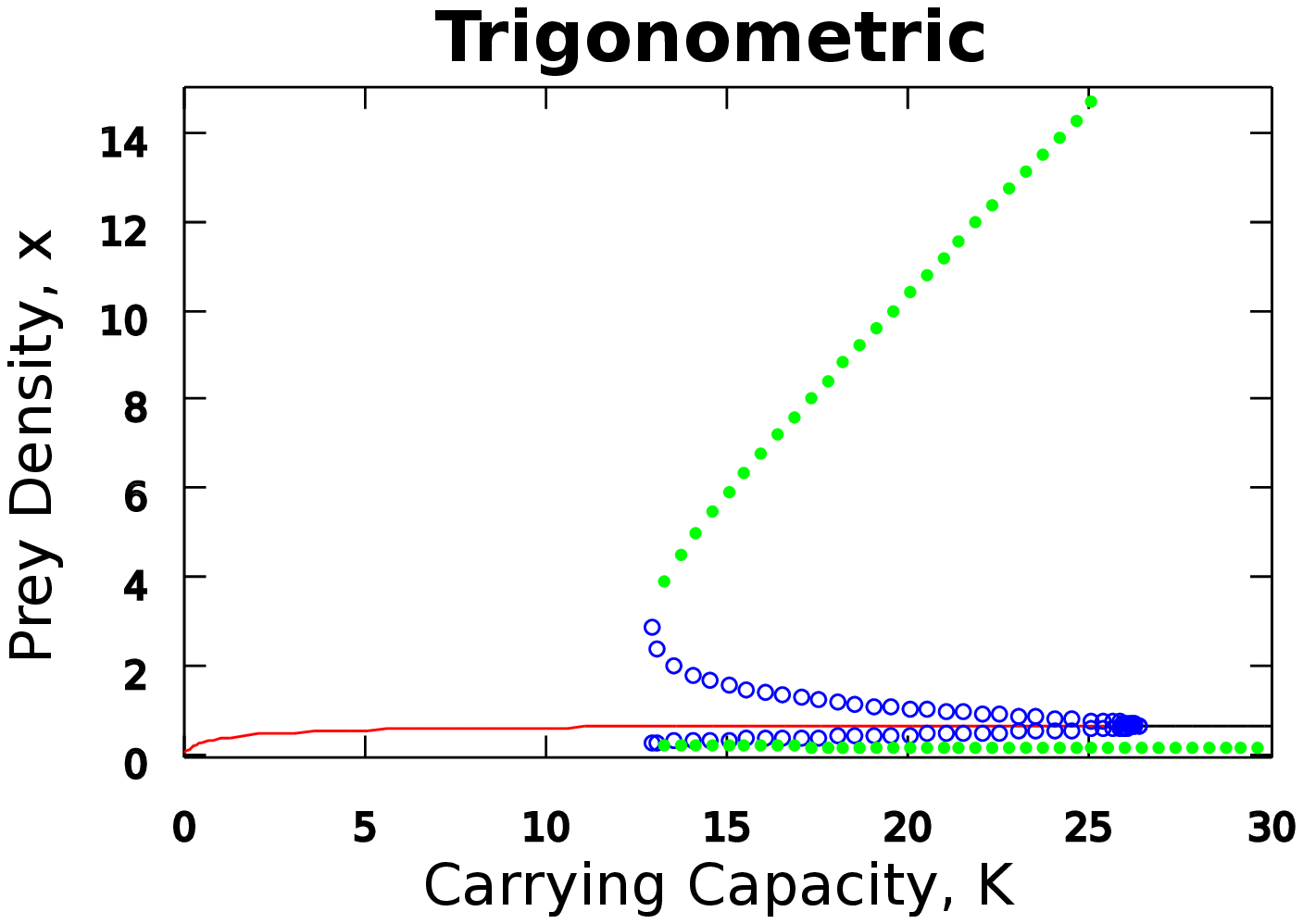}
        \subcaption{Leslie-Gower-May model with trigonometric fixed}
    \end{subfigure}
    \\
    \includegraphics[width=\linewidth]{Pictures/FIGLeg2.eps}
    \caption{Bifurcation diagrams for the Rosenzweig-MacArthur and  Leslie-Gower-May models with least squares applied to the functional responses over the prey density intervals $[0,0.1]$ and $[0.3,1]$, respectively. The parameter values used are given in Table \ref{suptable:ParametersLS}.}
    \label{supfig:CoexSSBifur}
\end{figure}

%-----------------------------------------------------------------
The bifurcation diagrams corresponding to the fitted functional response cases are given in Figure \ref{supfig:CoexSSBifur}. Both the Rosenzweig-MacArthur and the Leslie-Gower-May models show convergence of model behaviour towards the fixed model.

\section{Piecewise Identical Functional Response: Other Fits}
\subsection{Functional Responses}
\begin{figure}[ht]
    \centering
    \begin{subfigure}{\linewidth}
        \includegraphics[width=0.325\linewidth]{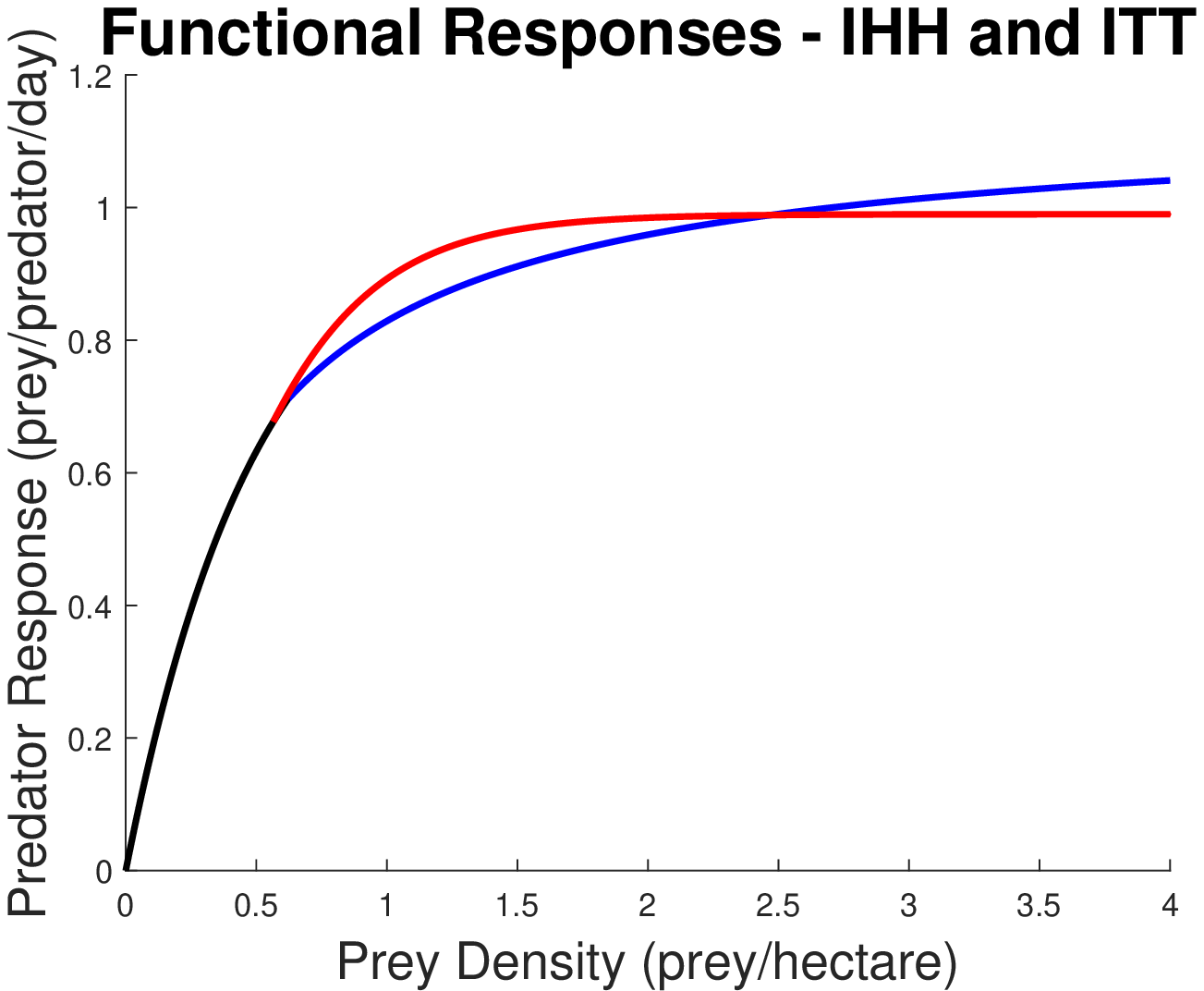}
        \includegraphics[width=0.325\linewidth]{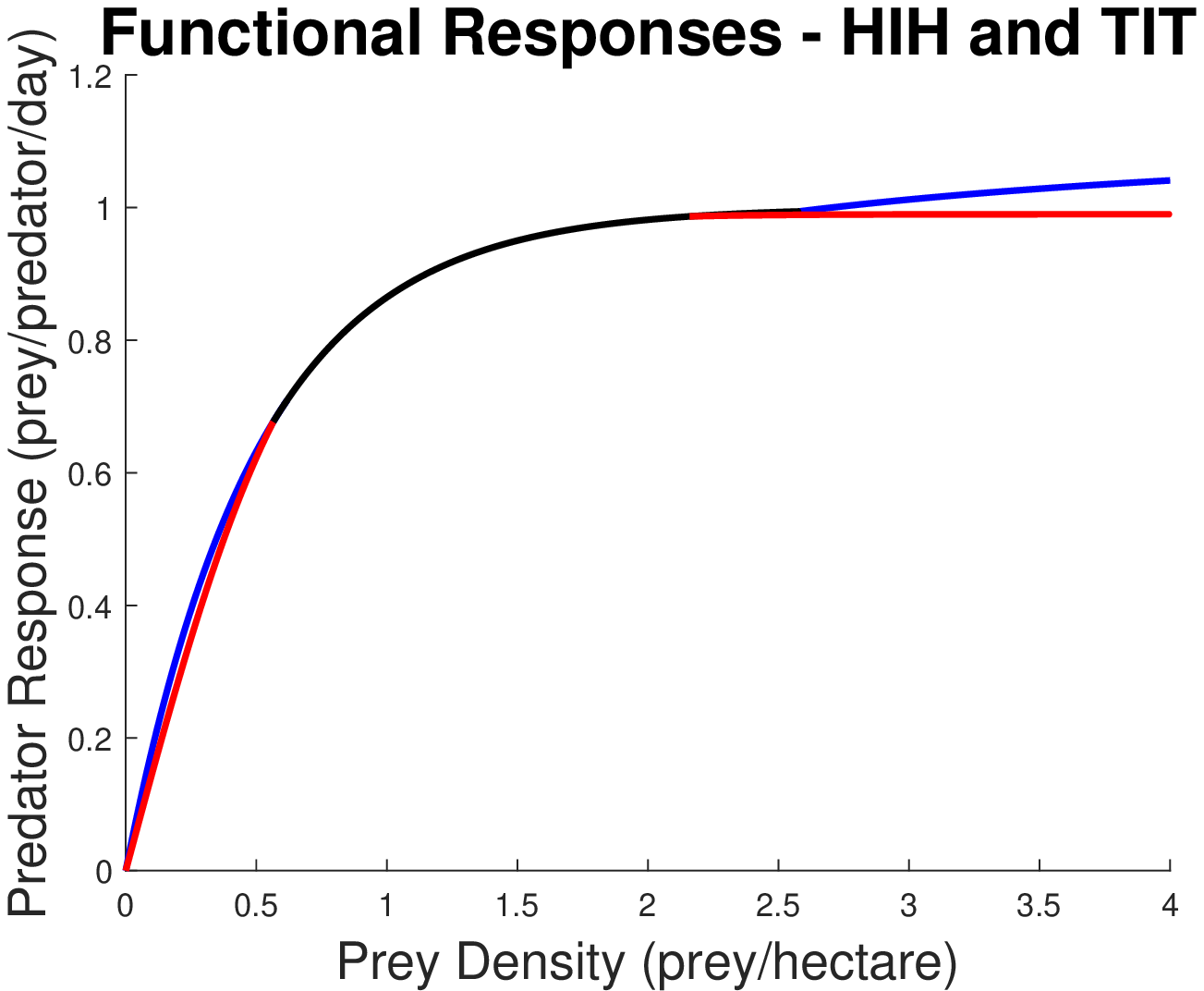}
        \includegraphics[width=0.325\linewidth]{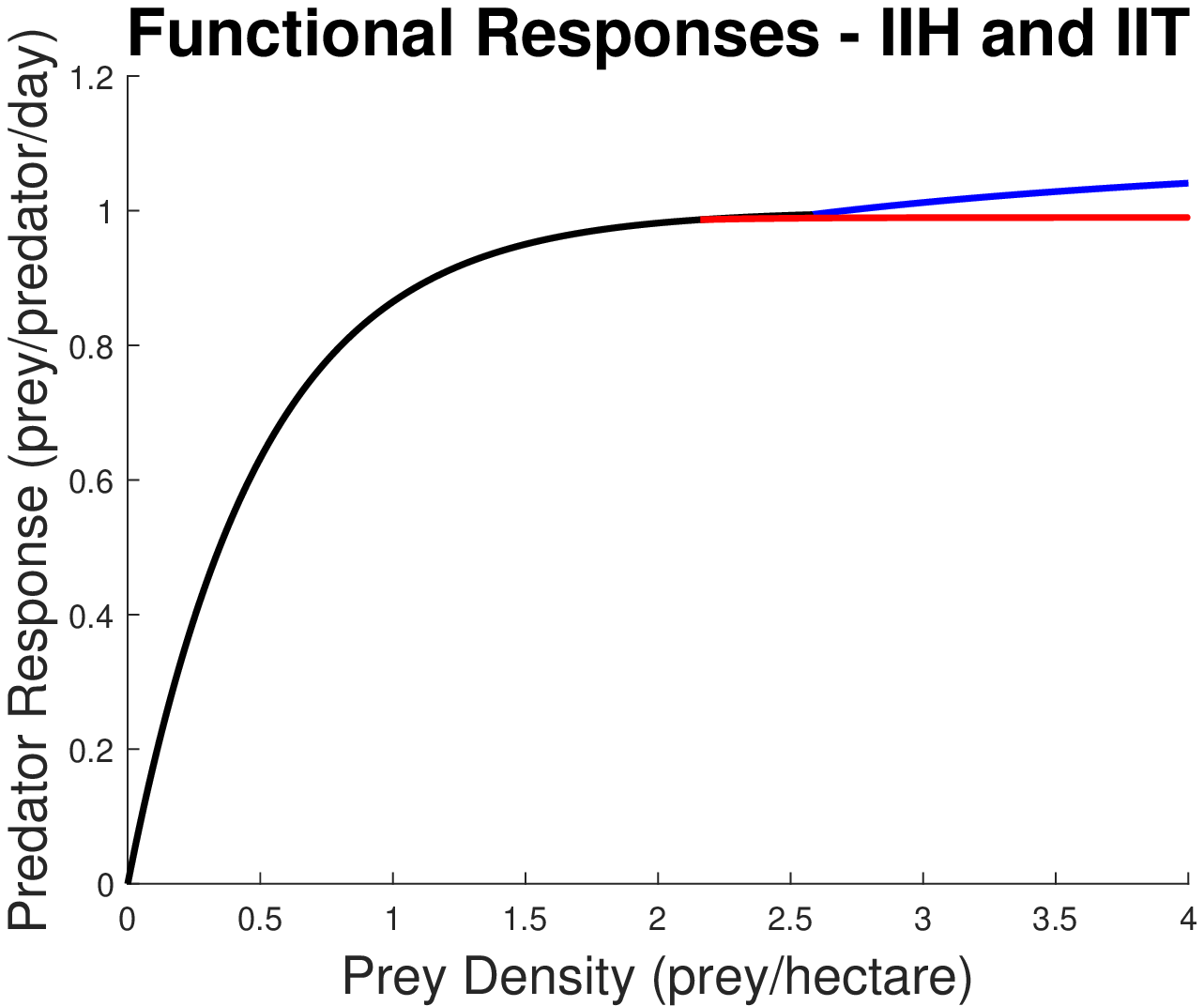}
        \\
        \includegraphics[width=0.325\linewidth]{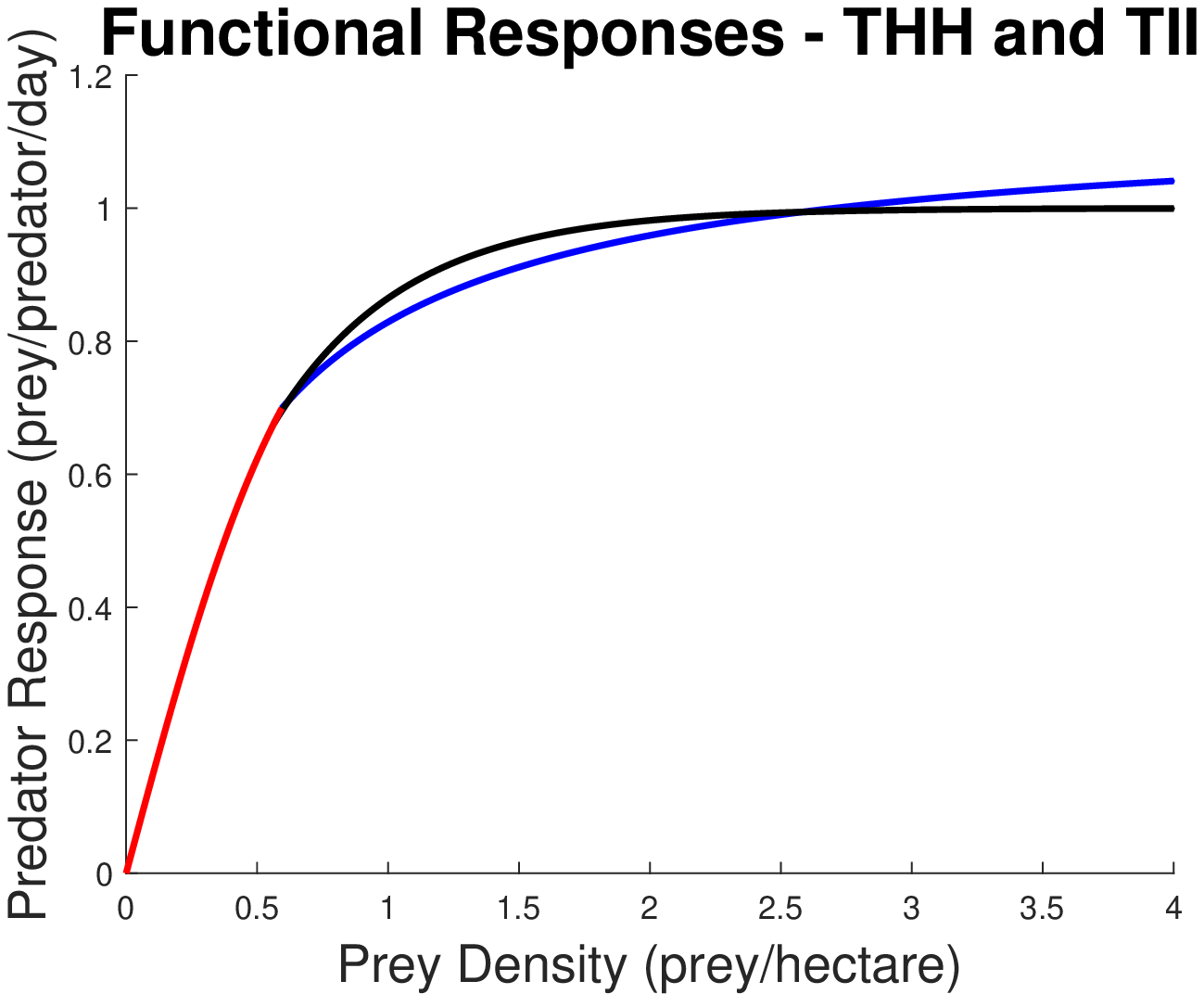}
        \includegraphics[width=0.325\linewidth]{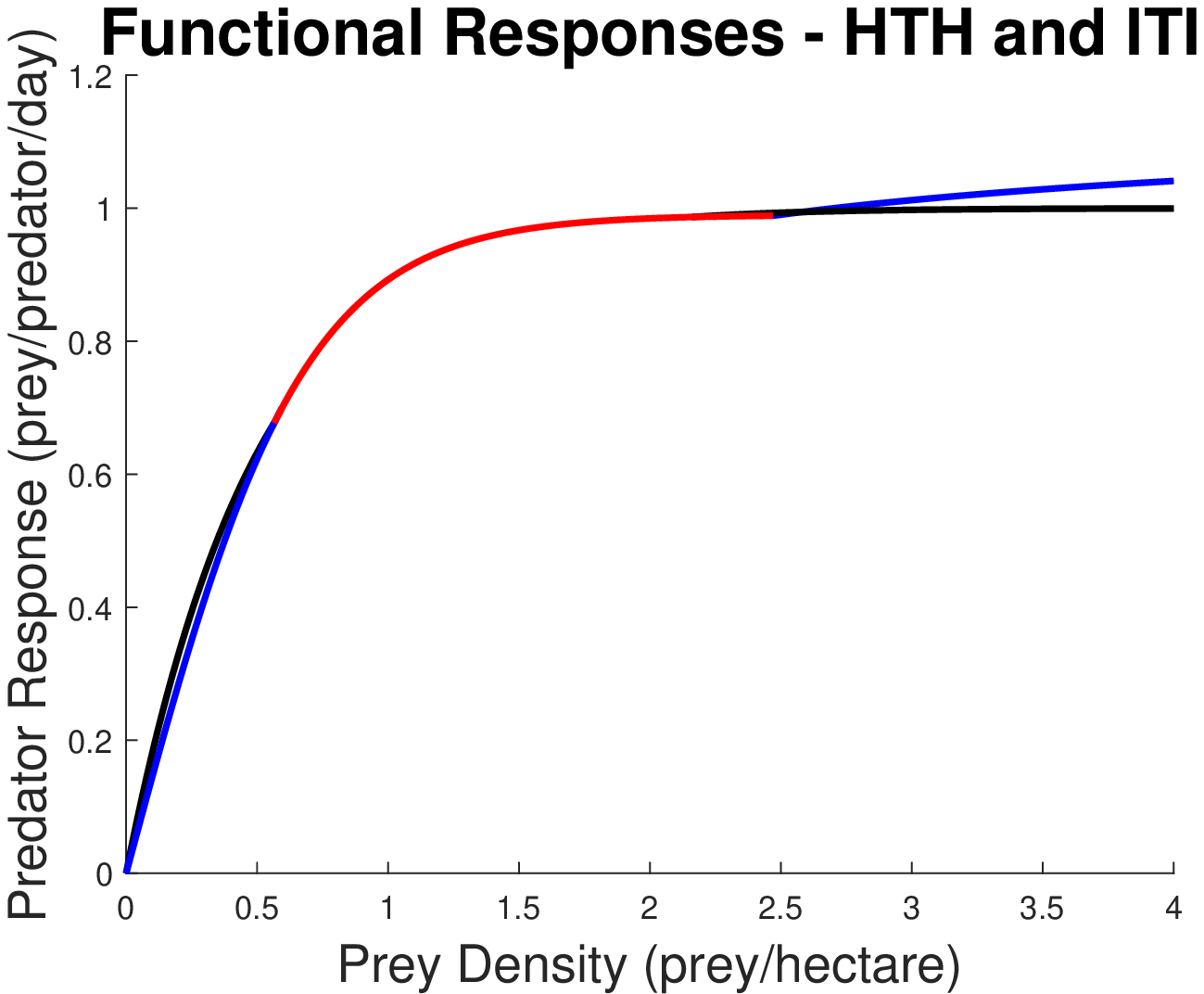}
        \includegraphics[width=0.325\linewidth]{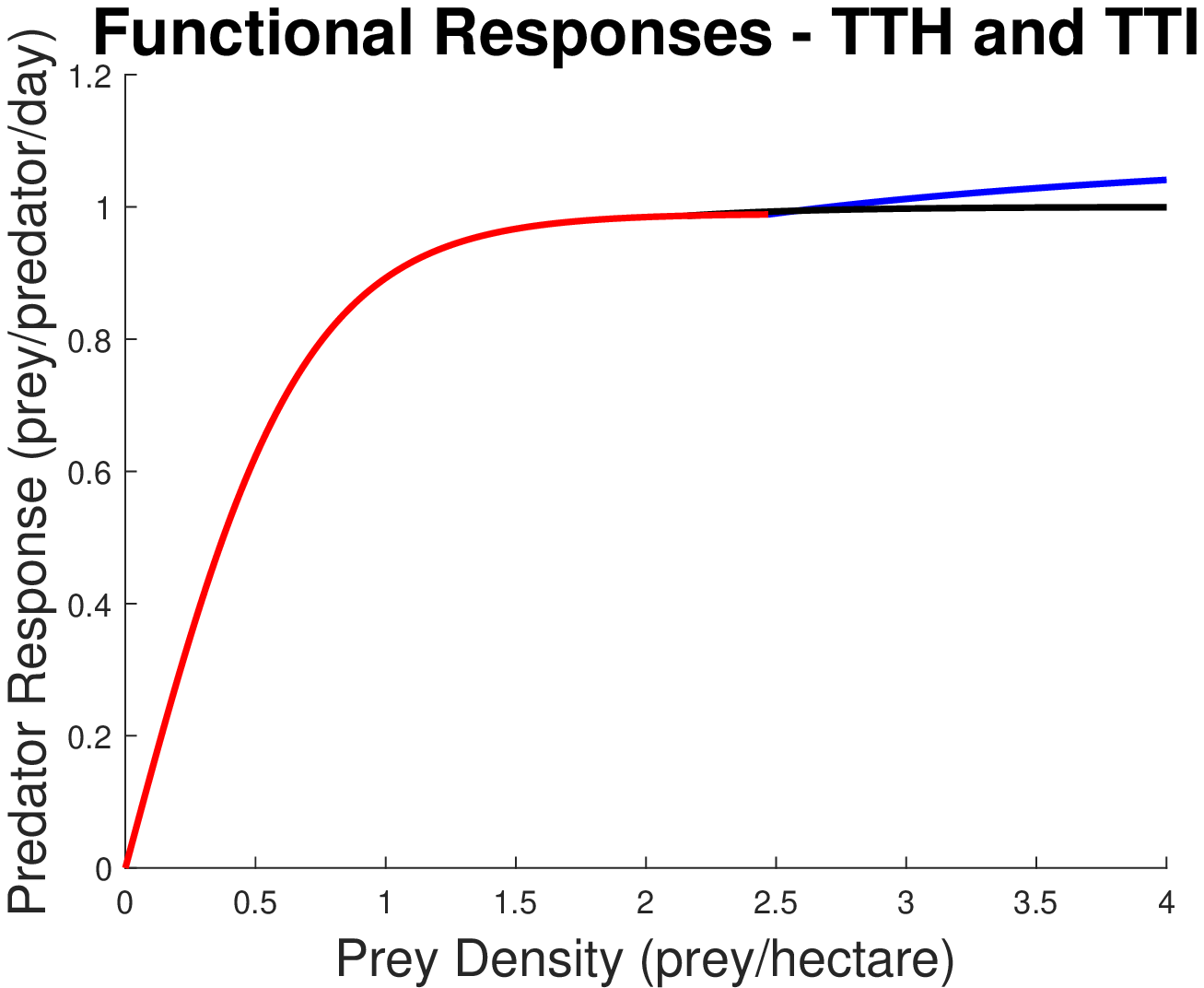}
        \subcaption{Rosenzweig-MacArthur model}
    \end{subfigure}
    \begin{subfigure}{\linewidth}
        \includegraphics[width=0.325\linewidth]{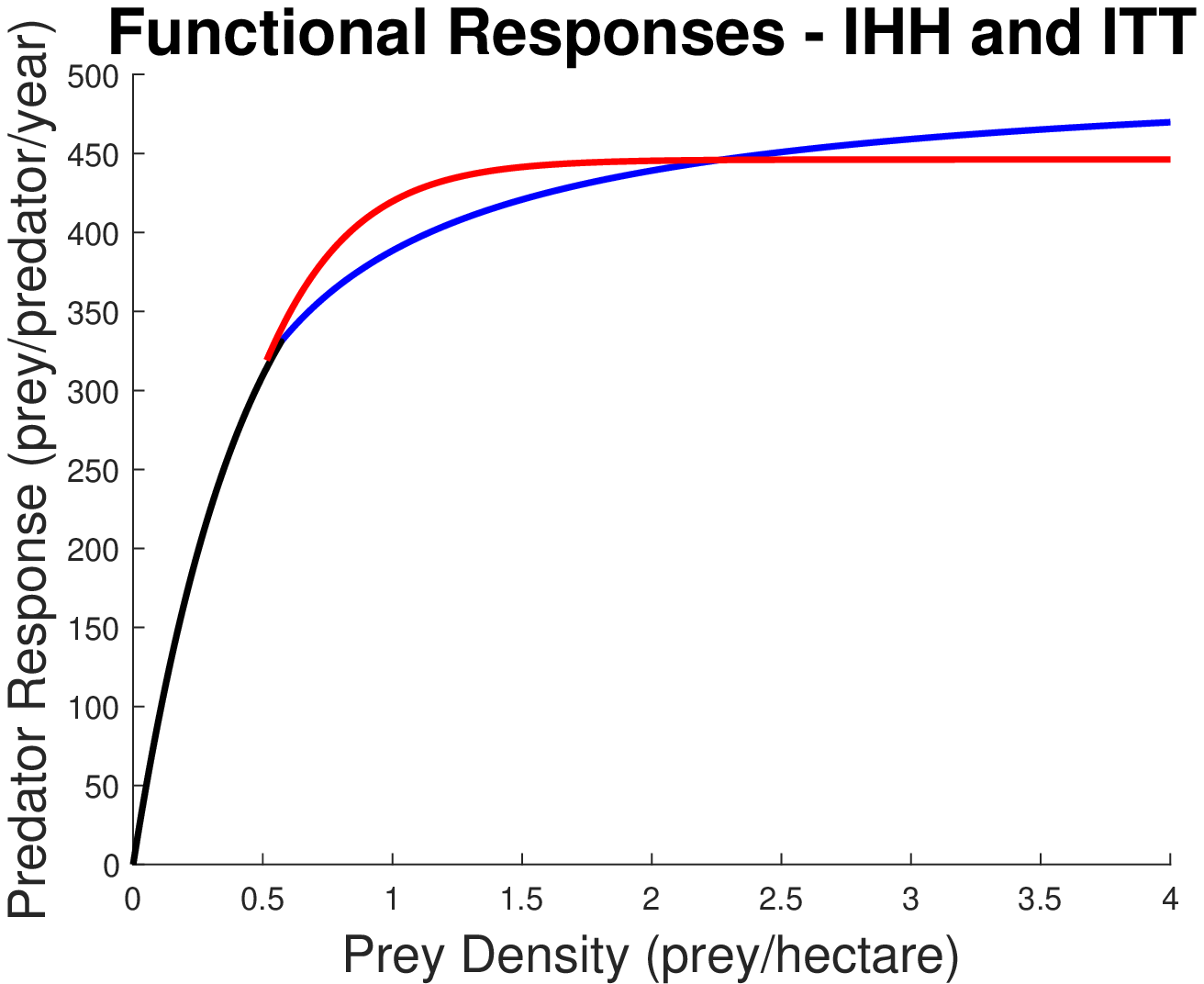}
        \includegraphics[width=0.325\linewidth]{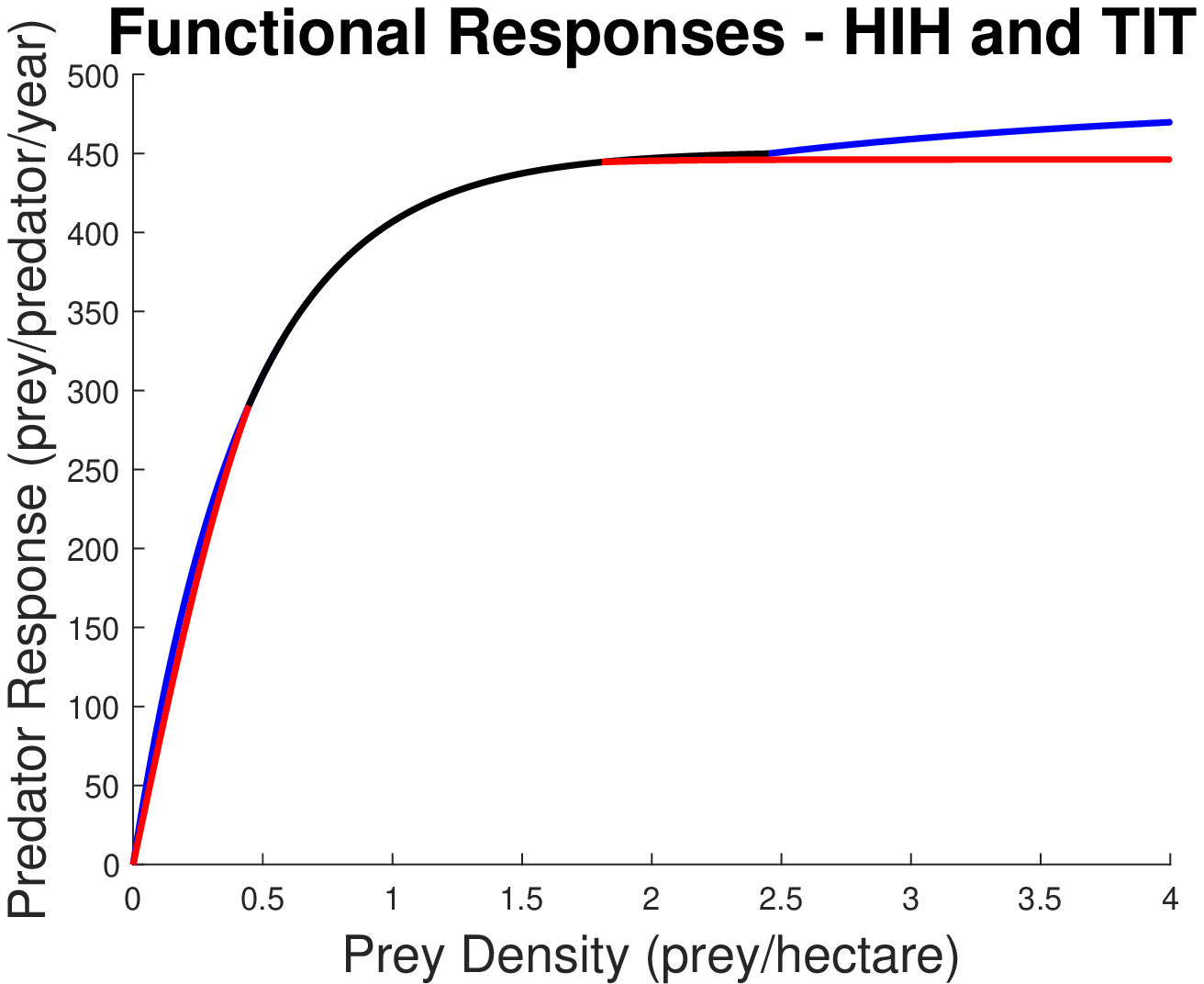}
        \includegraphics[width=0.325\linewidth]{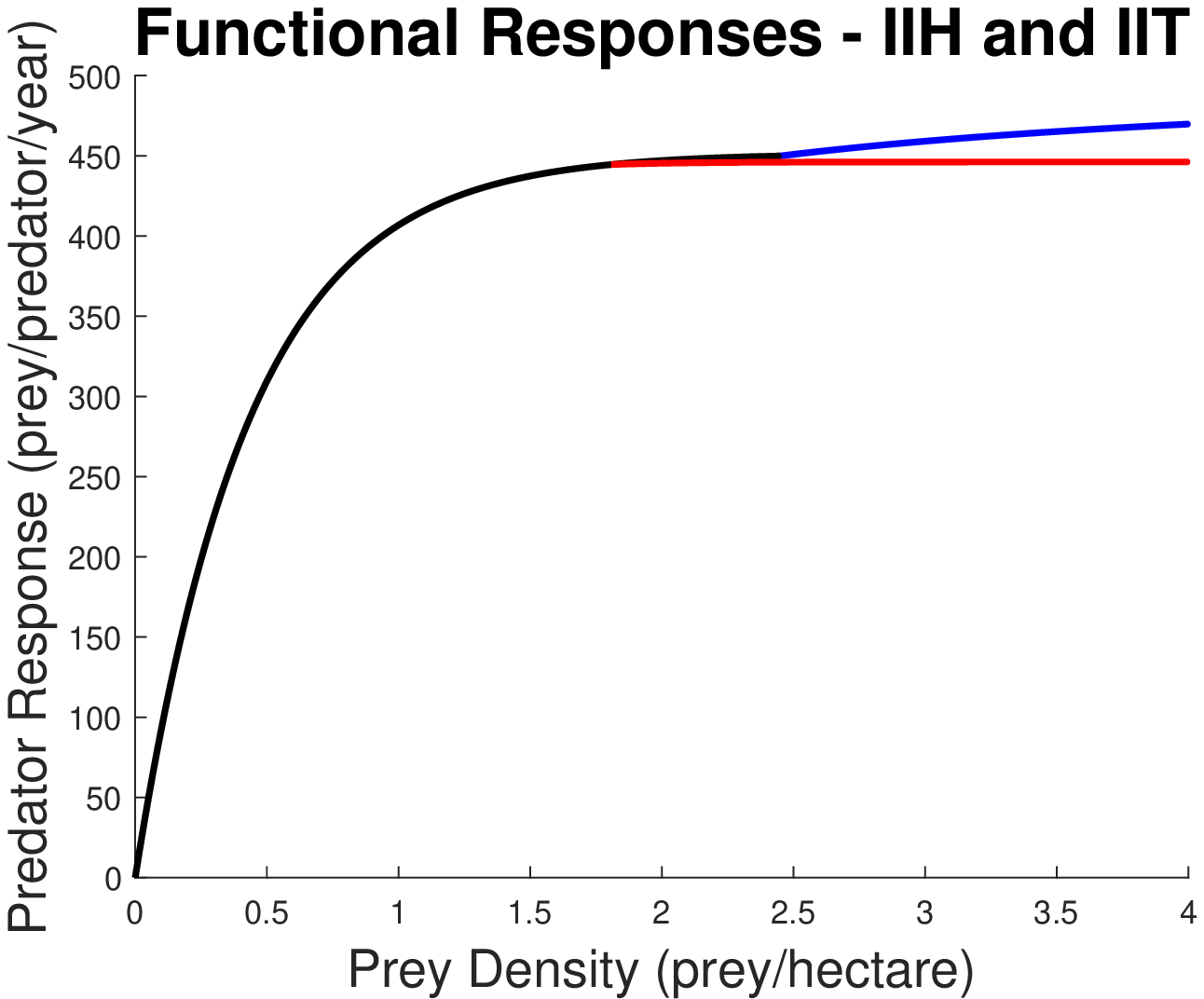}
        \\
        \includegraphics[width=0.325\linewidth]{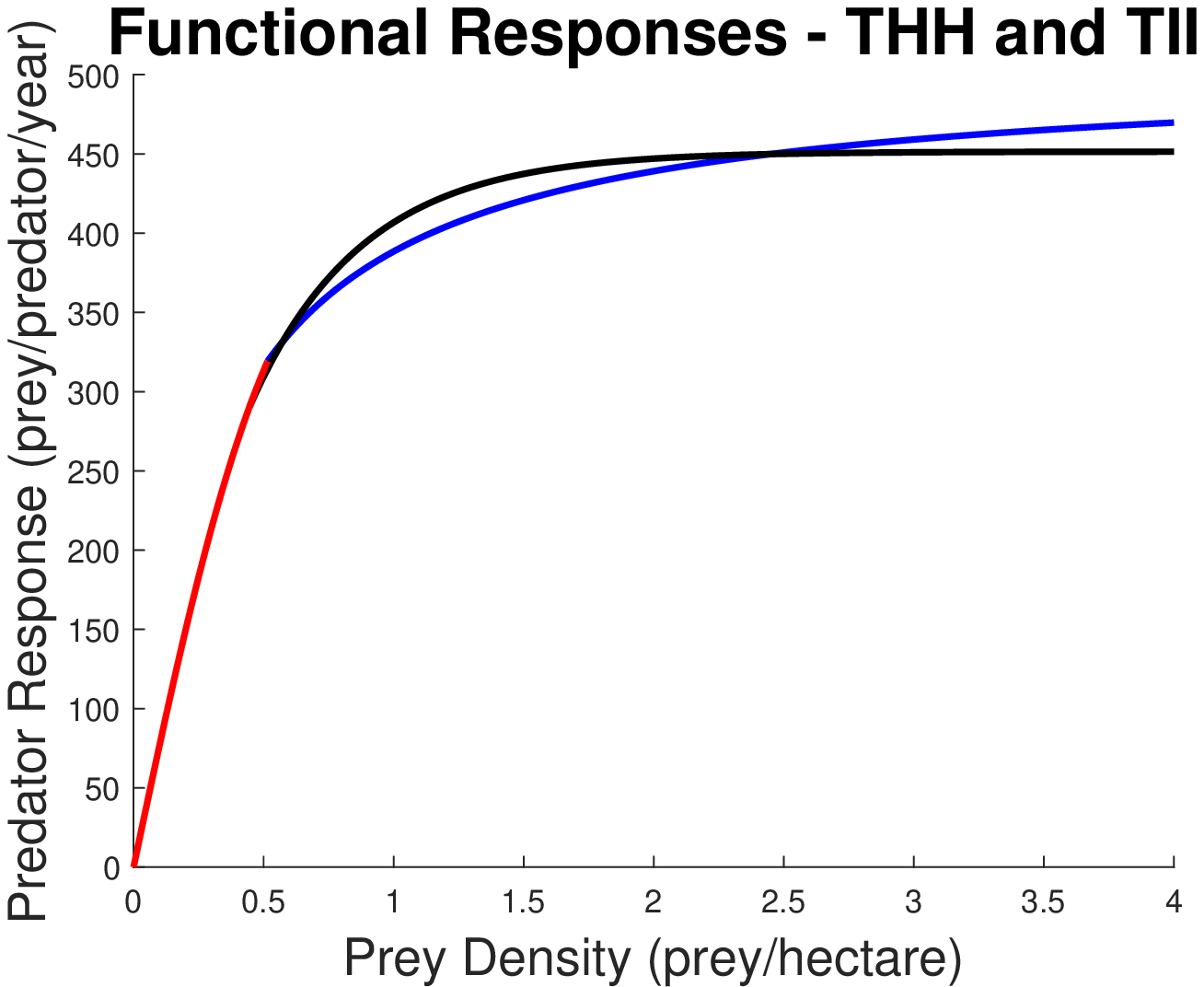}
        \includegraphics[width=0.325\linewidth]{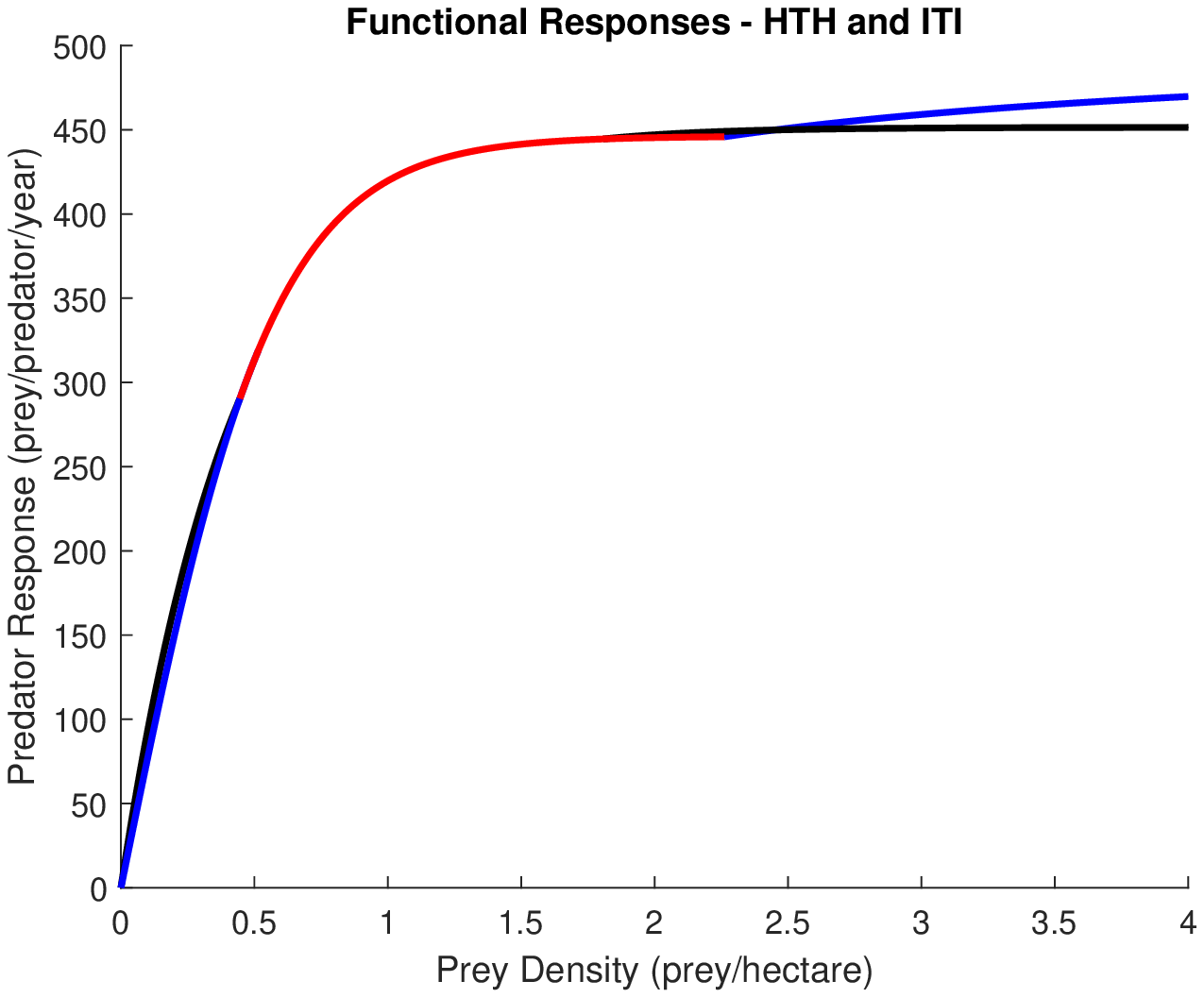}
        \includegraphics[width=0.325\linewidth]{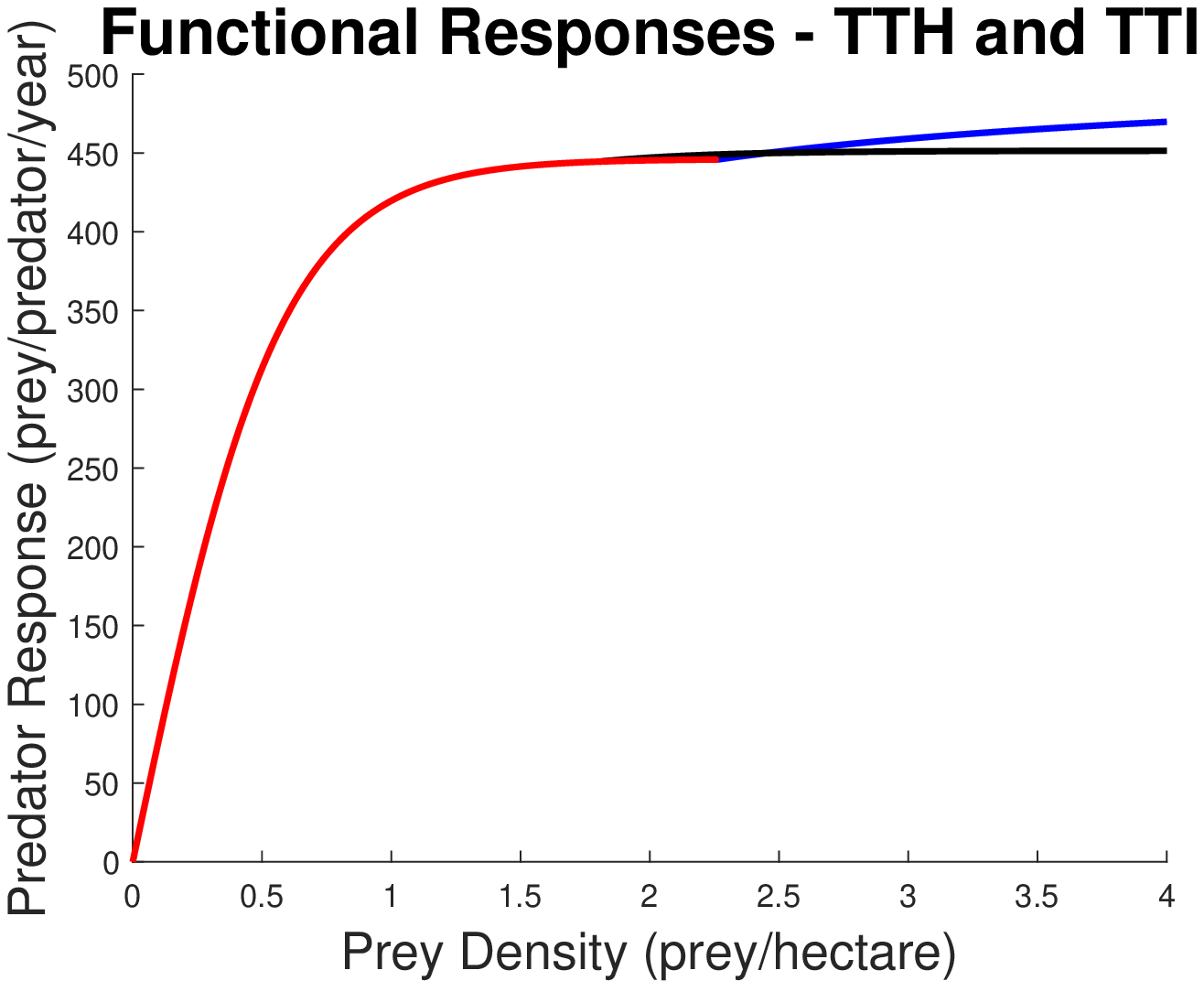}
        \subcaption{Leslie-Gower-May model}
    \end{subfigure}
    \includegraphics[width=0.5\linewidth]{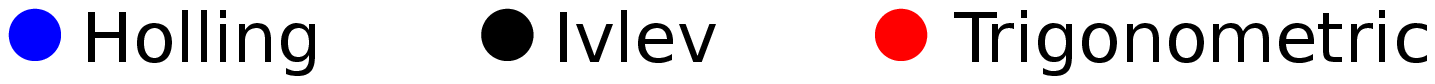}
    \caption{Functional responses for the Rosenzweig-MacArthur model defined piecewise from two functional responses and switching at the corresponding intersection. Parameters and intersection values are as given in the manuscript.}
    \label{supfig:SwitchingCurves}
\end{figure}

%---------------------------------------------------------
Figure \ref{supfig:SwitchingCurves} gives the piecewise identical functional responses.

\subsection{Bifurcation Diagrams}
%---------------------------------------------------------

\begin{figure}[ht]
    \centering
    \includegraphics[width=0.325\linewidth]{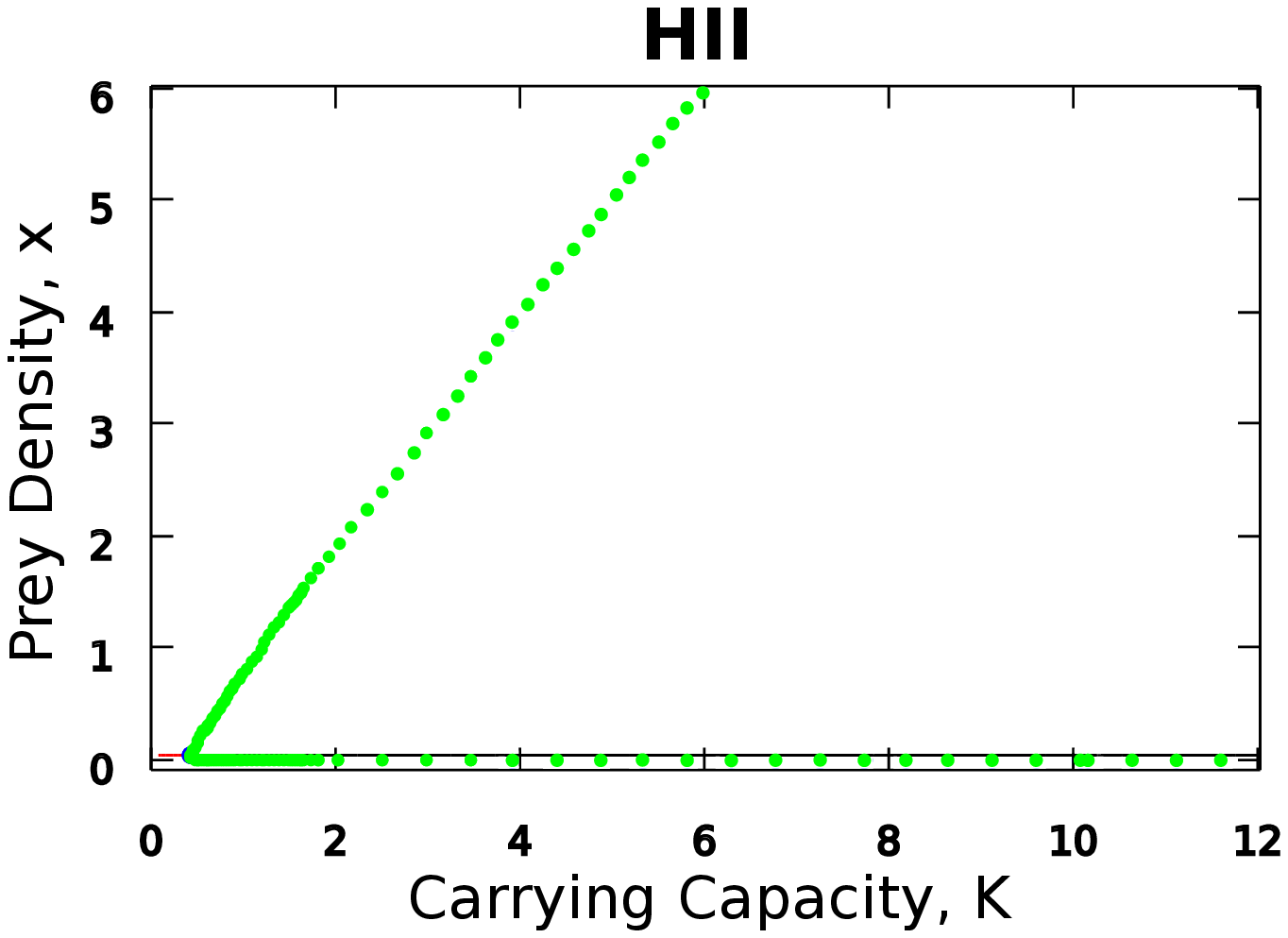}
    \includegraphics[width=0.325\linewidth]{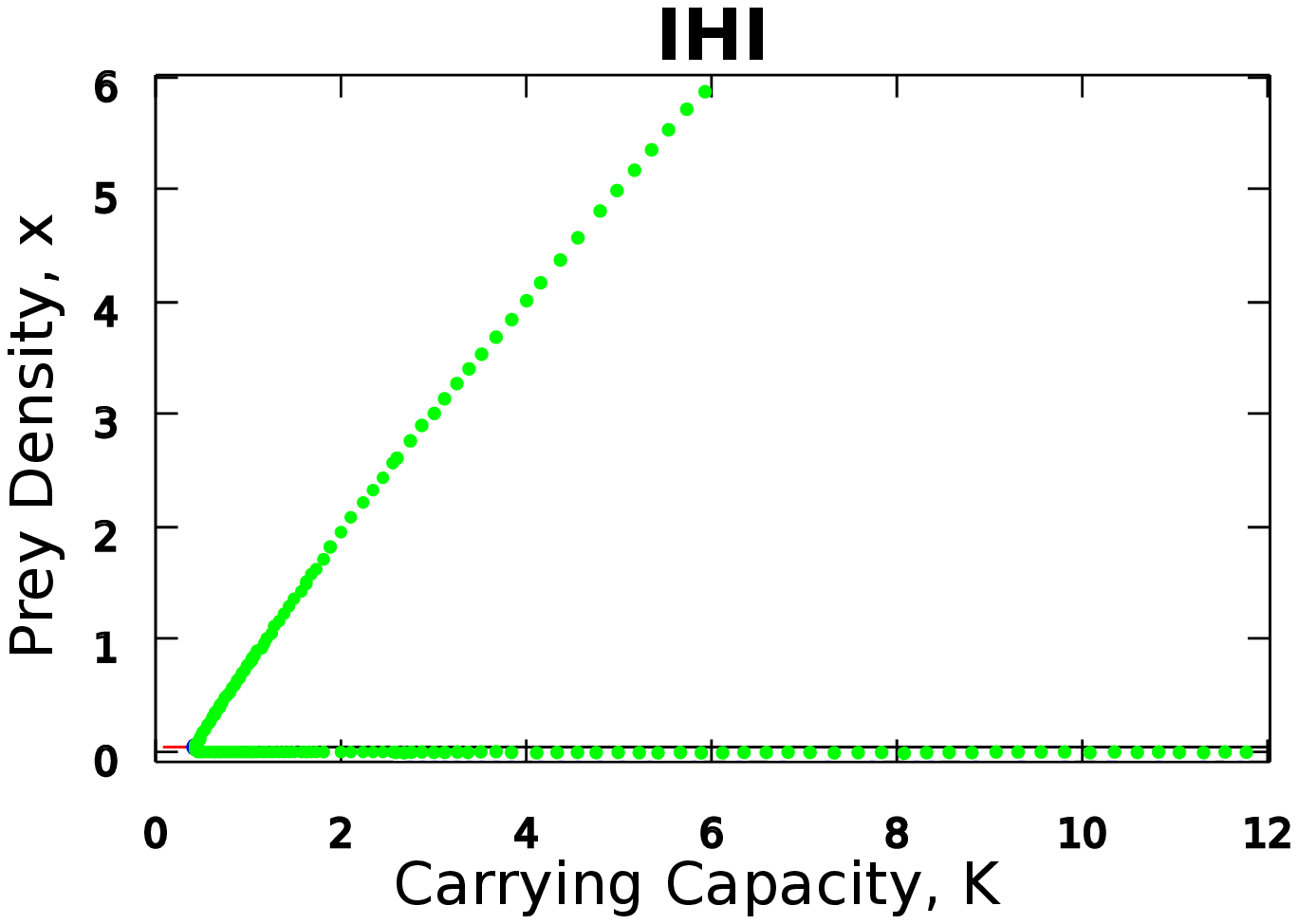}
    \includegraphics[width=0.325\linewidth]{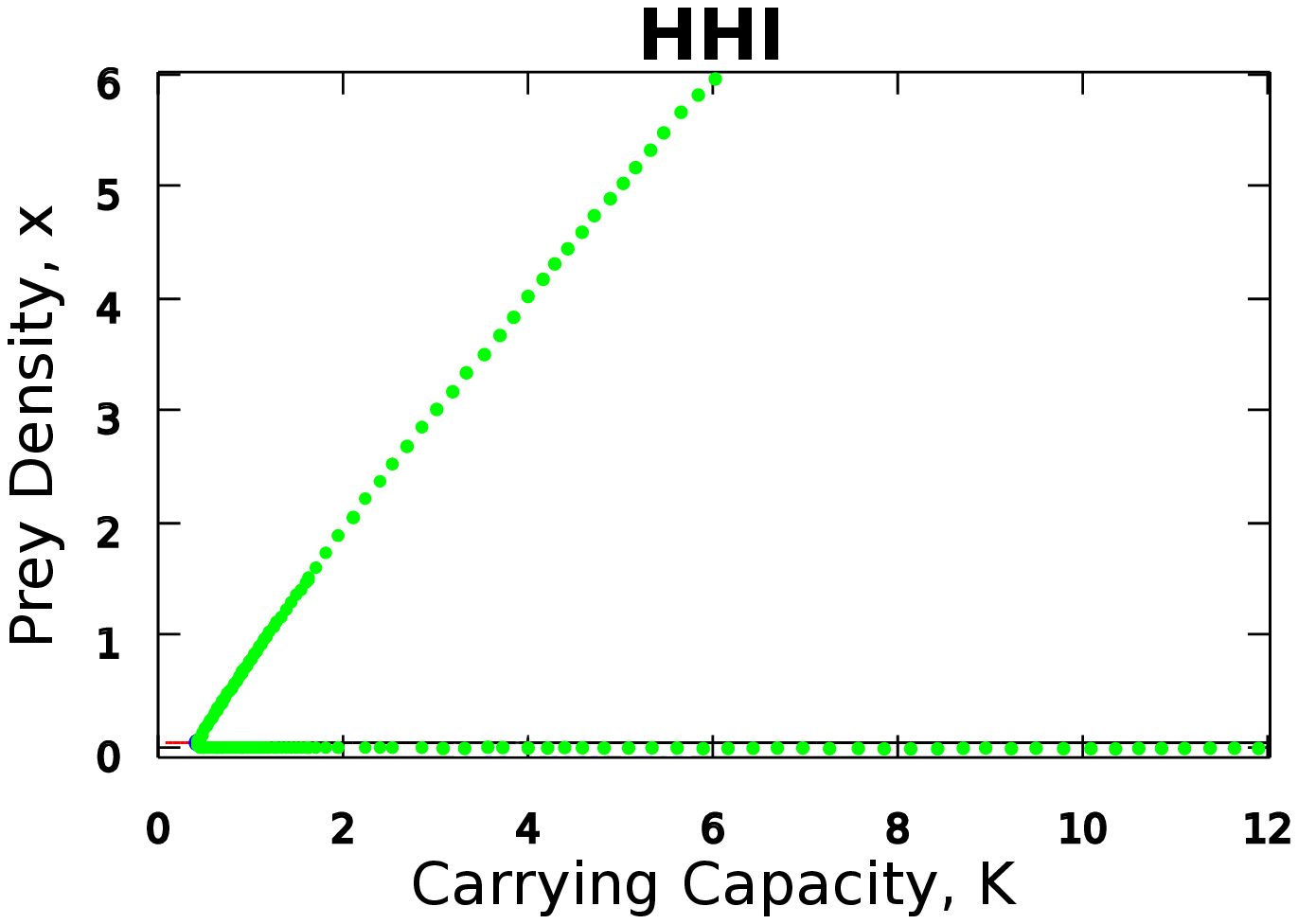}
    \\
    
    \includegraphics[width=0.325\linewidth]{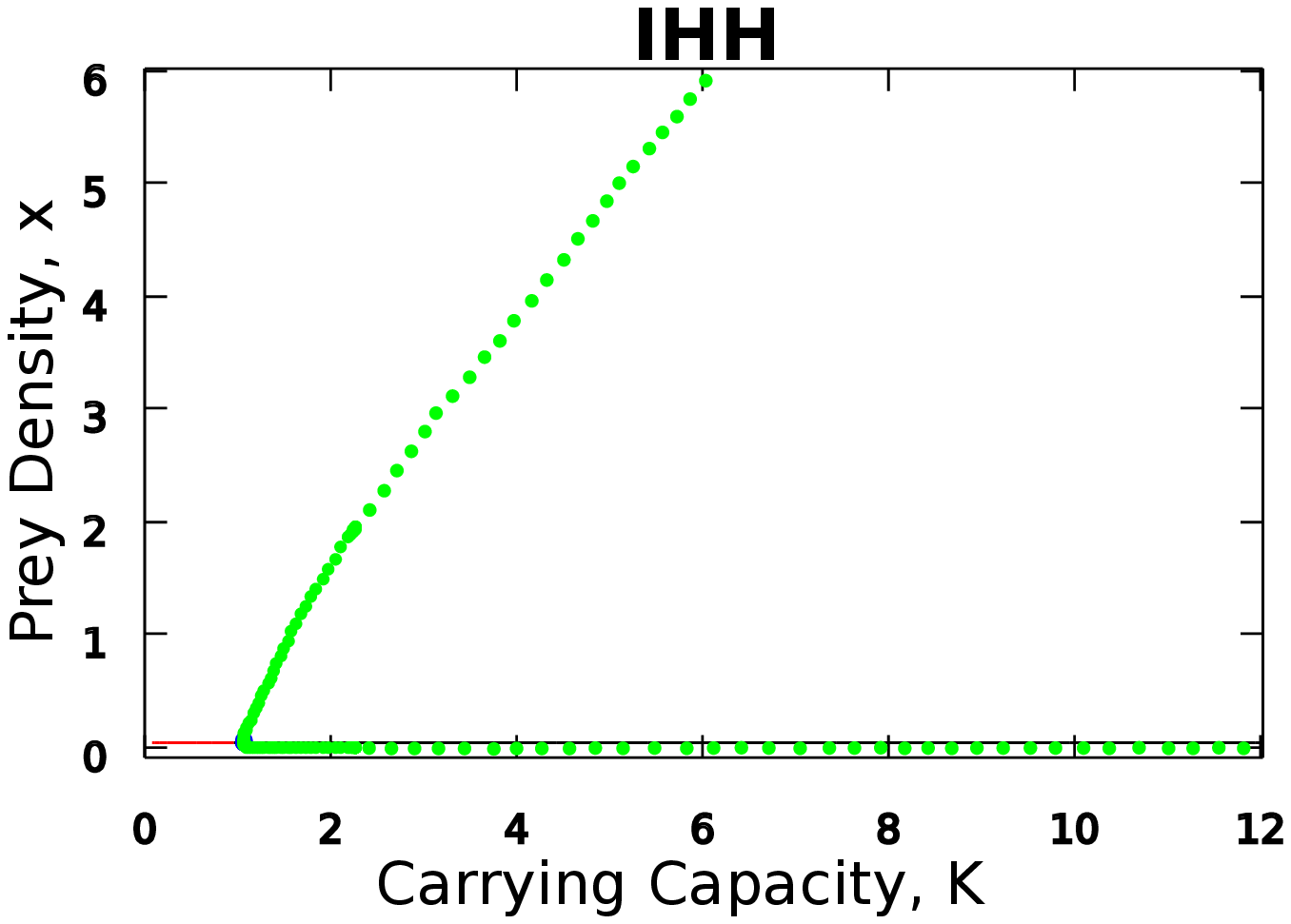}
    \includegraphics[width=0.325\linewidth]{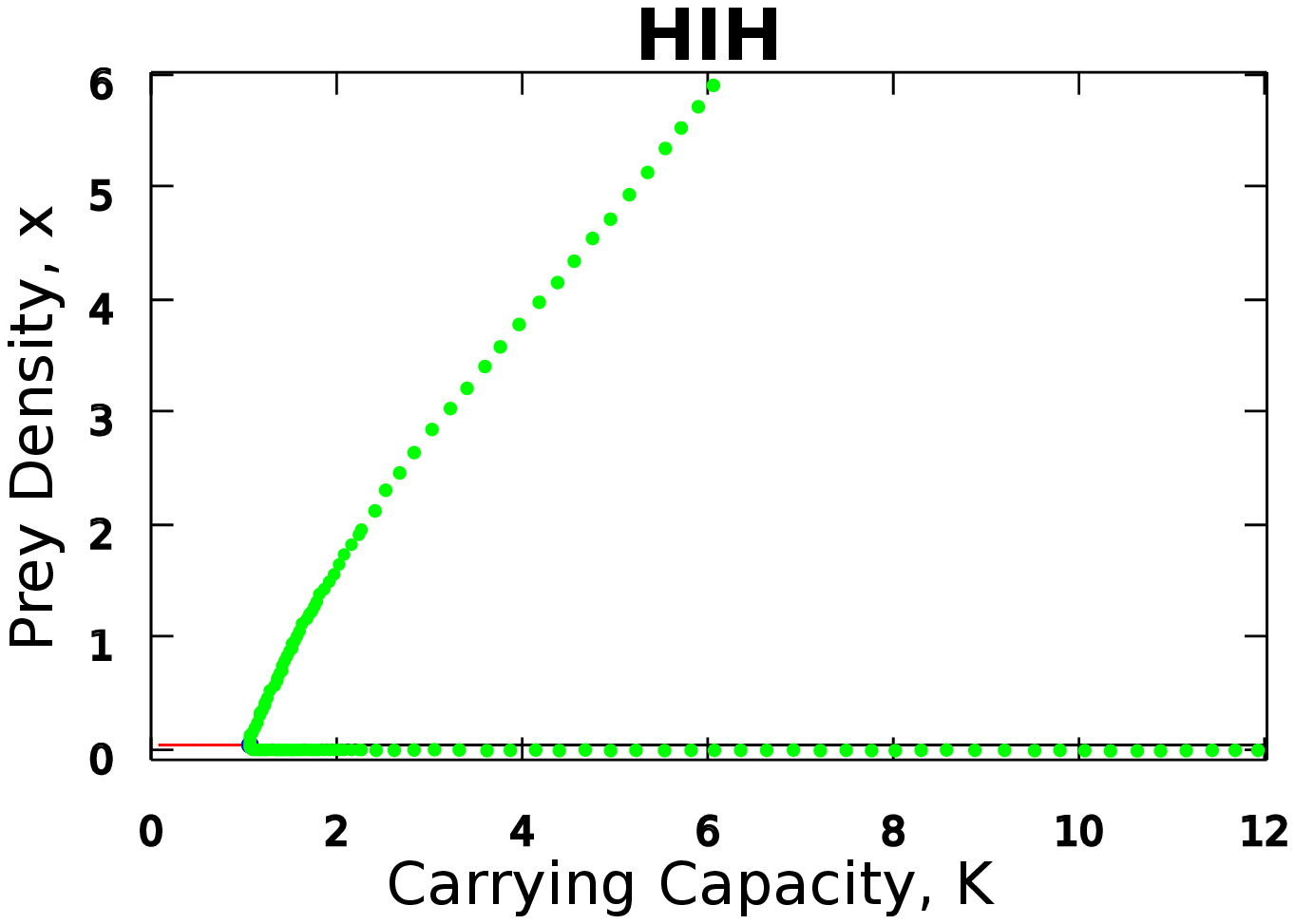}
    \includegraphics[width=0.325\linewidth]{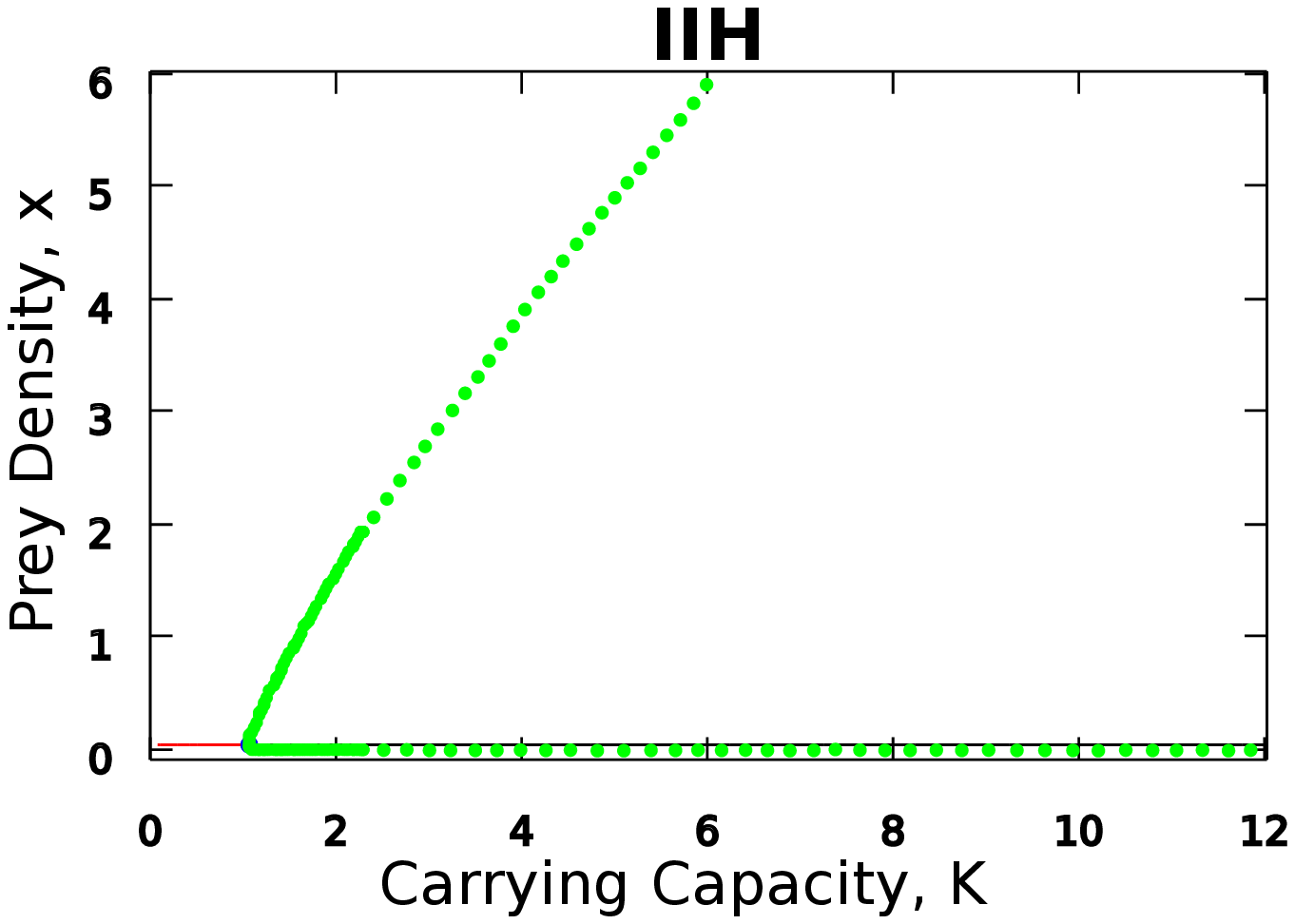}
    \\
    \includegraphics[width=0.325\linewidth]{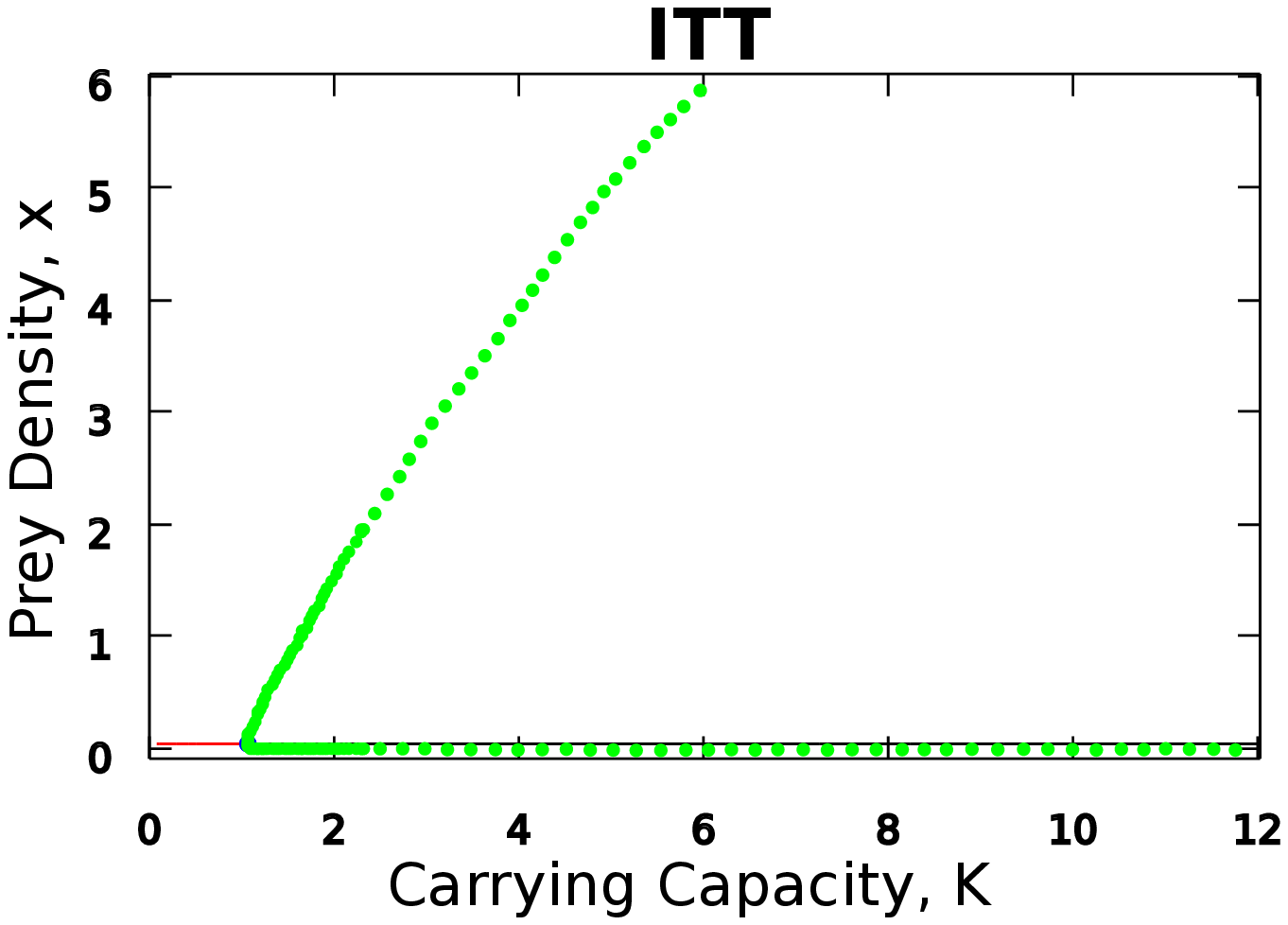}
    \includegraphics[width=0.325\linewidth]{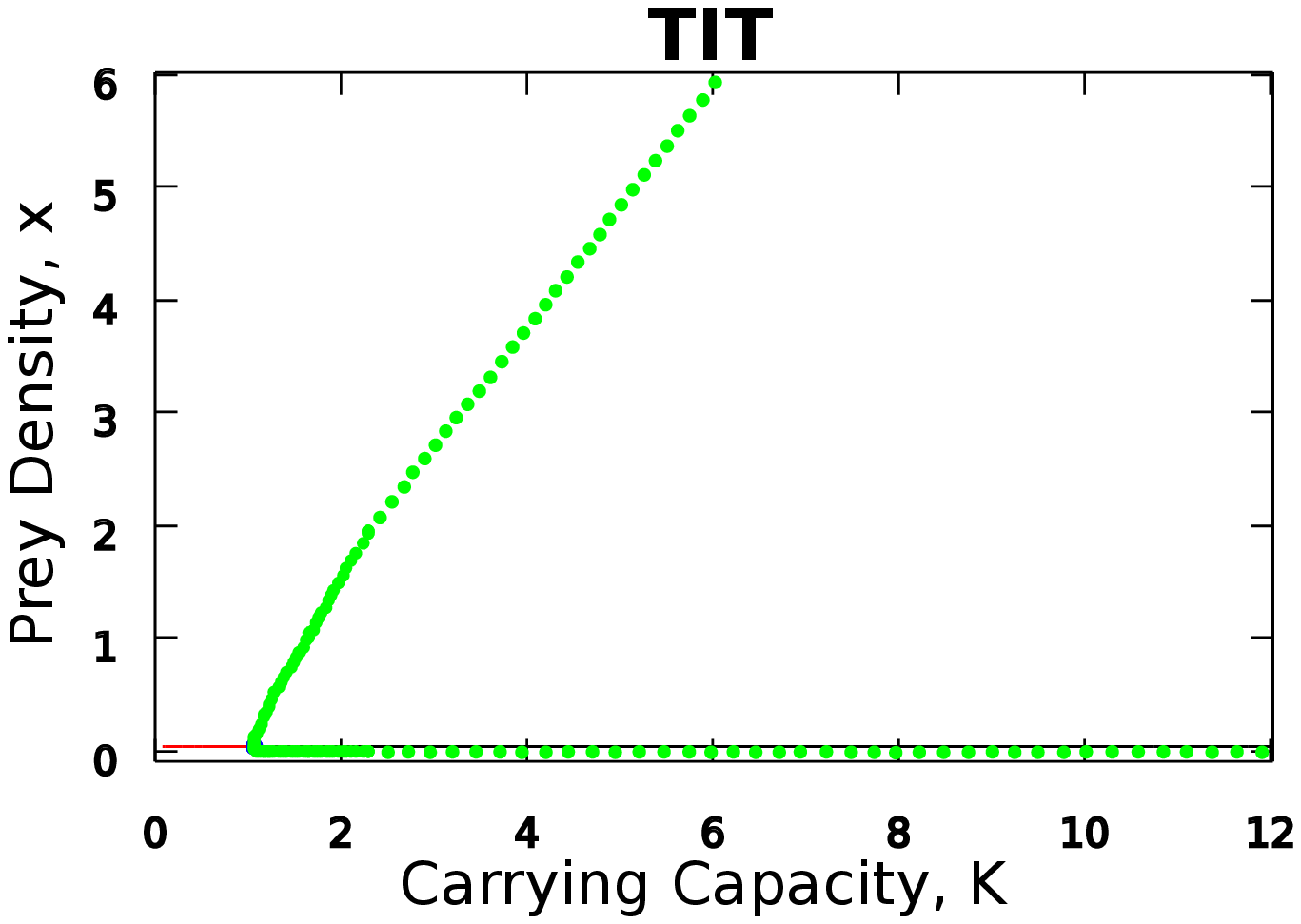}
    \includegraphics[width=0.325\linewidth]{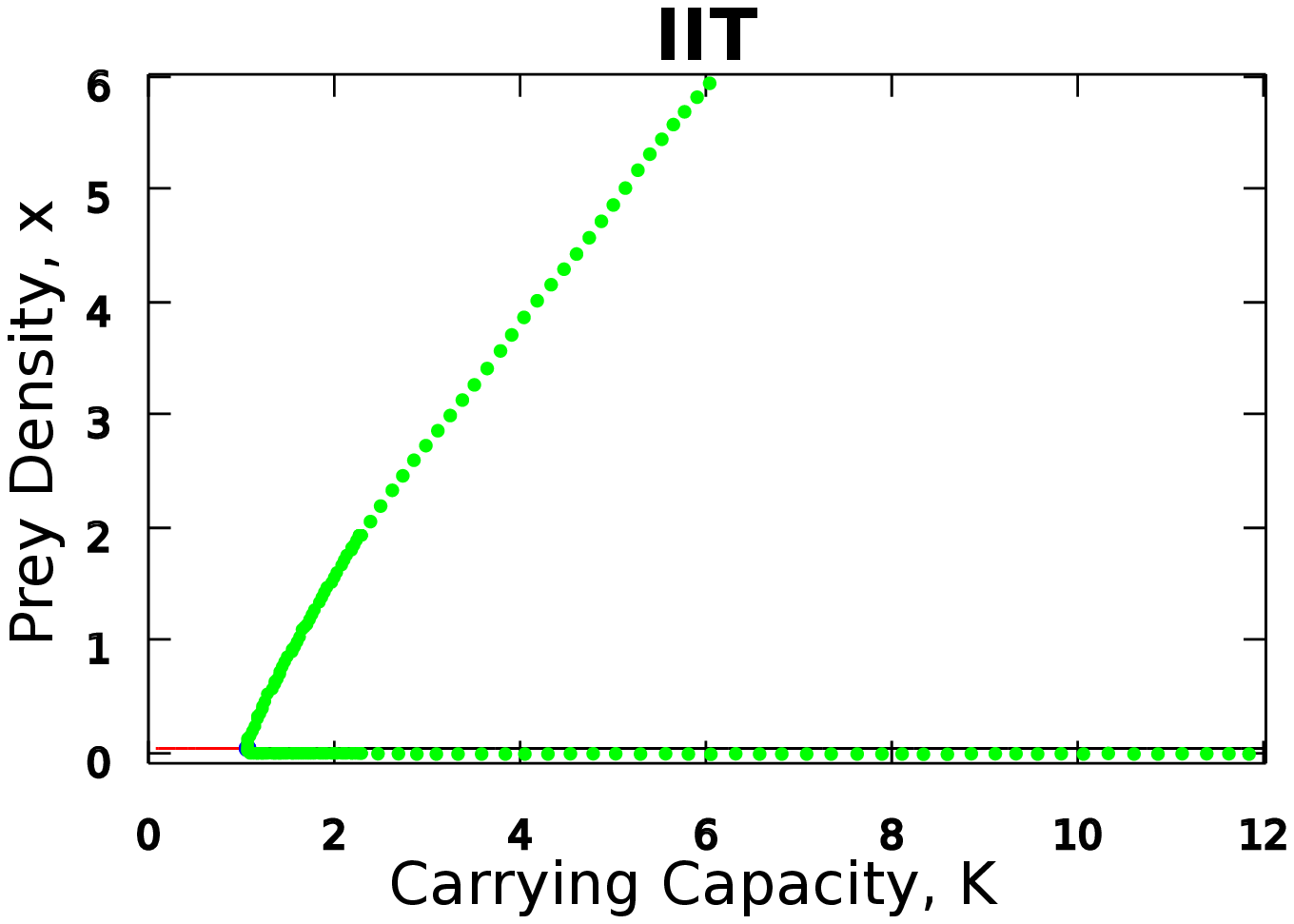}    
    \\
    
    \includegraphics[width=0.325\linewidth]{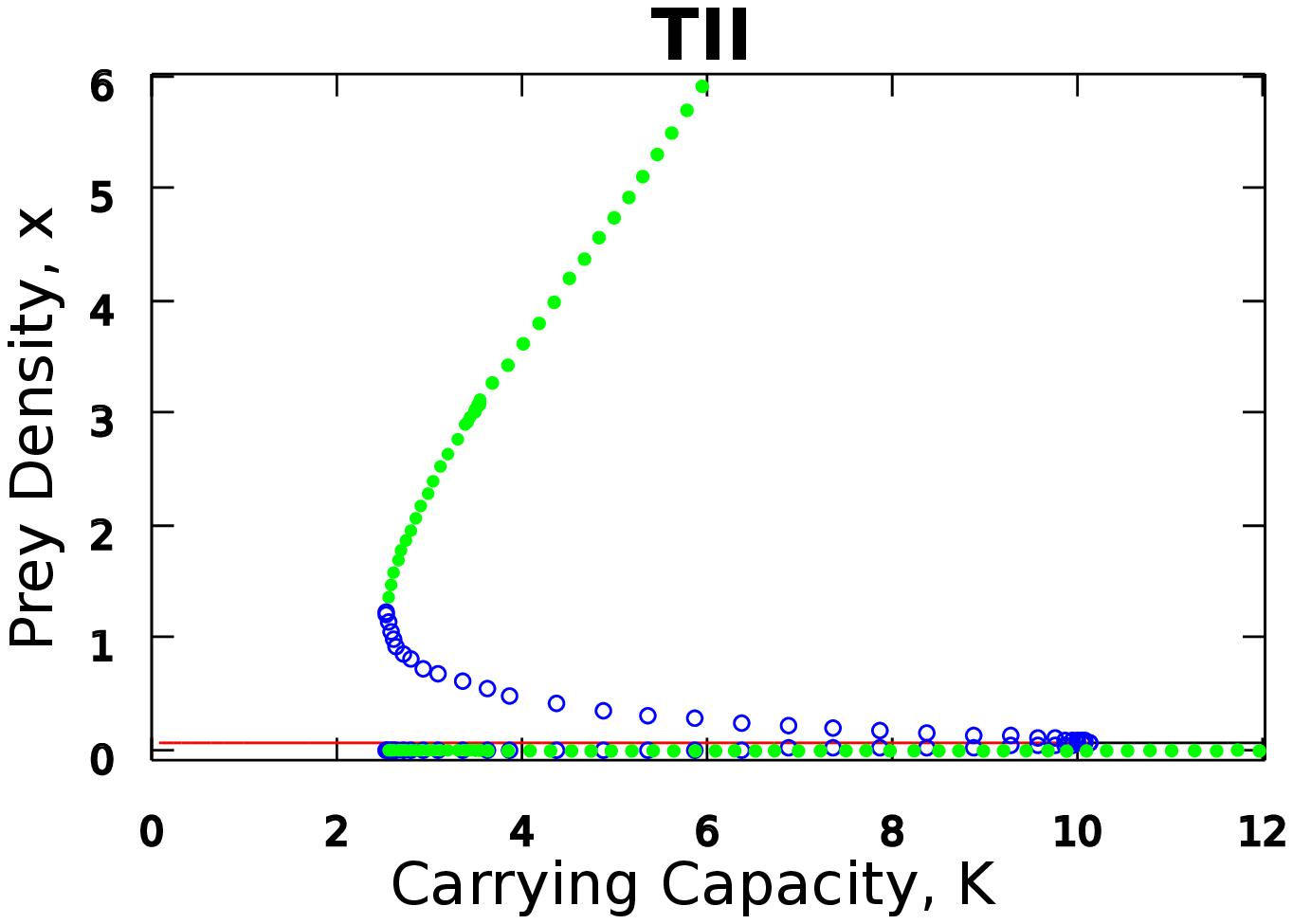}
    \includegraphics[width=0.325\linewidth]{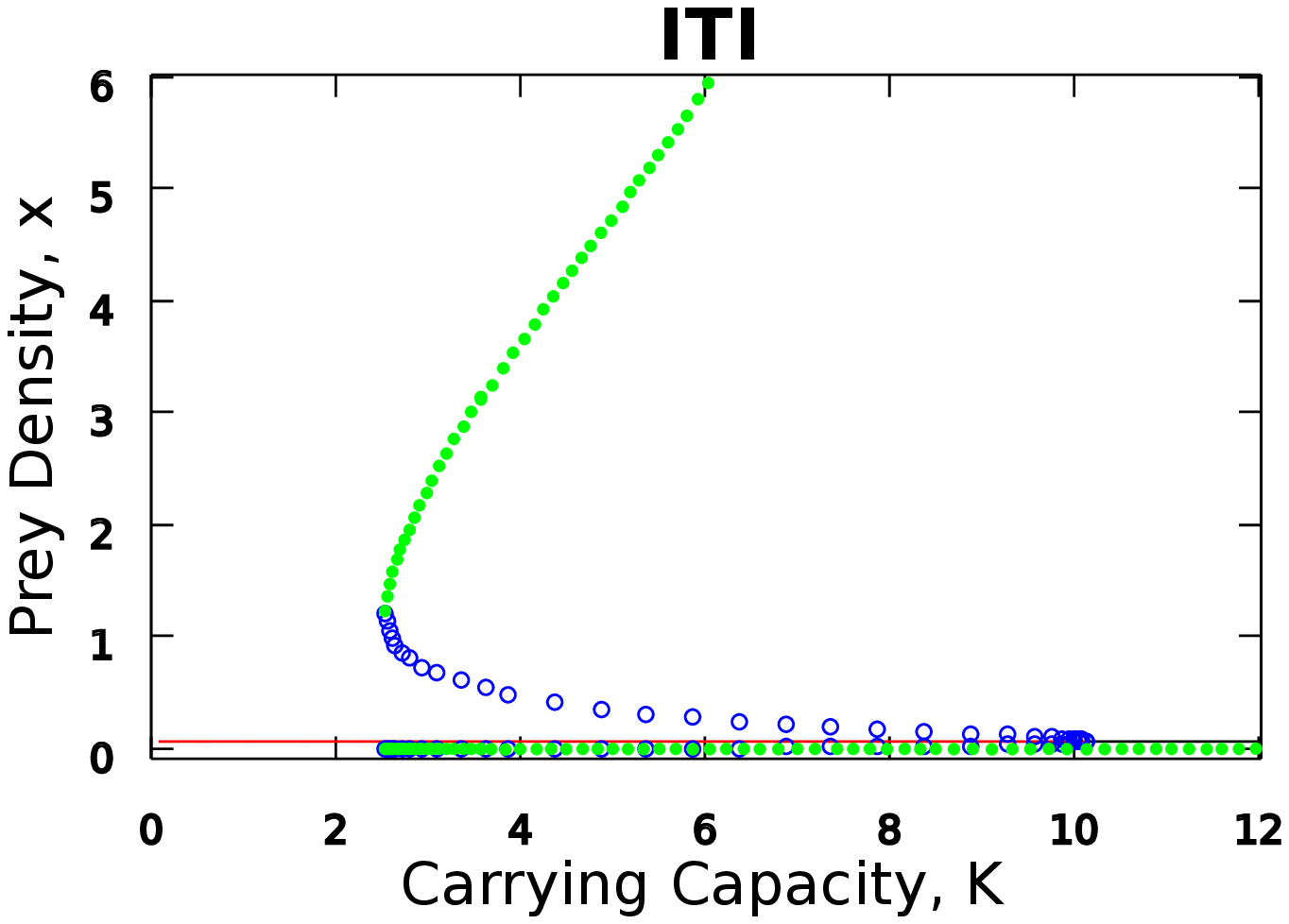}
    \includegraphics[width=0.325\linewidth]{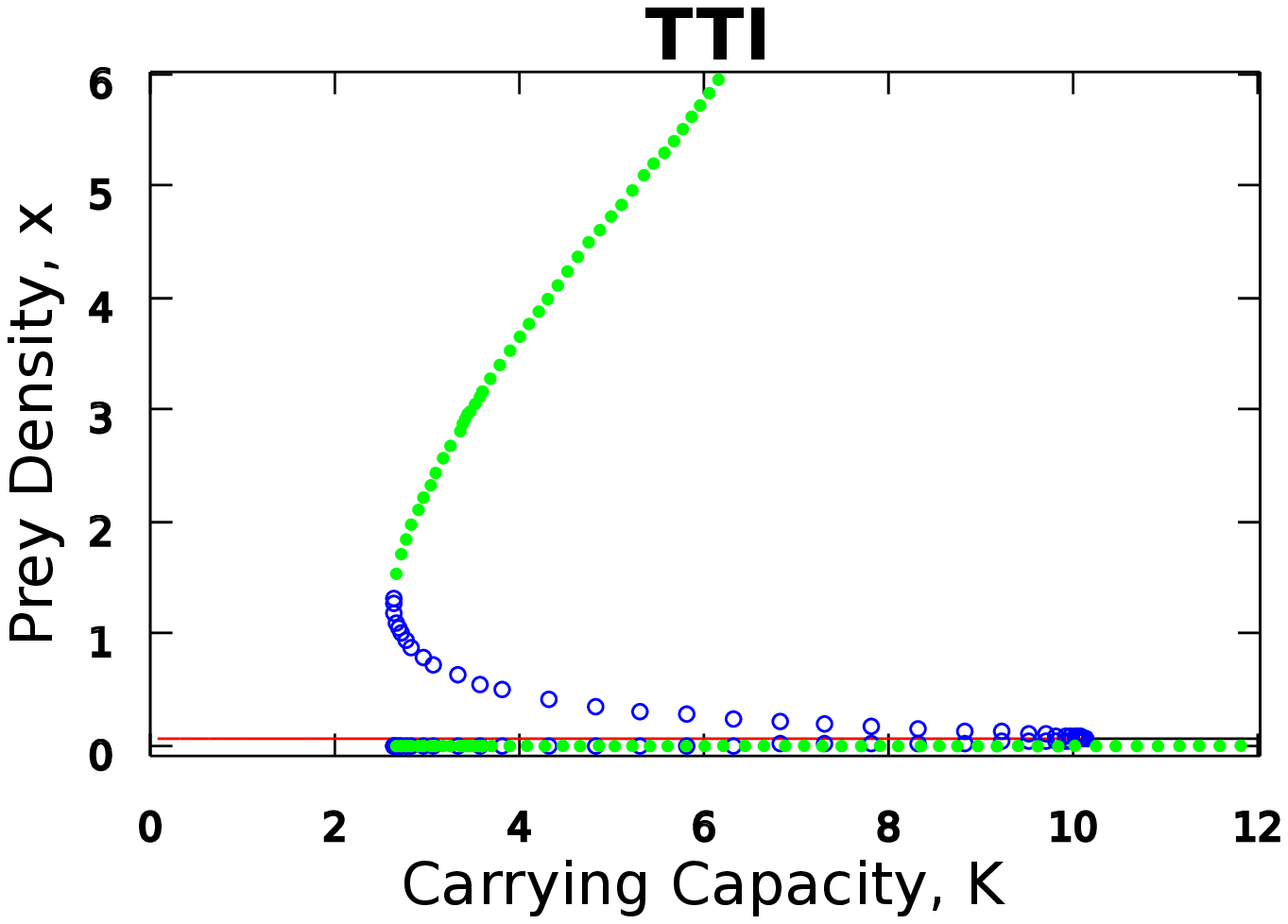}
    \\
    
    \includegraphics[width=\linewidth]{Pictures/FIGLeg2.eps}
    \caption{Bifurcation diagrams for the Rosenzweig-MacArthur model created using functional responses defined piecewise from two functional responses and switching at the corresponding intersections. Values for the intersections are given in the manuscript and bifurcation values are given in Table \ref{suptable:BifurValues}.}
    \label{supfig:SwitchingBifurRM}
\end{figure}

\begin{figure}[ht]
    \centering
    \includegraphics[width=0.325\linewidth]{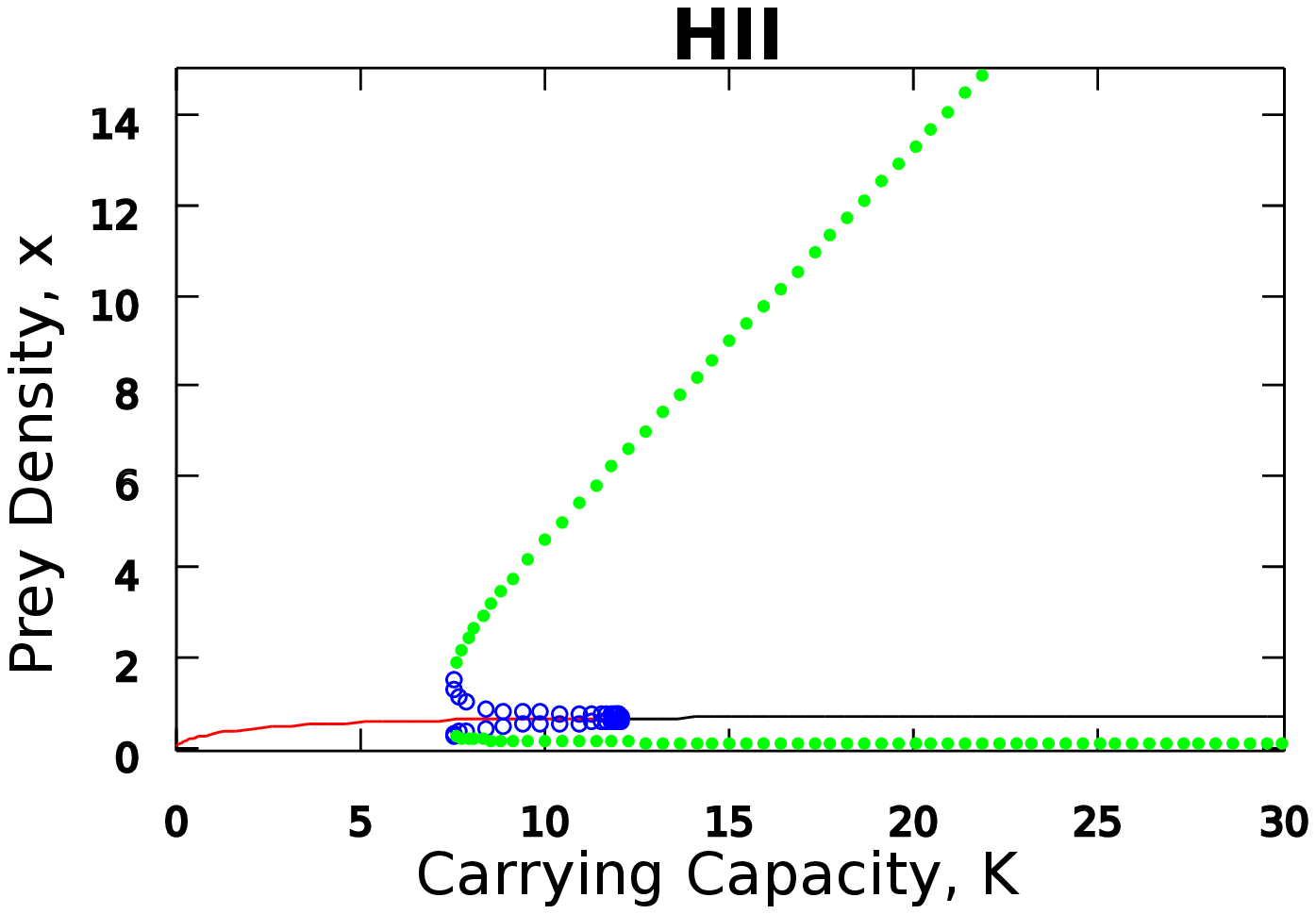}
    \includegraphics[width=0.325\linewidth]{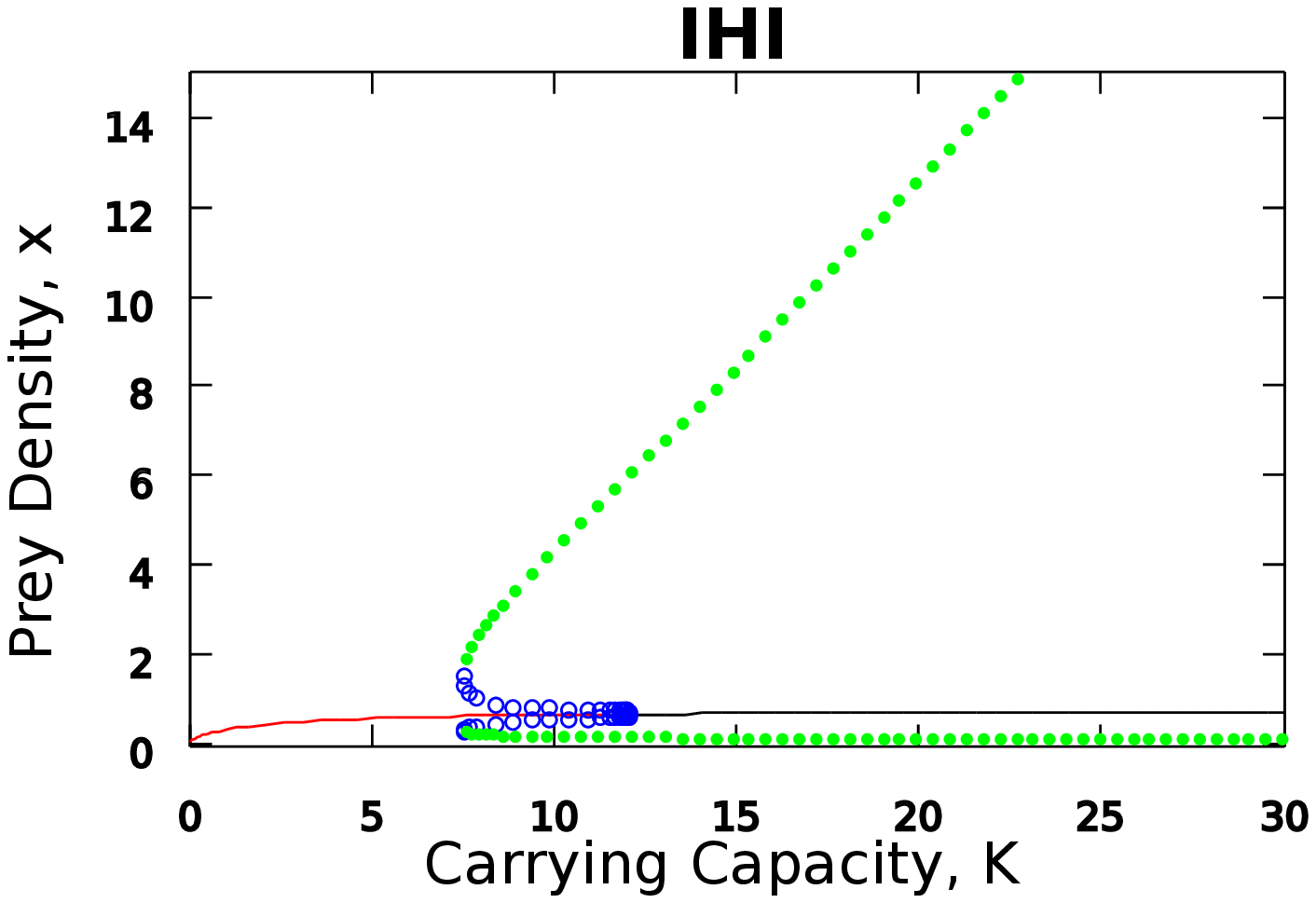}
    \includegraphics[width=0.325\linewidth]{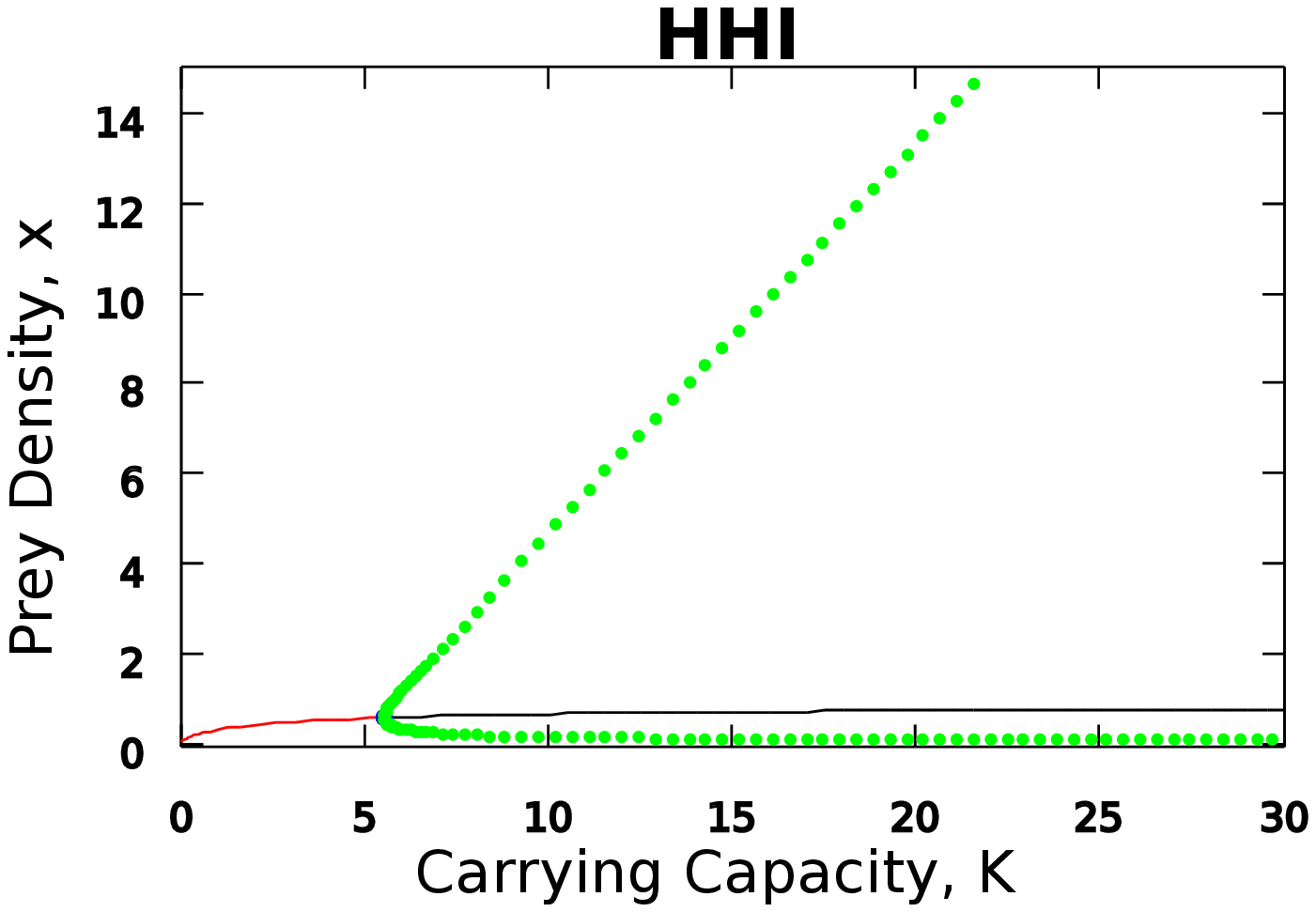}
    \\
      
    \includegraphics[width=0.325\linewidth]{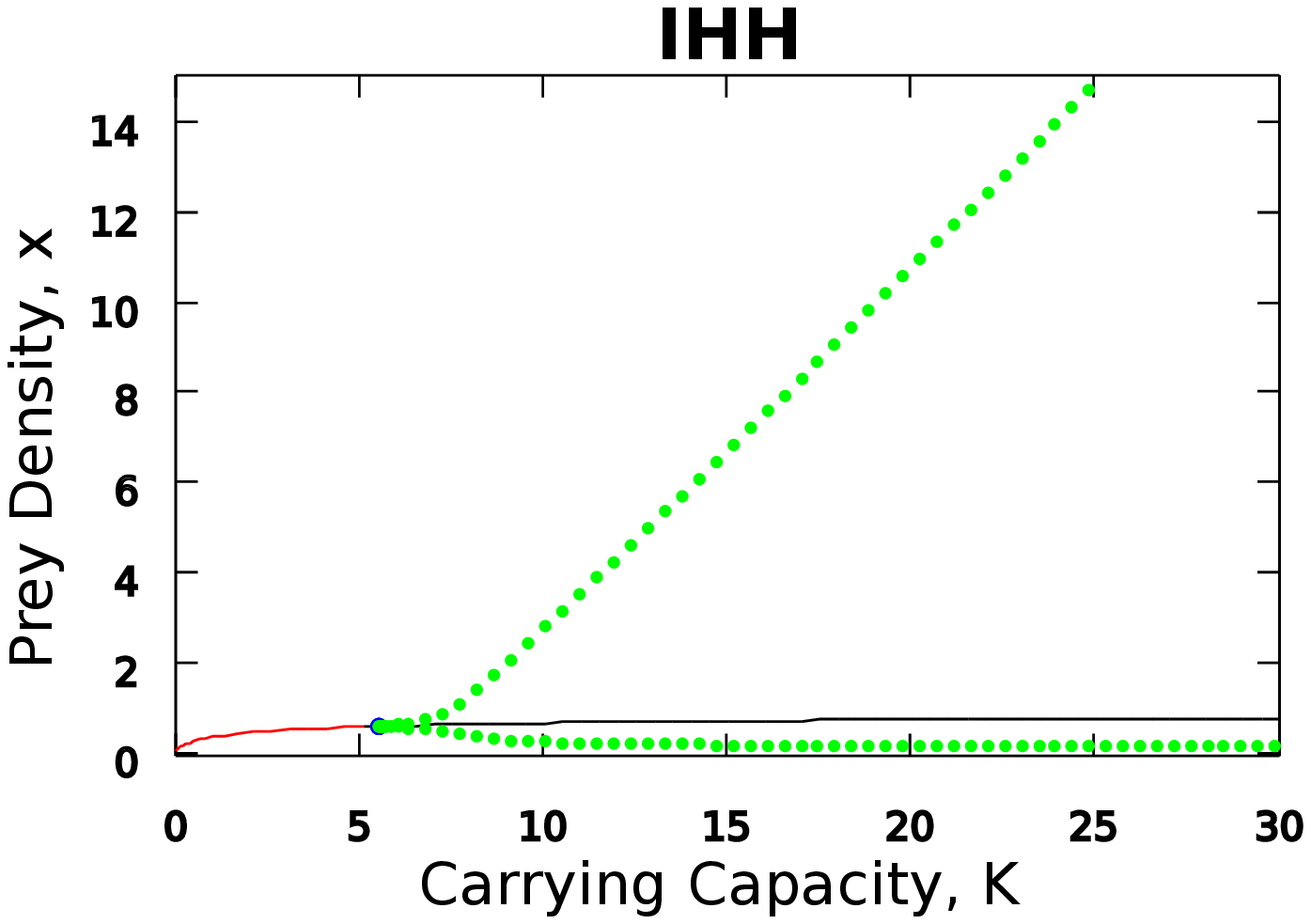}
    \includegraphics[width=0.325\linewidth]{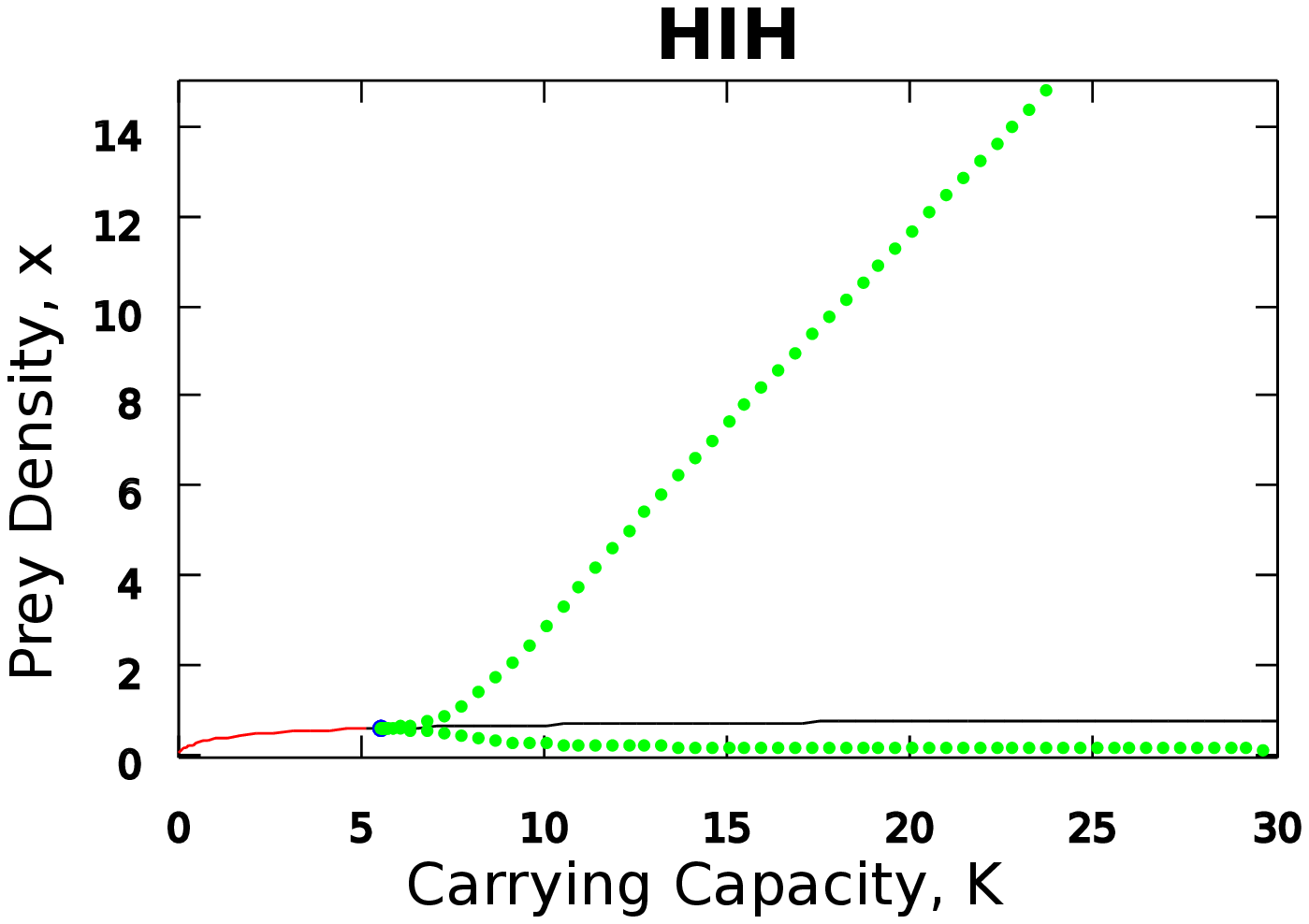}
    \includegraphics[width=0.325\linewidth]{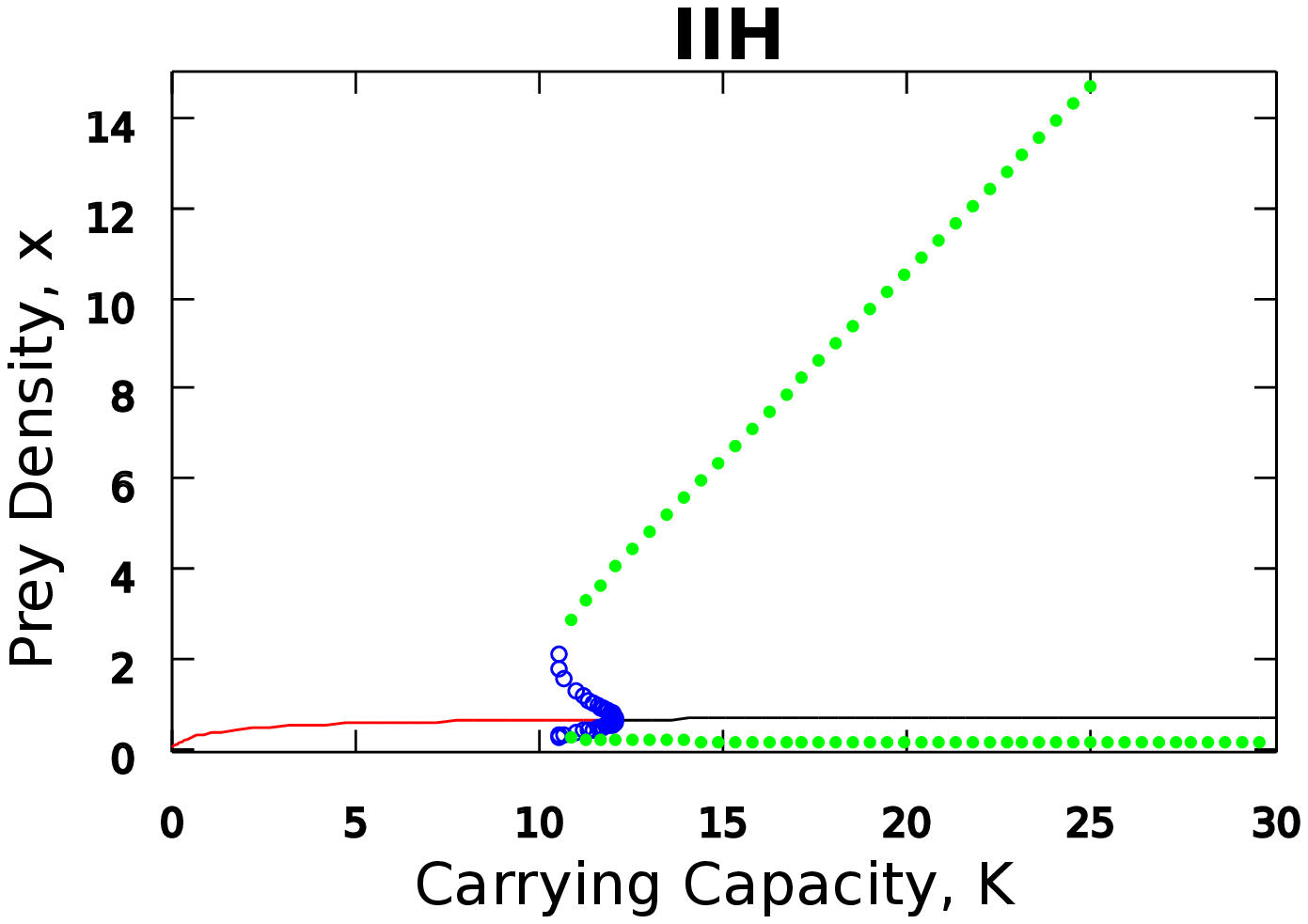}
    \\
    \includegraphics[width=0.325\linewidth]{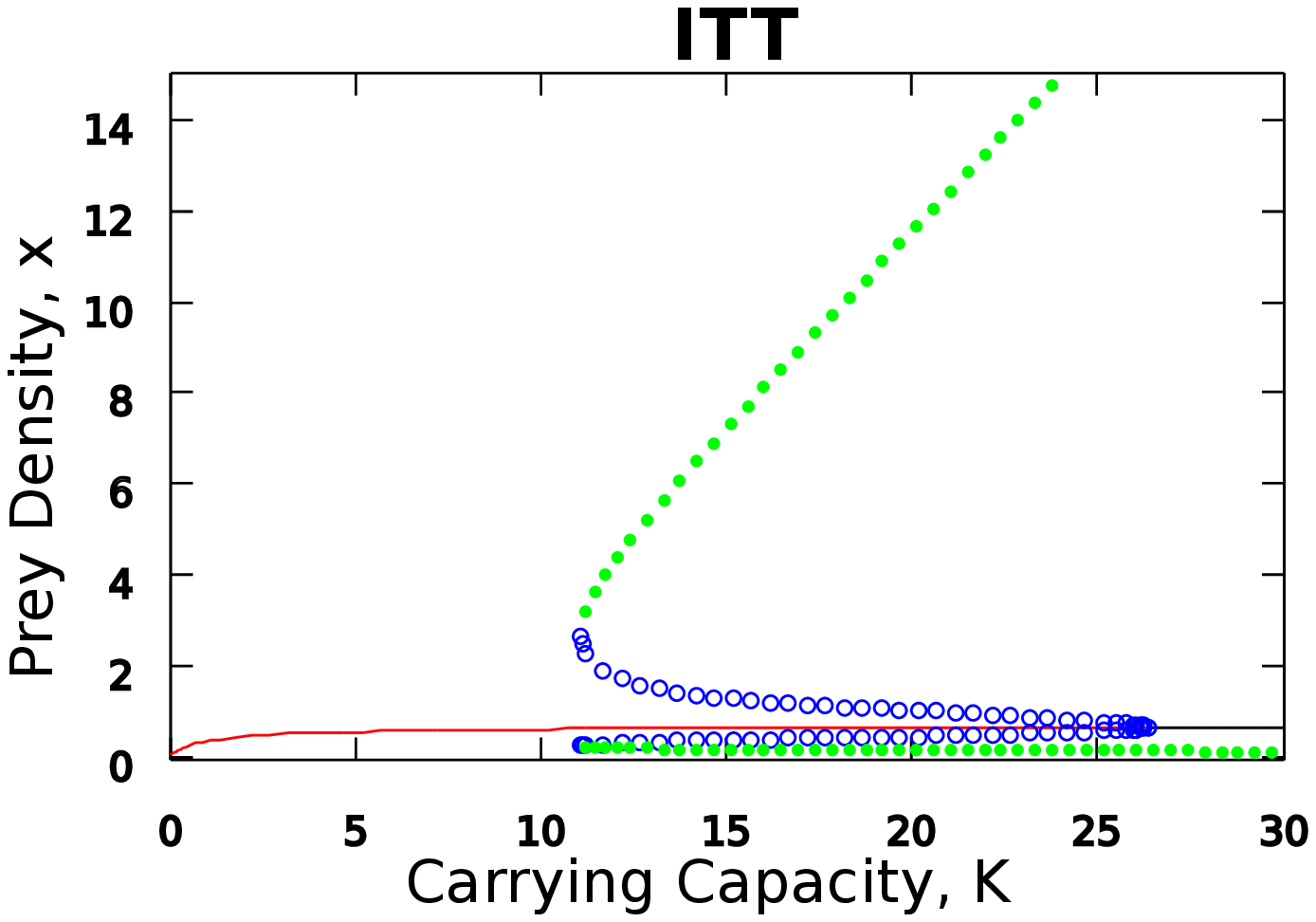}
    \includegraphics[width=0.325\linewidth]{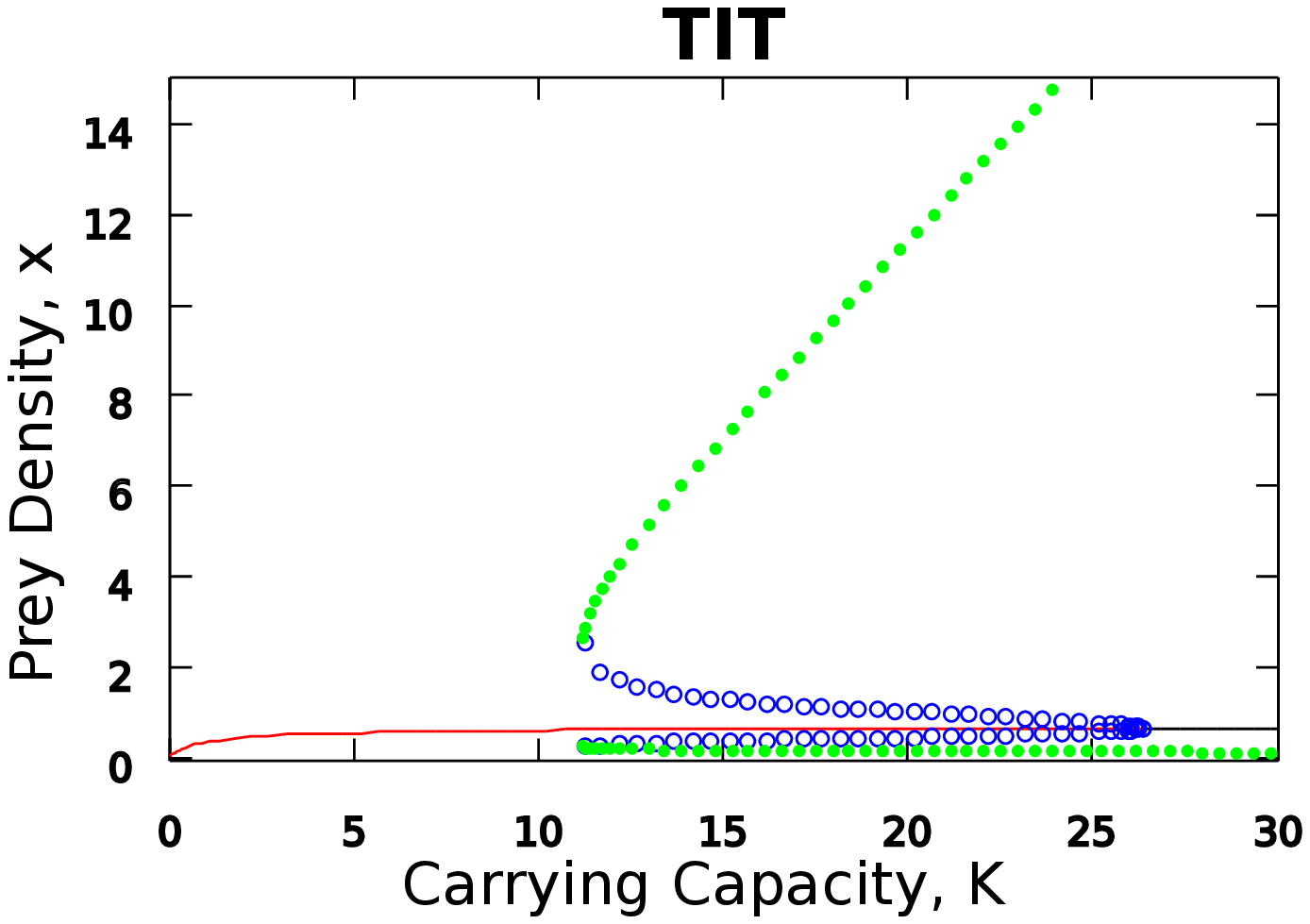}
    \includegraphics[width=0.325\linewidth]{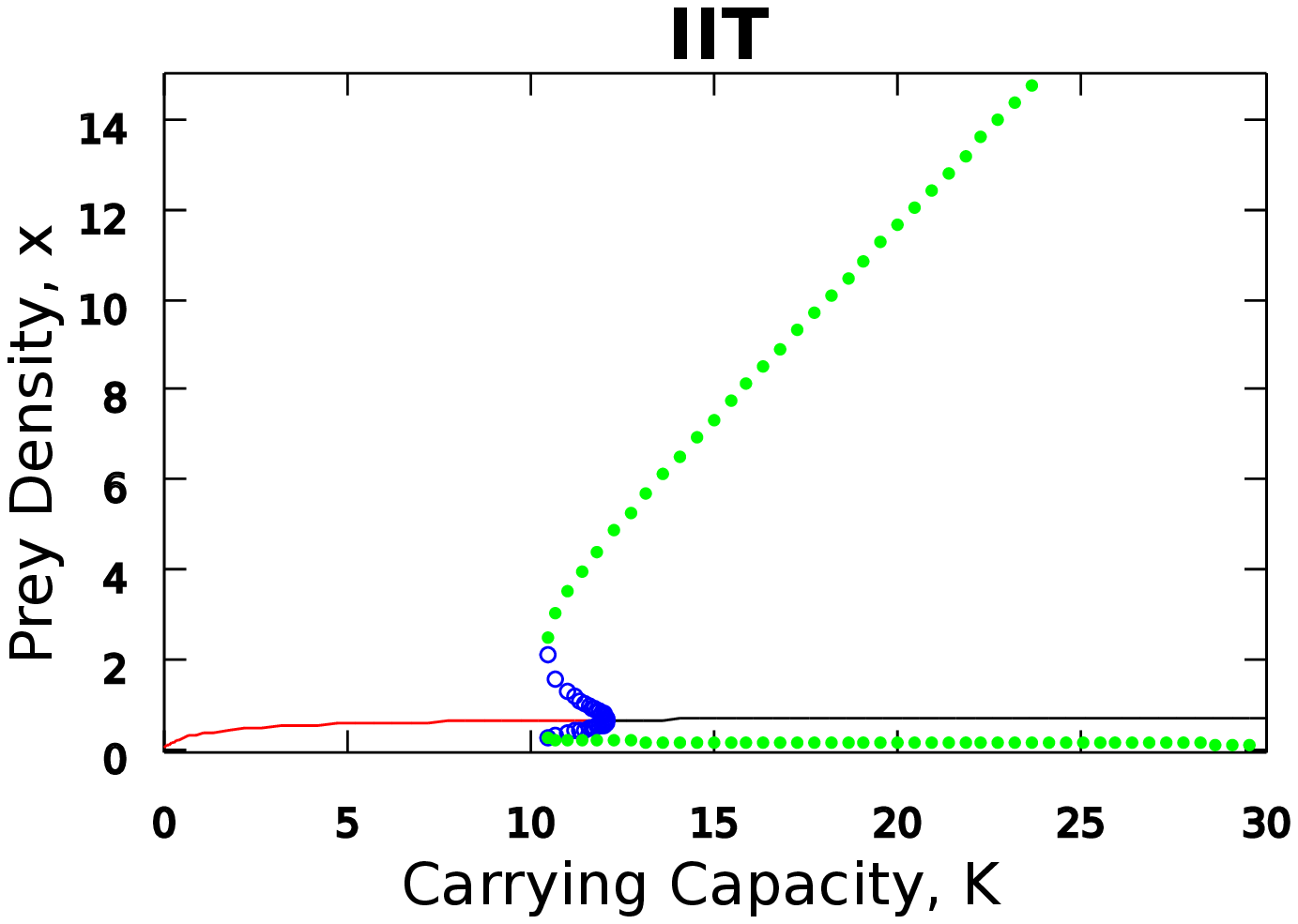}
    \\
    
    \includegraphics[width=0.325\linewidth]{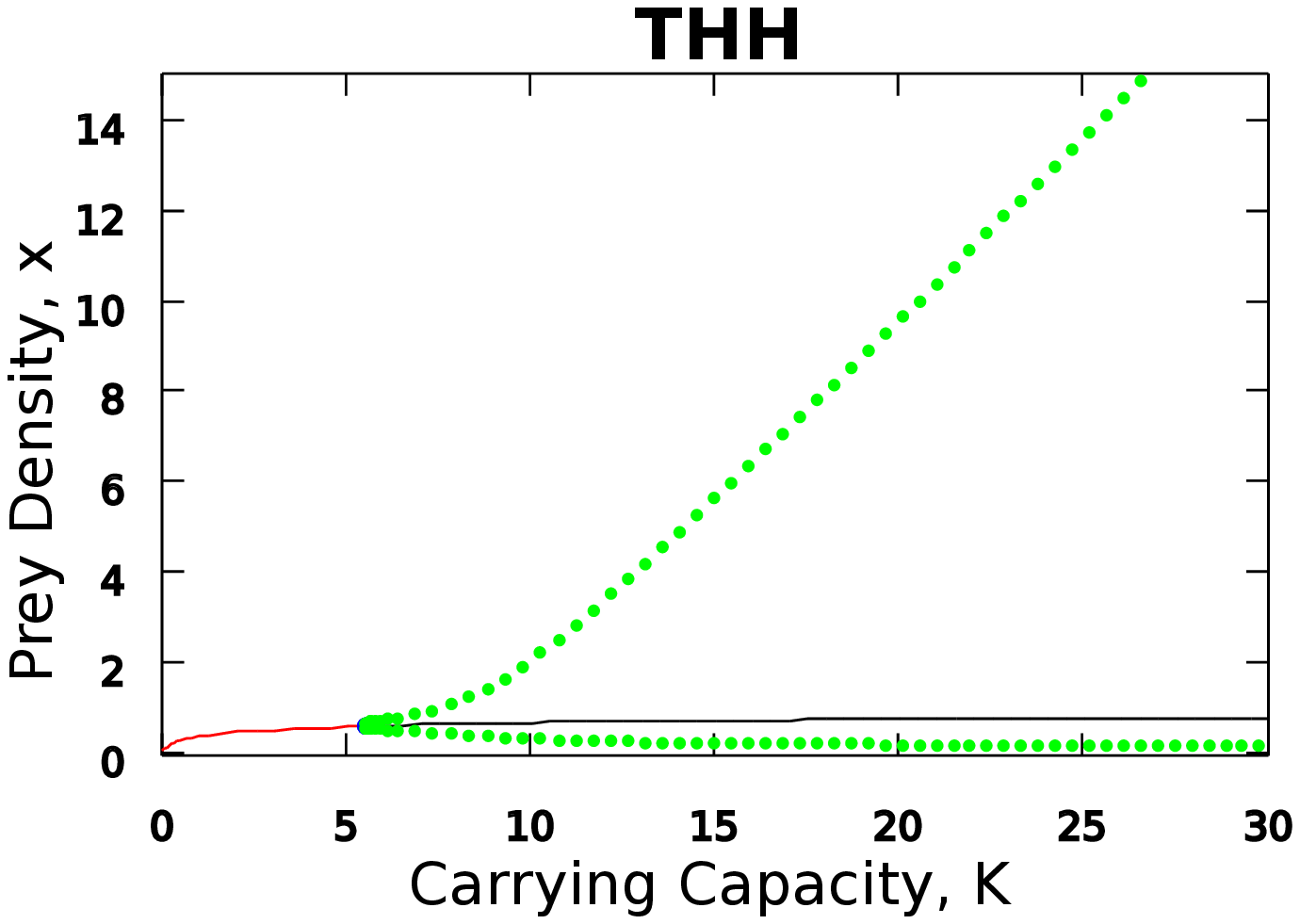}
    \includegraphics[width=0.325\linewidth]{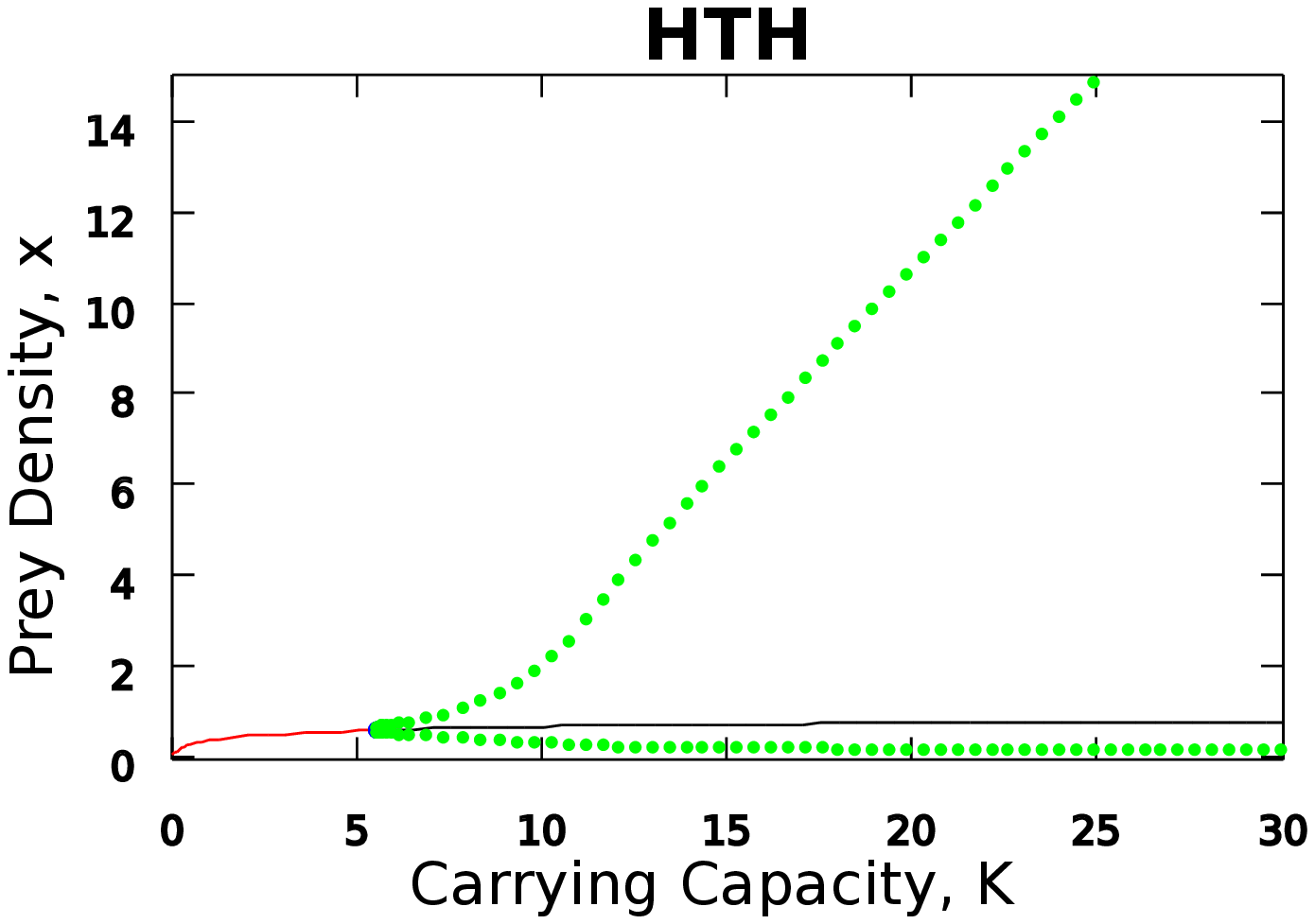}
    \includegraphics[width=0.325\linewidth]{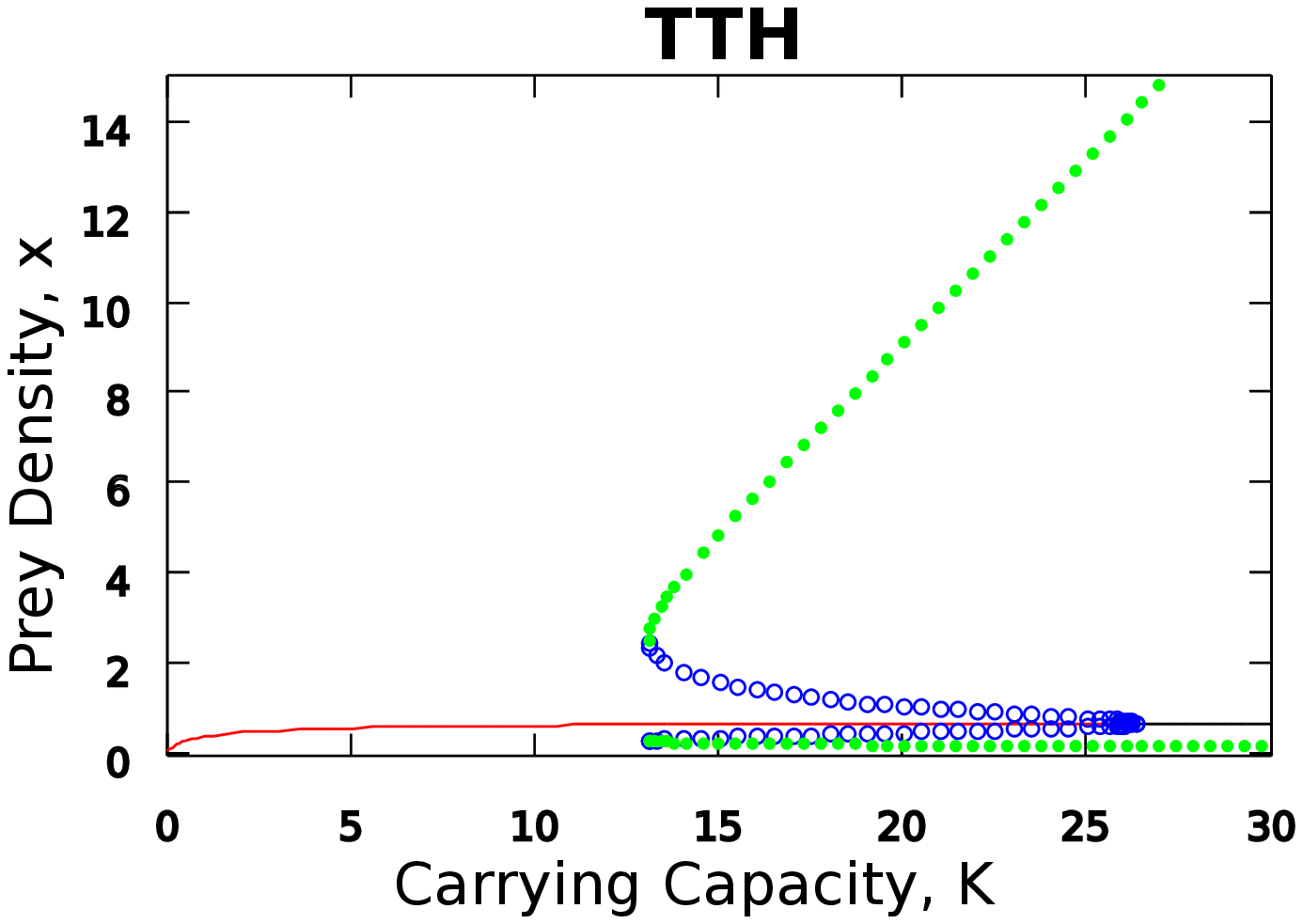}
    \\
    \includegraphics[width=0.325\linewidth]{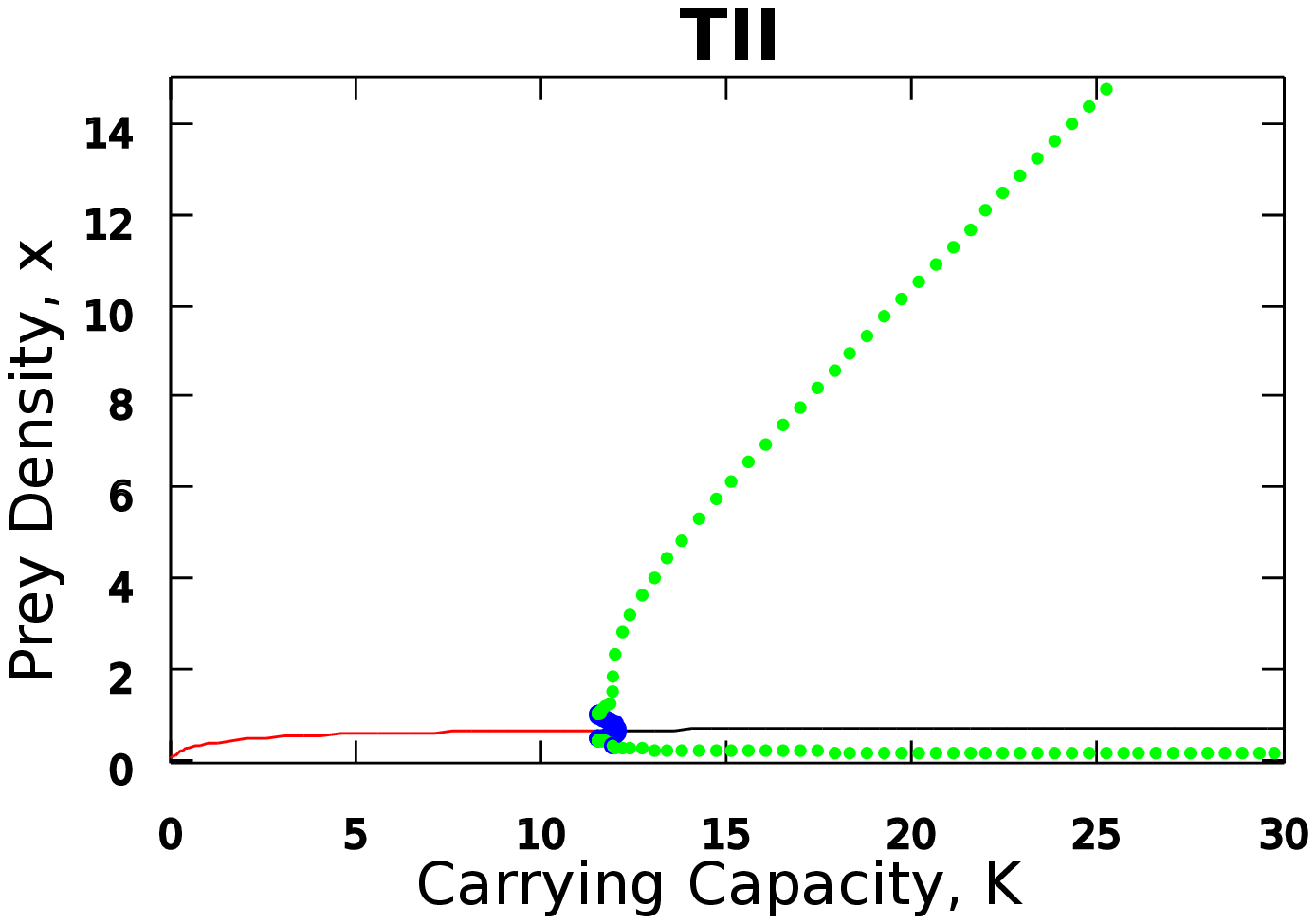}
    \includegraphics[width=0.325\linewidth]{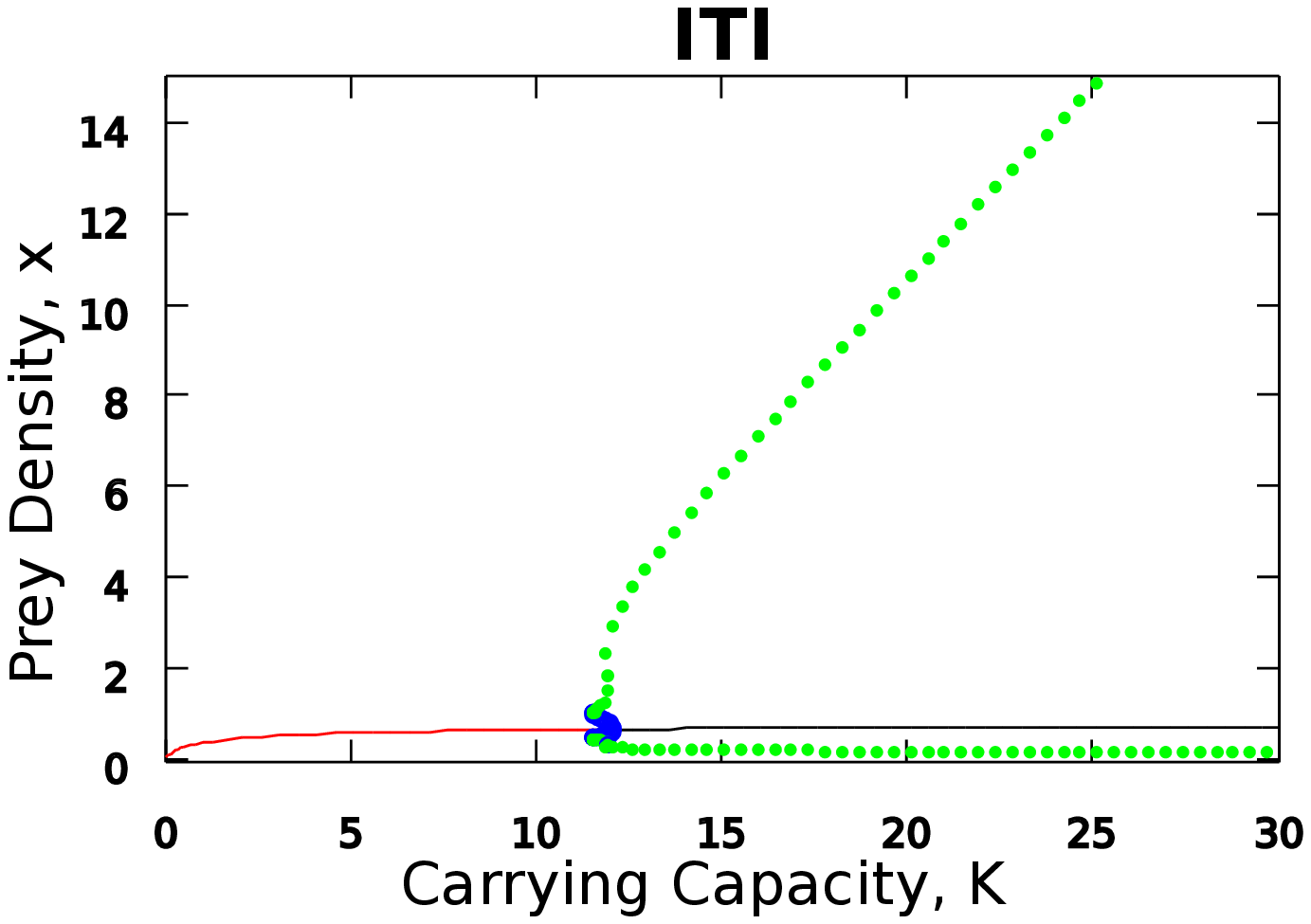}
    \includegraphics[width=0.325\linewidth]{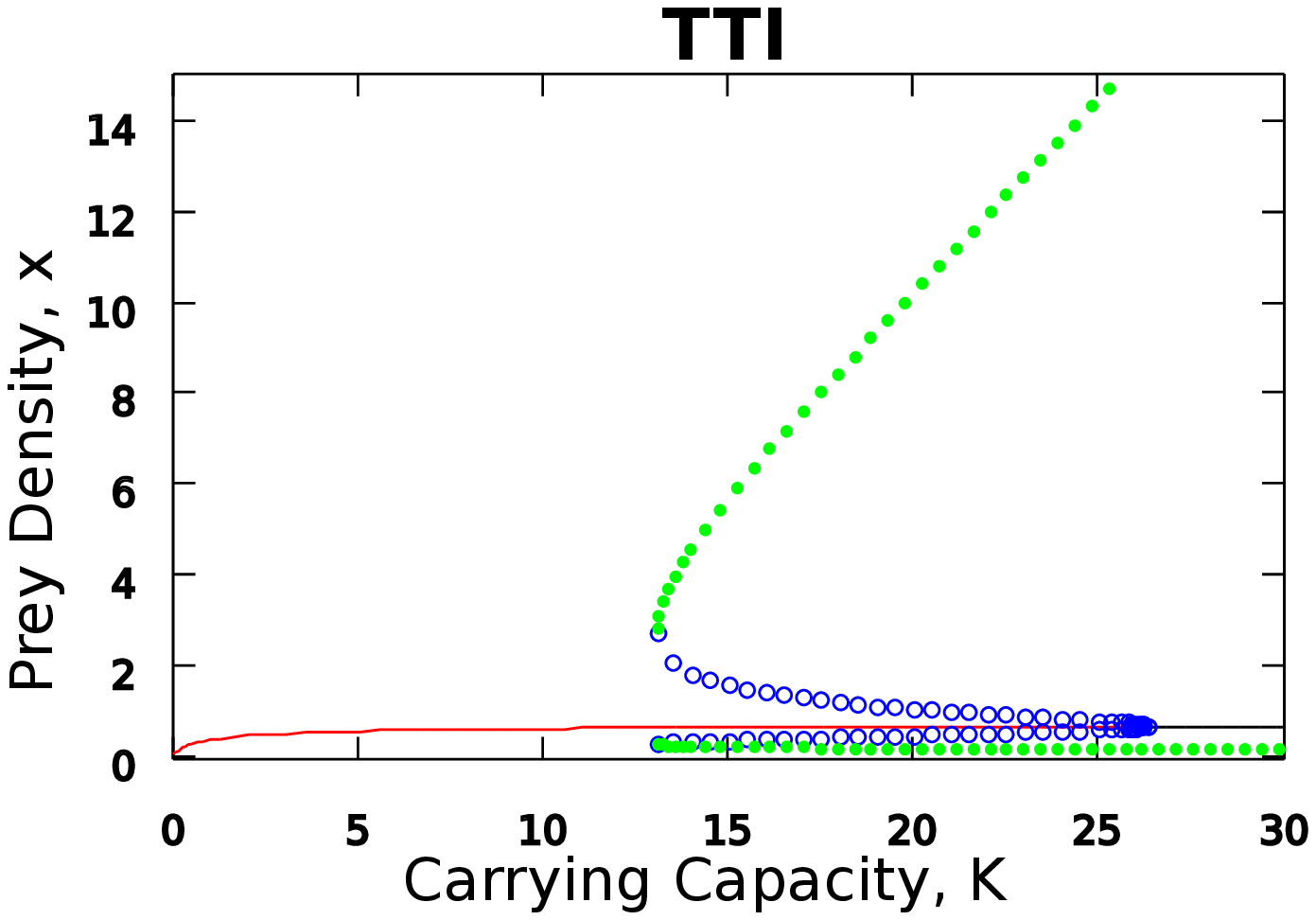}    
    \\
    
    \includegraphics[width=\linewidth]{Pictures/FIGLeg2.eps}
    \caption{Bifurcation diagrams for the Leslie-Gower-May model created using functional responses defined piecewise from two functional responses and switching at the corresponding intersections. Values for the intersections are given in the manuscript and bifurcation values are given in Table \ref{suptable:BifurValues}. }
    \label{supfig:SwitchingBifurMay}
\end{figure}

The bifurcation diagrams corresponding to the Rosenzweig-MacArthur and Leslie-Gower-May piecewise identical models are given in Figures \ref{supfig:SwitchingBifurRM} and \ref{supfig:SwitchingBifurMay}, respectively. As in the manuscript, none of the models show consistent decrease in structural sensitivity. The Holling and Ivlev cases for the Rosenzweig-MacArthur model (first three rows of Figure \ref{supfig:SwitchingBifurRM}) show behaviour corresponding to the functional responses used to describe low and medium prey densities. On the other hand, the trigonometric Rosenzweig-MacArthur model (last row of Figure \ref{supfig:SwitchingBifurRM}) and the Leslie-Gower-May models show model behaviour consistent with the functional response that describes more regions of prey density.

\begin{table}[ht]
\centering
\begin{tabular}{|c|c|c|c|c|c|c|}
    \hline 
   & \multicolumn{3}{|c|}{Rosenzweig-MacArthur} & \multicolumn{3}{|c|}{Leslie-Gower-May} \\
    \hline
    
   \textbf{\large Original Model (M)} & \textbf{Holling} & \textbf{Ivlev} & \textbf{trig} & \textbf{Holling} & \textbf{Ivlev} & \textbf{trig} \\
    \hline
    Hopf & 0.4452 & 1.071 & 10.12 & 5.554 & 12.09 & 26.36 \\
    Limit Cycle Saddle & - & - & 2.644 & - & 10.51 & 12.94\\
    \hline
    
    \textbf{\large Holling Fixed (M)} & \textbf{Holling} & \textbf{Ivlev} & \textbf{trig} & \textbf{Holling} & \textbf{Ivlev} & \textbf{trig} \\
    \hline
    Hopf & 0.4452 & 0.4632 & 0.7729 & 5.554 & 5.857 & 6.052 \\
    Limit Cycle Saddle & - & - & 0.5057 & - & - & 5.964\\
    \hline
    
    \textbf{\large Ivlev Fixed (\ref{supfig:CoexSSBifur})} & \textbf{Holling} & \textbf{Ivlev} & \textbf{trig} & \textbf{Holling} & \textbf{Ivlev} & \textbf{trig} \\
    \hline
    Hopf & 1.063 & 1.071 & 1.343 & 12.60 & 12.09 & 11.83 \\
    Limit Cycle Saddle & - & - & 0.8089 & - & 10.51 & 9.625\\
    \hline
    
    \textbf{\large Trigonometric Fixed (\ref{supfig:CoexSSBifur})} & \textbf{Holling} & \textbf{Ivlev} & \textbf{trig} & \textbf{Holling} & \textbf{Ivlev} & \textbf{trig} \\
    \hline
    Hopf & 10.36 & 10.36 & 10.12 & 31.69 & 27.69 & 26.36 \\
    Limit Cycle Saddle & - & - & 2.644 & 23.45 & 14.91 & 12.94\\
    \hline

    \textbf{\large Pw-I Holling (M, \ref{supfig:SwitchingBifurRM},  \ref{supfig:SwitchingBifurMay})} & \textbf{HII} & \textbf{IHI} & \textbf{HHI} & \textbf{HII} & \textbf{IHI} & \textbf{HHI} \\
    \hline
    Hopf & 0.4452 & 0.4452 & 0.4452 & 12.09 & 12.09 & 5.554 \\
    Limit Cycle Saddle & - & - & -  & 7.518 & 7.518 & - \\
    \hline
    & \textbf{HTT} & \textbf{THT} & \textbf{HHT} & \textbf{HTT} & \textbf{THT} & \textbf{HHT} \\
    \hline
    Hopf & 0.4451 & 0.4451 & 0.4451 & 26.36 & 26.36 & 5.554 \\
    Limit Cycle Saddle & - & - & - & 7.95 & 7.95 & - \\
    \hline
    
    \textbf{\large Pw-I Ivlev (\ref{supfig:SwitchingBifurRM}, \ref{supfig:SwitchingBifurMay})} & \textbf{IHH} & \textbf{HIH} & \textbf{IIH} & \textbf{IHH} & \textbf{HIH} & \textbf{IIH} \\
    \hline
    Hopf & 1.071 & 1.071 & 1.071 & 5.554 & 5.554 & 12.09 \\
    Limit Cycle Saddle & - & - & - & - & - & 10.51 \\
    \hline
    & \textbf{ITT} & \textbf{TIT} & \textbf{IIT} & \textbf{ITT} & \textbf{TIT} & \textbf{IIT} \\
    \hline
    Hopf & 1.071 & 1.071 & 1.071 & 26.36 & 26.36 & 12.09 \\
    Limit Cycle Saddle & - & - & - & 11.10 & 11.24 & 10.47 \\
    \hline
    
    \textbf{\large Pw-I Trigonometric (M, \ref{supfig:SwitchingBifurRM}, \ref{supfig:SwitchingBifurMay})} & \textbf{THH} & \textbf{HTH} & \textbf{TTH} & \textbf{THH} & \textbf{HTH} & \textbf{TTH} \\
    \hline
    Hopf & 10.12 & 10.12 & 10.12 & 5.554 & 5.554 & 26.36 \\
    Limit Cycle Saddle & 2.431 & 2.431 & 2.644 & - & - & 13.13 \\
    \hline
    & \textbf{TII} & \textbf{ITI} & \textbf{TTI} & \textbf{TII} & \textbf{ITI} & \textbf{TTI} \\
    \hline
     Hopf & 10.12 & 10.12 & 10.12 & 12.09 & 12.09 & 26.36\\
    Limit Cycle Saddle & 2.551 & 2.551 & 2.644 & 11.54 & 11.54 & 13.14\\
    \hline
\end{tabular}
    \caption{Values of the bifurcation parameter, the carrying capacity, for each bifurcation shown in the respective figures. The values in parentheses give the figure number of the corresponding bifurcation diagrams and M corresponds to figures given in the manuscript. Pw-I stands for piecewise identical, H stands for Holling, I stands for Ivlev, and T stands for trigonometric.}
    \label{suptable:BifurValues}
\end{table}

%\bibliographystyle{spmpscinat.bst}
%\bibliography{Sources.bib}